\documentclass[12pt]{amsart}
\usepackage[top=80pt,bottom=80pt,left=70pt,right=70pt]{geometry}
\usepackage{graphicx} % Required for inserting images
\usepackage{amsmath}
\usepackage{enumitem}
\usepackage{amssymb}
\usepackage{amsthm}
\usepackage{hyperref}
\usepackage{tikz}
\usepackage{tikz-cd} 

\usepackage[utf8]{inputenc}

\title[Higher rank partially hyperbolic actions on nilmanifolds]{Global rigidity for some partially hyperbolic abelian actions with 1-dimensional center}
\author{Sven Sandfeldt}
\email{svensan@kth.se}

\date{December 2023}

\newcommand{\norm}[1]{\left\lVert#1\right\rVert}

\newcommand{\intd}{{\rm d}}

\makeatletter
\newtheorem*{rep@theorem}{\rep@title}
\newcommand{\newreptheorem}[2]{%
\newenvironment{rep#1}[1]{%
 \def\rep@title{#2 \ref{##1}}%
 \begin{rep@theorem}}%
 {\end{rep@theorem}}}
\makeatother

\counterwithin{equation}{section}

\newtheorem{mainTheorem}{Theorem}

\newreptheorem{mainTheorem}{Theorem}

\newreptheorem{mainCorollary}{Corollary}

\newtheorem{theorem}{Theorem}[section]

\newtheorem{lemma}{Lemma}[section]

\newtheorem{definition}{Definition}[section]
\newtheorem*{definition*}{Definition}

\theoremstyle{definition}

\newtheorem{remark}{Remark}

\begin{document}

\begin{abstract}
We obtain a global rigidity result for abelian partially hyperbolic higher rank actions on certain $2-$step nilmanifolds $X_{\Gamma}$. We show that, under certain natural assumptions, all such actions are $C^{\infty}-$conjugated to an affine model. As a consequence, we obtain a centralizer rigidity result, classifying all possible centralizers for any $C^{1}-$small perturbation of an irreducible, affine partially hyperbolic map on $X_{\Gamma}$. Along the way, we prove two results of independent interest. We describe fibered partially hyperbolic diffeomorphisms on $X_{\Gamma}$ and we show that topological conjugacies between partially hyperbolic actions and higher rank affine actions are $C^{\infty}$.
\end{abstract}

\maketitle

%\tableofcontents

%\input{RunLargeCenterRigidity}

\section{Introduction}

Rigidity of $\mathbb{Z}^{k}-$actions on tori (and nilmanifolds) with some hyperbolicity have been studied extensively. The general philosophy is: large abelian actions with some hyperbolicity should be \textbf{globally rigid}, i.e. smoothly conjugated to algebraic models. A big breakthrough result in this direction was obtained by Katok and Spatzier \cite{KatokSpatzier1994} where they prove that all perturbations of certain algebraic Anosov actions are smoothly conjugated back to the corresponding algebraic models. There they also outline a rigidity program for abelian actions with hyperbolicity. Since the paper by Katok and Spatzier, a lot of results have been obtained for large abelian actions close to some algebraic model with some hyperbolicity \cite{DamjanovicKatok2010,DamjanovicKatok2011,DamjanovicFayad2019,Vinhage2015,VinhageWang2019,Wang2010:1,Wang2010:2,Wang2022,EinsiedlerFisher2007,DamjanovicQinbo2022}. Even earlier than the result by Katok and Spatzier, Katok and Lewis \cite{KatokLewis1991} proved a global rigidity statement for Anosov action on tori. In \cite{KatokLewis1991}, the authors show that a certain class of Anosov $\mathbb{Z}^{d}-$actions on the torus $\mathbb{T}^{d+1}$ is, necessarily, smoothly conjugated to an algebraic model, even though the action might not be close to the algebraic model. A crucial assumption in \cite{KatokLewis1991} is that the $\mathbb{Z}^{d}-$action contains \textbf{many Anosov elements}. The property of having many Anosov elements was removed by F. Rodriguez Hertz in \cite{Rodriguez-Hertz2007}, where Rodriguez Hertz only assumes that the action contains one Anosov element. Nilmanifolds are natural generalizations of tori. There has been a lot of work studying the global rigidity of higher rank\footnote{See Definition \ref{Def:HigherRankActions}.} Anosov actions on nilmanifolds since the paper by Katok and Lewis \cite{KalininSadovskaya2006,KalininSadovskaya2007,KalininSpatzier2007,FisherKalininSpatzier2011}. The culmination of these works is the result by F. Rodriguez Hertz, Z. Wang and D. Fisher, B. Kalinin, R. Spatzier \cite{Rodriguez-HertzWang2014,FisherKalininSpatzier2013} proving that abelian higher rank Anosov actions on (infra-)nilmanifolds are smoothly conjugated to algebraic models, completely resolving the question of global rigidity of abelian higher rank Anosov actions on these manifolds. Relaxing the Anosov assumption leads to the question: 
\begin{center}
\textit{When are partially hyperbolic $\mathbb{Z}^{k}-$actions on nilmanifolds by affine maps?}
\end{center}
These actions have been remarkably resistant to classification results. Even local rigidity for algebraic partially hyperbolic actions on (non-toral) nilmanifolds has been open for a long time until recent advances by Z. J. Wang \cite{Wang2022}. In the global Anosov setting, any higher rank action is topologically conjugated to an affine action by the global topological rigidity for Anosov diffeomorphisms \cite{Franks1969,Manning1974}. Because of this, the key problem to solve is to upgrade the regularity of an already existing topological conjugacy. In contrast, for partially hyperbolic diffeomorphisms there is no global topological rigidity. In fact, there are nilmanifolds on which a generic partially hyperbolic diffeomorphism can not be conjugated to an affine model (for example, the manifolds considered in Theorem \ref{Thm:MainTheorem1}). The reason for this is that affine partially hyperbolic diffeomorphisms always have \textbf{isometric center} whereas this is a very special property for general partially hyperbolic systems. As a consequence, when we study the global rigidity of partially hyperbolic actions, a significant problem to solve is the existence of a topological conjugacy. However, once the existence of a topological conjugacy is proved it turns out that upgrading to a smooth conjugacy follows similarly to the Anosov setting (see Theorem \ref{Thm:ImprovedRegularity1} and Section \ref{Sec:ImprovedRegularity}). In this paper, we produce an initial result towards answering the global rigidity question for partially hyperbolic actions.  

In Theorem \ref{Thm:MainTheorem1} we prove the first global rigidity result for higher rank abelian actions on nilmanifolds with one partially hyperbolic element. The nilmanifolds under consideration are products between tori and Heisenberg nilmanifolds. These manifolds are, in a sense, the simplest for studying partially hyperbolic diffeomorphisms (with $1-$dimensional center). The reason for this is that the affine models of partially hyperbolic systems on these manifolds simultaneously satisfy three important properties, (i) they have $1-$dimensional center direction, (ii) the center leaves are compact and (iii) these systems are \textbf{accessible} (see Section \ref{Sec:BackgroundPartiallyHyperbolic}). This allows us to use significant portions of the partially hyperbolic theory when studying these systems. In contrast, when studying Anosov actions the simplest manifolds are tori. However, there are no accessible affine partially hyperbolic systems on tori, so studying global rigidity for partially hyperbolic actions on tori (that is not derived from Anosov) becomes very difficult because many of the tools from partially hyperbolic theory are not applicable.

The questions of local and global rigidity of higher rank actions can also be studied for different types of actions, either dropping the assumption that the action is abelian, or dropping the assumption that the action should have some hyperbolicity. Local rigidity has been obtained for large abelian parabolic actions, with no hyperbolicity \cite{DamjanovicKatok2011Parabolic,DamjanovicFayadSaprykina2023,Wang2019,DamjanovicQinboPetkovic2022}. Removing the assumption that the action should be abelian, we can study the actions of higher rank lattices in semi-simple Lie groups, see for example \cite{Fisher2019} and the references therein. In fact, the rigidity result for Abelian actions in \cite{KatokLewis1991} was used to obtain local rigidity of ${\rm SL}(n,\mathbb{Z})-$action on tori. When considering rigidity of Anosov lattice actions on nilmanifolds, there are also global results, see for example \cite{BrownRodriguez-HertzWang2017}. A key point in \cite{BrownRodriguez-HertzWang2017} is that any conjugacy between the hyperbolic lattice action and the algebraic model also conjugates the action of a large abelian subgroup to some algebraic abelian action. So, the results of \cite{Rodriguez-HertzWang2014,FisherKalininSpatzier2013} can be applied to improve the regularity of the conjugacy. Considering the main results of this paper, see Theorem \ref{Thm:MainTheorem1}, a natural question is: are higher rank partially hyperbolic lattice actions on nilmanifolds always by affine maps? Conjecturally this question should have an affirmative answer (see \cite[Conjecture 5]{Gorodnik2007}). In \cite{LeeSandfeldt2024} H. Lee and the author show, by using results from this paper, that the question has an affirmative answer when the manifold is a Heisenberg nilmanifold.

\subsection{Global rigidity of partially hyperbolic actions}

Let $G$ be a simply connected $\ell-$step nilpotent Lie group. That is, the lower central series $G^{(1)} = G$, $G^{(j+1)} = [G,G^{(j)}]$, terminate at $\ell+1$, $G^{(\ell+1)} = e$. Given a lattice $\Gamma\leq G$ we define the associated \textbf{compact nilmanifold} as the quotient $X_{\Gamma} = \Gamma\setminus G$. Compact nilmanifolds have associated groups of \textbf{automorphisms} and \textbf{affine maps}
\begin{align}
& {\rm Aut}(X_{\Gamma}) = \{L\in{\rm Aut}(G)\text{ : }L\Gamma = \Gamma\}, \\
& {\rm Aff}(X_{\Gamma}) = \{f_{0}(x) = L(x)g^{-1}\text{ : }L\in{\rm Aut}(X_{\Gamma}),\text{ }g\in G\}.
\end{align}
By automorphism rigidity of nilpotent lattices \cite{DiscreteSubgroupsOfLieGroups1972} we can, equivalently, define ${\rm Aut}(X_{\Gamma}) = {\rm Aut}(\Gamma) = {\rm Aut}(\pi_{1}X_{\Gamma})$. An automorphism $\rho:\mathbb{Z}^{k}\to{\rm Aut}(X_{\Gamma})$ is said to have a rank$-1$ factor if there is some quotient $\Hat{X}_{\Gamma}$ of $X_{\Gamma}$ such that $\rho$ descends to $\Hat{X}_{\Gamma}$ and the induced map $\Hat{\rho}:\mathbb{Z}^{k}\to{\rm Aut}(\Hat{X}_{\Gamma})$ factor through a map $\mathbb{Z}\to{\rm Aut}(\Hat{X}_{\Gamma})$. A homomorphism $\rho:\mathbb{Z}^{k}\to{\rm Aut}(X_{\Gamma})$ is \textbf{higher rank} if it has no rank$-1$ factor. More generally, given a smooth action $\alpha:\mathbb{Z}^{k}\times X_{\Gamma}\to X_{\Gamma}$ we have an induced map $\alpha_{*}:\mathbb{Z}^{k}\to{\rm Aut}(\pi_{1}X_{\Gamma})$. We say that $\alpha$ is higher rank if the induced map $\alpha_{*}$ is higher rank.

A diffeomorphism $f:X_{\Gamma}\to X_{\Gamma}$ is \textbf{partially hyperbolic} if there is a $Df-$invariant splitting $TX_{\Gamma} = E^{s}\oplus E^{c}\oplus E^{u}$ such that $E^{s}$ is exponentially contracted, $E^{u}$ is exponentially expanded and the behaviour along $E^{c}$ is dominated by the behaviour of $Df$ along $E^{s}$ and $E^{u}$ (for a precise definition, see Section \ref{Sec:BackgroundPartiallyHyperbolic}). We will need a technical assumption on the \textbf{stable and unstable foliations} $W^{s}$ and $W^{u}$, tangent to $E^{s}$ and $E^{u}$. We say that $W^{s}$ and $W^{u}$ are \textbf{quasi-isometric in the universal cover} if the metric along the leaves of the foliations, after we lift them to the universal cover, is comparable to the ambient metric (see Definition \ref{Def:QuasiIsometricLeaves}). If the \textbf{center bundle} $E^{c}$ is the trivial bundle then $f$ is Anosov. The main result of this paper is an extension of the results of \cite{Rodriguez-HertzWang2014} to certain nilmanifolds by weakening the assumption that $\alpha$ is Anosov. Instead, we assume that the action $\alpha$ contains a partially hyperbolic element.
\begin{mainTheorem}\label{Thm:MainTheorem1}
Let $G$ be a $2-$step nilpotent Lie group with $\dim[G,G] = 1$, $\Gamma\leq G$ a lattice and $X_{\Gamma} = \Gamma\setminus G$ the associated nilmanifold. Let $\alpha:\mathbb{Z}^{k}\times X_{\Gamma}\to X_{\Gamma}$ be a smooth higher rank action with $\mathbf{n}_{0}\in\mathbb{Z}^{k}$ such that $f = \alpha^{\mathbf{n}_{0}}$ is partially hyperbolic and satisfying
\begin{enumerate}[label = (\roman*)]
    \item $f$ has $1-$dimensional center,
    \item the stable and unstable foliations $W^{s}$, $W^{u}$ are quasi-isometric in the universal cover.
\end{enumerate}
Then $\alpha$ is $C^{\infty}-$conjugated to some affine action $\alpha_{0}:\mathbb{Z}^{k}\to{\rm Aff}(X_{\Gamma})$.
\end{mainTheorem}
\begin{remark}
Condition $(ii)$ in the theorem is a technical assumption that guarantees that $f$ is \textbf{fibered} (see Theorem \ref{Thm:PropertiesOfPHdiffeos}). Motivated by results in \cite{Hammerlindl2013:1,Hammerlindl2013:2,HammerlindlPotrie2015}, it seems plausible that $(ii)$ is always satisfied. Moreover, in \cite[Lemma A2]{LeeSandfeldt2024} it is shown that if $f$ is topologically conjugated (or more generally, leaf conjugated) to some affine $f_{0}$ then condition $(ii)$ is satisfied. It follows that $\alpha$ satisfies the assumptions of the theorem if $\alpha$ satisfies the conclusion of the theorem, so the theorem is sharp in this sense. It also follows that if there are exotic partially hyperbolic systems that do not satisfy $(ii)$, then there exists no global rigidity result for these systems.
\end{remark}
\begin{remark}
With $G$ as in Theorem \ref{Thm:MainTheorem1} the group $G$ can be written as $G = H^{n}\times\mathbb{R}^{m}$, $n\neq0$, where $H^{n}$ is a Heisenberg group. The nilmanifold $X_{\Gamma}$ is also a product of a Heisenberg nilmanifold and a torus. On these manifolds there exists no Anosov actions since the derived subgroup, $[G,G]$, is isometric for any automorphism. So Theorem \ref{Thm:MainTheorem1} is the only global rigidity result on these manifolds, since \cite{Rodriguez-HertzWang2014,FisherKalininSpatzier2013} do not apply. In fact, to the author's knowledge, Theorem \ref{Thm:MainTheorem1} is the first global rigidity result for abelian actions assuming only one partially hyperbolic element.
\end{remark}
\begin{remark}
In principle the proof of Theorem \ref{Thm:MainTheorem1} should work for $\ell-$step $G$ with $\ell > 2$ as long as $\dim G^{(\ell)} = 1$. In this case the quasi-isometric assumption, assumption $(ii)$, would have to be changed. This is a work in progress.
\end{remark}
\begin{remark}
If $\alpha$ is topologically conjugated to some affine action, then the methods from \cite{Rodriguez-HertzWang2014} generalize to partially hyperbolic systems, see Theorem \ref{Thm:ImprovedRegularity1}. So, the main novelty of Theorem \ref{Thm:MainTheorem1} is that we produce a topological conjugacy from $\alpha$ to an affine action $\alpha_{0}$.
\end{remark}

\subsection{Applications to centralizer classification and centralizer rigidity}

Given a diffeomorphism $f:M\to M$ on a closed manifold we define its \textbf{smooth centralizer} as the group of diffeomorphisms that commute with $f$. That is, we define
\begin{align}
Z^{\infty}(f) = \{g\in{\rm Diff}^{\infty}(M)\text{ : }fg = gf\}.
\end{align}
We are interested in two questions about the group $Z^{\infty}(f)$:
\begin{enumerate}[label = (\roman*)]
    \item What are possible groups that arise as $Z^{\infty}(f)$ for some $f\in{\rm Diff}^{\infty}(M)$?
    \item If $Z^{\infty}(f)$ is large (compared to the conjecturally generic size $\mathbb{Z}$, \cite{Smale1991,Smale1998}) what can be said about $f$?
\end{enumerate}
In this level of generality, questions $(i)$ and $(ii)$ are difficult (or possibly impossible) to answer. Instead, we fix $f_{0}\in{\rm Aff}(X_{\Gamma})$ for some $X_{\Gamma} = \Gamma\setminus G$, and consider questions $(i)$ and $(ii)$ for those $f\in{\rm Diff}^{\infty}(X_{\Gamma})$ that are $C^{1}-$close to $f_{0}$. We call $(i)$ the question of \textbf{local centralizer classification} and $(ii)$ the question of \textbf{local centralizer rigidity} around $f_{0}$. These questions were raised and addressed by Damjanović, Wilkinson and Xu in \cite{DamjanovicWilkinsonXu2021} where the authors study perturbations of time$-t_{0}$ map of geodesic flows on negatively curved manifolds and trivial circle extensions of hyperbolic automorphisms. In \cite{BarthelemeGogolev2021} the authors study local centralizer rigidity of time$-1$ maps of Anosov flows on $3-$manifolds, generalizing results from \cite{DamjanovicWilkinsonXu2021} in the context of $3-$manifolds. Another generalization of results from \cite{DamjanovicWilkinsonXu2021} was obtained by W. Wang in \cite{Wang2023}, where semi-simple Lie groups of higher rank were studied instead of rank$-1$ simple groups. For ergodic toral automorphisms, Gan, Xu, Shi and Zhang studied partially hyperbolic diffeomorphisms on $\mathbb{T}^{3}$ homotopic to an hyperbolic automorphism \cite{GanShiXuZhang2022}. In \cite{Sandfeldt2023} the author studies local centralizer classification and rigidity for some partially hyperbolic, irreducible\footnote{An automorphism $A\in{\rm GL}(d,\mathbb{Z})$ is irreducible if the characteristic polynomial $p_{A}(t)$ is irreducible in $\mathbb{Q}[t]$.} toral automorphisms.

If $f$ is partially hyperbolic with (uniquely integrable) center foliation $W^{c}$, then we obtain a normal subgroup $Z_{c}^{\infty}(f)\subset Z^{\infty}(f)$, the \textbf{center fixing centralizer}:
\begin{align}\label{Eq:CenterFixingCentralizer}
Z_{c}^{\infty}(f) := \{g\in Z^{\infty}(f)\text{ : }gx\in W^{c}(x),\text{ }x\in X_{\Gamma}\}.
\end{align}
From \cite[Theorem 5]{DamjanovicWilkinsonXu2023}, if $Z_{c}^{\infty}(f)$ is sufficiently big and $f$ is fibered (see \cite[Definition 1]{DamjanovicWilkinsonXu2023}) then $f$ is smoothly conjugated to an isometric extension of an Anosov map. Combining this with Theorem \ref{Thm:MainTheorem1} we completely classify the centralizers of diffeomorphisms $C^{1}-$close to affine partially hyperbolic maps.

Let $G$ be the $(d+1)-$dimensional Heisenberg group and $X_{\Gamma} = \Gamma\setminus G$ a compact Heisenberg nilmanifold. We have a natural fibration
\begin{align}
\pi:X_{\Gamma}\to\mathbb{T}^{d},
\end{align}
where we refer to $\mathbb{T}^{d}$ as the \textbf{base} of $X_{\Gamma}$. Any affine map on $X_{\Gamma}$ descends to an affine map on $\mathbb{T}^{d}$ and from the group relations in $G$ any automorphism $L\in{\rm Aut}(X_{\Gamma})$ induces an element of ${\rm Sp}(d,\mathbb{Z})$ on $\mathbb{T}^{d}$. Conversely, any $L\in{\rm Sp}(d,\mathbb{Z})$ defines an element of ${\rm Aut}(X_{\Gamma})$, see Section \ref{Sec:BackgroundNilmanifolds}. Given $f_{0}\in{\rm Aff}(X_{\Gamma})$ we denote by $L_{su}\in{\rm Sp}(d,\mathbb{Z})$ the induced automorphism on $\mathbb{T}^{d}$. Before stating the theorem we define for any $f_{0}\in{\rm Aff}(X_{\Gamma})$, with irreducible, hyperbolic induced map on the base $L_{su}\in{\rm Sp}(d,\mathbb{Z})$, the natural number
\begin{align}
r_{0}(f_{0}) = {\rm rank}\left(\frac{Z^{\infty}(f_{0})}{Z_{c}^{\infty}(f_{0})}\right) = {\rm rank}\left(Z_{{\rm Sp}(d,\mathbb{Z})}(L_{su})\right).
\end{align}
Lemma \ref{L:RankOfCentralizerInSymplecticGroup} explicitly calculates the number $r_{0}(f_{0})$, and if $d\geq 6$, then $r_{0}(f_{0}) > 1$.
\begin{mainTheorem}\label{Thm:MainTheorem2}
Let $X_{\Gamma}$ be a compact Heisenberg nilmanifold and let $f_{0}\in{\rm Aff}(X_{\Gamma})$ be partially hyperbolic with $1-$dimensional center and $L_{su}$ irreducible. If $f\in{\rm Diff}^{\infty}(X_{\Gamma})$ is $C^{1}-$close to $f_{0}$ then one of the following holds
\begin{enumerate}[label = (\roman*)]
    \item either $Z^{\infty}(f)$ is virtually trivial,
    \item or $Z^{\infty}(f)$ is virtually $\mathbb{Z}\times\mathbb{T}$ in which case $f$ is an isometric extension of some Anosov diffeomorphism on $\mathbb{T}^{d}$,
    \item or $Z^{\infty}(f)$ is virtually $\mathbb{Z}^{r_{0}}\times\mathbb{T}$ and if $r_{0} > 1$ then $f$ is $C^{\infty}-$conjugate to some (possibly different) affine map $\Tilde{f}_{0}\in{\rm Aff}(X_{\Gamma})$.
\end{enumerate}
\end{mainTheorem}
\begin{remark}
All cases $(i)$, $(ii)$ and $(iii)$ occur, so Theorem \ref{Thm:MainTheorem2} completely classifies the centralizer of $f\in{\rm Diff}^{\infty}(X_{\Gamma})$ close to a partially hyperbolic $f_{0}\in{\rm Aff}(X_{\Gamma})$. Case $(i)$ holds generically \cite{BonattiCroivisierWilkinson2009}. Case $(ii)$ can be produced by fixing some irreducible, hyperbolic $L\in{\rm Sp}(d,\mathbb{Z})$ and defining $f$ on $G\cong\mathbb{R}^{d}\times\mathbb{R}$ by $f(x,t) = (Lx,t + \beta(x))$ where the second coordinate is identified with $[G,G]\cong\mathbb{R}$ and $\beta:\mathbb{T}^{d}\to\mathbb{R}$ is a cocycle over $L$ that is not cohomologous to a constant. The last case holds when $f$ is $C^{\infty}-$conjugate to some affine $\Tilde{f}_{0}$, so in particular when we take the trivial perturbation $f = f_{0}$.
\end{remark}
\begin{remark}
Similar results as Theorem \ref{Thm:MainTheorem2} have been obtained independently by Damjanović, Wilkinson and Xu using different methods with additional assumptions \cite{DamjanovicWilkinsonXu(PrivateCommunication)}.
\end{remark}

\subsection{Partially hyperbolic maps on nilmanifolds}

When proving Theorems \ref{Thm:MainTheorem1} and \ref{Thm:MainTheorem2}, we use a description of partially hyperbolic diffeomorphisms on the nilmanifolds considered in Theorem \ref{Thm:MainTheorem1}. The main property that we show is that, under the assumptions of Theorem \ref{Thm:MainTheorem1}, the system $f$ is \textbf{fibered} in the terminology of \cite{AvilaVianaWilkinson2022}:
\begin{theorem}\label{Thm:PropertiesOfPHdiffeos}
Let $G = H^{n}\times\mathbb{R}^{m}$ (where we allow $n = 0$), $\Gamma\leq G$ a lattice with associated nilmanifold $X_{\Gamma} = \Gamma\setminus G$. Let $f\in{\rm Diff}^{\infty}(X_{\Gamma})$ be partially hyperbolic and satisfy $(i)$, $(ii)$ from Theorem \ref{Thm:MainTheorem1}. If $G$ is abelian we assume, in addition, that the induced map $f_{*}:H_{1}(X_{\Gamma})\to H_{1}(X_{\Gamma})$ has at least one rational eigenvalue. The following holds
\begin{enumerate}[label = (\roman*)]
    \item $f$ is dynamically coherent with global product structure,
    \item all foliations $W^{\sigma}$, $\sigma = s,c,u,cs,cu$, are uniquely integrable,
    \item the center foliation $W^{c}$ have compact oriented circle leaves,
    \item $f$ is fibered over some hyperbolic $L_{su}\in{\rm GL}(d,\mathbb{Z})$ in the sense that there is some Hölder $\Phi:X_{\Gamma}\to\mathbb{T}^{d}$ such that $\Phi(fx) = L_{su}\Phi(x)$, $W^{c}(x) = \Phi^{-1}(\Phi(x))$ and $\Phi$ is homotopic to the projection $\pi:X_{\Gamma}\to\mathbb{T}^{d}$,
    \item there is a finite index subgroup $Z_{\rm fix}^{\infty}(f)\leq Z^{\infty}(f)$ such that if $g\in Z_{\rm fix}^{\infty}(f)$ we have $\Phi(gx) = B\Phi(x)$ where $B\in{\rm GL}(d,\mathbb{Z})$ is defined by $B\Phi_{*} = \Phi_{*}g_{*}$, $B$ is the induced map on homology if $G$ is non-abelian,
\end{enumerate}
moreover, if $G$ is not abelian then
\begin{enumerate}[label = (\roman*)]
    \setcounter{enumi}{5}
    \item $f$ is accessible.
\end{enumerate}
\end{theorem}
\begin{remark}
The assumption that $f_{*}$ has at least one rational eigenvalue is to remove derived-from-Anosov examples since these examples are not fibered.
\end{remark}
\begin{remark}
Properties $(i)$ and $(ii)$ follow from \cite{Brin2003}.
\end{remark}
\begin{remark}
This Theorem is similar to the classification of partially hyperbolic diffeomorphisms on $3-$dimensional manifolds by Hammerlindl and Hammerlindl-Potrie \cite{Hammerlindl2013:1,Hammerlindl2013:2,HammerlindlPotrie2015}. In fact, in dimension $3$, using \cite{BrinBuragoIvanov2009}, Theorem \ref{Thm:PropertiesOfPHdiffeos} essentially reduces to the main results of \cite{Hammerlindl2013:1,Hammerlindl2013:2} (in \cite{Hammerlindl2013:2} we must make the extra assumption that the linearization $L\in{\rm GL}(3,\mathbb{Z})$ has at least one rational eigenvalue). 
\end{remark}

\subsection{Improved regularity of topological conjugacies between higher rank actions}

The conjugacy in Theorems \ref{Thm:MainTheorem1} and \ref{Thm:MainTheorem2}, case $(iii)$, is produced in two steps. First, we construct a topological conjugacy and second we show that the topological conjugacy is $C^{\infty}$. The second step is the content of the following theorem, that may be of independent interest.
\begin{theorem}\label{Thm:ImprovedRegularity1}
Let $X_{\Gamma}$ be a nilmanifold and $\alpha_{0}:\mathbb{Z}^{k}\to{\rm Aff}(X_{\Gamma})$ a homomorphism. Suppose that $\alpha_{0}$ is higher rank. If $\alpha:\mathbb{Z}^{k}\to{\rm Diff}^{\infty}(X_{\Gamma})$ is bi-Hölder conjugate to $\alpha_{0}$ by $H:X_{\Gamma}\to X_{\Gamma}$ and there is some $\mathbf{n}_{0}\in\mathbb{Z}^{k}$ such that $f = \alpha^{\mathbf{n}_{0}}$ is partially hyperbolic and accessible with center $\dim(E^{c}) = \dim(E_{\alpha_{0}^{\mathbf{n}_{0}}}^{c})$, then $H$ is a $C^{\infty}-$diffeomorphism.
\end{theorem}
Theorem \ref{Thm:ImprovedRegularity1} is a generalization of the global rigidity result by F. Rodriguez Hertz and Z. Wang \cite{Rodriguez-HertzWang2014} to some higher rank partially hyperbolic actions. In fact, large parts of the results in \cite{Rodriguez-HertzWang2014} generalize immediately to partially hyperbolic actions. One exception is that Rodriguez Hertz and Wang use a characterization of Anosov diffeomorphisms due to Mañé \cite{Mane1977}, to show that many elements of the action are Anosov. This characterization can not be applied in the partially hyperbolic setting. We also change some technical aspects of the proof, removing the use of Pesin theory.

\subsection{Description of proofs}

Let $G$ and $\alpha:\mathbb{Z}^{k}\times X_{\Gamma}\to X_{\Gamma}$ be as in Theorem \ref{Thm:MainTheorem1} and $f = \alpha^{\mathbf{n}_{0}}$ the partially hyperbolic element. By considering the Lie algebra of $G$, $\mathfrak{g}$, it is immediate that $G$ takes the form $G = H^{n}\times\mathbb{R}^{m}$, $n\neq0$, where $H^{n}$ is the $(2n+1)-$dimensional Heisenberg group (for $n = 0$ we will consider $H^{n}$ as the trivial group $1$). Moreover, under the assumptions of Theorem \ref{Thm:MainTheorem2} $G$ has to be a Heisenberg group $H^{n}$ for some $n$ (this follows since any lattice $\Gamma$ in $H^{n}\times\mathbb{R}^{m}$ is, virtually, a product lattice so irreducibility of $(f_{0})_{*}$ guarantee that either $n = 0$ or $m = 0$). The proof of Theorem \ref{Thm:MainTheorem1} is divided into three steps, first we show that any $f$ as in Theorem \ref{Thm:MainTheorem1} is fibered, then we show that any action $\alpha$ as in \ref{Thm:MainTheorem1} is topologically conjugated to some affine model and finally we show that the topological conjugacy can be improved to a smooth conjugacy. 

\subsubsection{Step 1}

The first part of the proof of Theorem \ref{Thm:MainTheorem1} is to show that any element $f$ as in the theorem has to fiber over a hyperbolic automorphism of the torus and that $f$ is accessible. Both of these properties are contents of Theorem \ref{Thm:PropertiesOfPHdiffeos}. First, we obtain a map $\Phi$ on $G$, which is a contender for being the map in Theorem \ref{Thm:PropertiesOfPHdiffeos}, this is done as in \cite{Franks1969}. Second, we show that $\Phi$ is injective on the lifted stable and unstable leaves $\Hat{W}^{\sigma}$, $\sigma = s,u$. By invariance of domain, this implies that the stable and unstable distributions of $f$ satisfy $\dim(E^{\sigma})\leq\dim(E_{0}^{\sigma})$, where $E_{0}^{\sigma}$ is the corresponding distribution for the linearization $L$ of $f$. Since $f$ has center $E^{c}$ of dimension $1$ by assumption we conclude that $\dim(E^{\sigma}) = \dim(E_{0}^{\sigma})$ for $\sigma = s,c,u$. This shows, see Lemma \ref{L:ExistenceFranksManning}, that $\Phi$ descends to a map $\Phi:X_{\Gamma}\to\mathbb{T}^{d}$. Showing that $\Phi$ gives a fiber bundle structure as in Theorem \ref{Thm:PropertiesOfPHdiffeos} is then similar to \cite{Brin2003}.

Proving accessibility uses a topological argument. Since $X_{\Gamma}$ does not have a virtually abelian fundamental group, and since the kernel of the induced map 
\begin{align}
\Phi_{*} = \pi_{*}:\pi_{1}X_{\Gamma} = \Gamma\to\mathbb{Z}^{d}
\end{align}
is the center of $\Gamma$, there can not exist a connected compact set $K\subset X_{\Gamma}$ such that $\Phi:X_{\Gamma}\supset K\to\mathbb{T}^{d}$ is a finite covering map. This implies, in particular, that $f$ can not have a compact $su-$leaf. So, the proof of accessibility reduces to proving that if $f$ has a non-open accessibility class, then there is a compact $su-$leaf. Obtaining a compact $su-$leaf is done by studying the holonomies between center leaves in the universal cover induced by the fundamental group, as in \cite[map $T_{n}$ defined on page 71]{Rodriguez-Hertz2005}. If $K_{x}\subset G\cap\Hat{W}^{c}(x)$ is the closed set such that $y\in K_{x}$ does not have open accessibility class, then there is an action of $\Gamma$ on $K_{x}$ defined by mapping $y$ to the unique intersection between the accessibility class of $\gamma y$ and $K_{x}$. We show that the induced $\Gamma-$action on the image of $K_{x}$ in $\Hat{W}^{c}(x)/[\Gamma,\Gamma]$ has a fixed point if $K_{x}$ is non-empty. This fixed point corresponds to a compact $su-$leaf, which gives a contradiction so $K_{x}$ must be empty.

\subsubsection{Step 2}

The remainder of the proof of Theorem \ref{Thm:MainTheorem1} follows an idea by Spatzier and Vinhage \cite{SpatzierVinhage2021}: instead of producing the conjugacy directly, we produce a homogeneous structure on $X_{\Gamma}$ that is compatible with $\alpha$. The homogeneous structure on $X_{\Gamma}$ is obtained as the action of a certain quotient of the \textbf{su-path group} (see Section \ref{Sec:suPathGroup}). We use the map $\Phi$ to define the $su$-path group, $\mathcal{P}$, and a natural action of $\mathcal{P}$ on $X_{\Gamma}$. The most technical part of the paper is the following theorem from Section \ref{Sec:InvariantStructureInCenter}.
\begin{theorem}\label{Thm:InvariantStructureInCenter}
The system $f$ has a unique measure of maximal entropy $\mu$ satisfying $\lambda_{\mu}^{c} = 0$.
\end{theorem}
Using Theorem \ref{Thm:InvariantStructureInCenter} we show that there is a normal subgroup, $\mathcal{N}$, of $\mathcal{P}$ such that $N = \mathcal{P}/\mathcal{N}$ is a Nilpotent Lie group that act transitively and freely on $G$. Moreover, the action of $N$ is constructed such that $\alpha$ is compatible with the $N-$action in the sense that the joint action of $N\rtimes\mathbb{Z}^{k}$ is through a semi-direct product. We then use the $N-$action to produce coordinates on $X_{\Gamma}$, in which $f$ is affine. These coordinates gives a bi-Hölder conjugacy $H$ from $\alpha$ to some affine action $\alpha_{0}$.

\subsubsection{Step 3}
 
We finish the proof of Theorem \ref{Thm:MainTheorem1} by proving Theorem \ref{Thm:ImprovedRegularity1}, improving the regularity of $H$ from bi-Hölder to $C^{\infty}$. The proof is similar to the proof in \cite{Rodriguez-HertzWang2014}. We begin the proof by using results from Wilkinson \cite{Wilkinson2008} to show that the conjugacy $H$ is smooth along the center $W^{c}$. The proof of Theorem \ref{Thm:ImprovedRegularity1} then follows \cite{Rodriguez-HertzWang2014} to show that the component of the conjugacy along some coarse exponent $[\chi]$ defining a chamber wall for the chamber that contains $\mathbf{n}_{0}$ (see Section \ref{Sec:BackgroundNilmanifolds}) is smooth. That is, we show that the $[\chi]-$component of the bi-Hölder conjugacy $H$ restricted to $W^{s}(x)$ and normalized by $x\mapsto e$, denoted $H_{x}^{[\chi]}:W^{s}(x)\to G^{[\chi]}$ where $G^{[\chi]}$ is the coarse group with Lie algebra $E_{0}^{[\chi]}$, is uniformly $C^{\infty}$. Once we know that $H_{x}^{[\chi]}$ is smooth, we study the map
\begin{align}
P:{\rm Gr}^{\ell}(E^{s})\to\mathbb{R},\quad P(x,V) = \det\left(D_{x}H_{x}^{[\chi]}|_{V}\right)
\end{align}
where $\ell = \dim(E_{0}^{[\chi]})$ and the determinant is calculated with respect to some background Riemannian metric. The main observation is that $P(x,V)$ can not vanish for all $V\in{\rm Gr}_{x}^{\ell}(E^{s})$ for any $x\in X_{\Gamma}$ (see Lemma \ref{L:ImprovedRegularity8}). This shows that $H_{x}^{[\chi]}$ is a submersion for every $x$, so its fibers form a $C^{\infty}-$foliation within $W^{s}$, denoted $W^{ss}$. Finally, we construct a $\alpha-$invariant distribution $E^{[\chi]}$ transverse to the distribution $E^{ss} = TW^{ss}$ by using a graph transform argument. Existence of the distributions $E^{[\chi]}$ and $E^{ss}$ allows us to produce new partially hyperbolic elements of the action $\alpha$ in a Weyl chamber adjacent to the Weyl chamber containing the first partially hyperbolic element $\mathbf{n}_{0}$. By induction we produce a partially hyperbolic element in every Weyl chamber. Using that $\alpha$ contains many partially hyperbolic elements it follows that $H_{x}^{[\chi]}$ is uniformly $C^{\infty}$ for every coarse exponent $[\chi]$, so $H$ is uniformly $C^{\infty}$ along $W^{s}$ and $W^{u}$. Since $H$ is uniformly $C^{\infty}$ along $W^{s}$, $W^{u}$, and $W^{c}$ we can apply Journé's lemma twice to show that $H$ is $C^{\infty}$.

Theorem \ref{Thm:MainTheorem2} follows from Theorem \ref{Thm:MainTheorem1} and results in \cite{DamjanovicWilkinsonXu2023}.

\subsection{Outline of paper}

In Section \ref{Sec:Background} we go through some of the background results, and basic definitions from partially hyperbolic dynamics and higher rank actions on nilmanifolds. In Section \ref{Sec:PropertiesOfPHdiffeos} we prove Theorem \ref{Thm:PropertiesOfPHdiffeos}. In Section \ref{Sec:suPathGroup} we introduce the $su-$path group, one of the main objects in this paper, and show some of its basic properties. In Section \ref{Sec:HigherRankInvariancePrinciple} we recall the suspension construction of an abelian action and use it, combined with results from \cite{AvilaViana2010}, to derive an invariance principle for higher rank actions on nilmanifolds. Section \ref{Sec:InvariantStructureInCenter} is the most technical part of the paper, here we prove Theorem \ref{Thm:InvariantStructureInCenter}. In Section \ref{Sec:CompatibleAlgebraicStructure} we prove that the action $\alpha$ in Theorem \ref{Thm:MainTheorem1} is topologically conjugated to some affine action. In Section \ref{Sec:ImprovedRegularity} we prove Theorem \ref{Thm:ImprovedRegularity1}, showing that the topological conjugacy is $C^{\infty}$. Finally, in Section \ref{Sec:ProofOfMainTheorems} we complete the proofs of Theorems \ref{Thm:MainTheorem1} and \ref{Thm:MainTheorem2}. We also include an appendix, Appendix \ref{Sec:Appendix}, proving some basic properties of higher rank, abelian algebraic actions on nilmanifolds.

\subsection{Acknowledgements} The author thanks Danijela Damjanović, Homin Lee, Kurt Vinhage, Amie Wilkinson and Disheng Xu for useful discussion.

\section{Background and definitions}\label{Sec:Background}

\subsection{Partially hyperbolic diffeomorphisms}\label{Sec:BackgroundPartiallyHyperbolic}

Let $M$ be a smooth closed manifold and $f\in{\rm Diff}^{\infty}(M)$ a diffeomorphism. We fix a smooth metric $g$ on $M$ inducing a norm $\norm{\cdot}$. We say that $f$ is (absolutely) \textbf{partially hyperbolic} if there is a continuous $Df-$invariant splitting
\begin{align}
T_{x}M = E^{s}(x)\oplus E^{c}(x)\oplus E^{u}(x)
\end{align}
and constants $\nu,\gamma,\widehat{\gamma},\widehat{\nu}\in(0,1)$, $n_{0}\in\mathbb{N}$ such that for $n\geq n_{0}$
\begin{align}
& \norm{Df^{n}|_{E^{s}}}\leq\nu^{n} < \gamma^{n}\leq\norm{\left(Df^{n}\right)^{-1}|_{E^{c}}}^{-1}, \\
& \norm{Df^{n}|_{E^{c}}}\leq\widehat{\gamma}^{-n} < \widehat{\nu}^{-n}\leq\norm{\left(Df^{n}|_{E^{u}}\right)^{-1}}^{-1}.
\end{align}
If we can choose the constants such that
\begin{align}
\nu < \gamma\widehat{\gamma}^{r},\quad\widehat{\nu} < \gamma^{r}\widehat{\gamma}
\end{align}
then we say that $f$ is $r-$\textbf{bunching}. The distributions $E^{s}$, $E^{c}$ and $E^{u}$ are the \textbf{stable, center} and \textbf{unstable} distributions respectively.

Let $f\in{\rm Diff}^{\infty}(M)$ be partially hyperbolic. The stable and unstable distributions are always uniquely integrable to foliations $W^{s}$ and $W^{u}$ with uniformly $C^{\infty}$ leaves, but the center distribution may fail to be integrable. A sufficient condition for $E^{c}$ being integrable is \textbf{dynamical coherence}. We say that $f$ is dynamically coherent if $E^{cs} = E^{c}\oplus E^{s}$ and $E^{cu} = E^{c}\oplus E^{u}$ are both integrable to foliations $W^{cs}$ and $W^{cu}$. In this case we obtain a foliation tangent to $E^{c}$ by intersecting $W^{c} = W^{cs}\cap W^{cu}$. We will denote the distance between two points $p,q\in W^{\sigma}(x)$, $\sigma = s,c,u,cs,cu$, in the leaf metric by $\intd_{\sigma}(p,q)$. Denote the ball about $x$ of radius $\varepsilon$ in $\intd_{\sigma}$ by $W_{\varepsilon}^{\sigma}(x)$. If $f$ is $r-$bunching and dynamically coherent then $W^{cs}$, $W^{cu}$ and $W^{c}$ have uniformly $C^{r}$ leaves \cite{HirschShubPugh1977} (or \cite[Theorem 7]{DamjanovicWilkinsonXu2021}).

Let $\Hat{M}$ be the universal cover of $M$, $\mathcal{F}$ a foliation on $M$, and $\Hat{\mathcal{F}}$ the lifted foliation to $\Hat{M}$.
\begin{definition}\label{Def:QuasiIsometricLeaves}
We say that a continuous foliation with $C^{1}-$leaves $\mathcal{F}$ of $M$ have quasi-isometric leaves in the universal cover if there is a constant $Q\geq 1$ such that
\begin{align}
\intd(x,y)\leq\intd_{\mathcal{F}}(x,y)\leq Q\intd(x,y),\quad x,y\in\Hat{\mathcal{F}}(p),
\end{align}
where $\intd_{\mathcal{F}}$ is the metric along $\mathcal{F}$.
\end{definition}
\begin{remark}
The inequality $\intd(x,y)\leq\intd_{\mathcal{F}}(x,y)$ is immediate since any path connecting $x$ and $y$ along $\Hat{\mathcal{F}}(p)$ also connect $x$ and $y$ in $\Hat{M}$.
\end{remark}
\begin{remark}
We could have asked $\intd_{\mathcal{F}}(x,y)\leq A\intd(x,y) + B$ in the definition, but this is equivalent to Definition \ref{Def:QuasiIsometricLeaves} since $\intd_{\mathcal{F}}$ and $\intd$ are comparable on small balls in $\mathcal{F}$ if $\mathcal{F}$ have uniformly $C^{1}-$leaves.
\end{remark}
In particular, if $f:M\to M$ is partially hyperbolic (and dynamically coherent), we can apply Definition \ref{Def:QuasiIsometricLeaves} to the lifted foliations $\Hat{W}^{\sigma}$, $\sigma = s,c,u,cs,cu$. We can also lift the distributions $E^{\sigma}$ on $M$ to distributions on $\Hat{M}$, also denoted $E^{\sigma}$.

Assume now that $f:M\to M$ is dynamically coherent. Since $E^{s}$ is uniformly transverse to $E^{cu}$, $W^{s}$ and $W^{cu}$ have a local product structure. Similarly, $W^{c}$ and $W^{u}$ subfoliate $W^{cu}$ and $E^{c}$ is transverse to $E^{u}$, so the foliations $W^{c}$ and $W^{u}$ have a local product structure in $W^{cu}$.
\begin{definition}
We say that $f$ have global product structure \cite{Hammerlindl2012} if
\begin{align}
& \#\Hat{W}^{cs}(x)\cap\Hat{W}^{u}(y) = 1,\quad x,y\in\Hat{M}, \\
& \#\Hat{W}^{s}(x)\cap\Hat{W}^{cu}(y) = 1,\quad x,y\in\Hat{M}, \\
& \#\Hat{W}^{s}(x)\cap\Hat{W}^{c}(y) = 1,\quad x,y\in\Hat{W}^{cs}(p), \\
& \#\Hat{W}^{c}(x)\cap\Hat{W}^{u}(y) = 1,\quad x,y\in\Hat{W}^{cu}(p).
\end{align}
\end{definition}
When $f$ has global product structure we define global holonomy maps in the universal cover $\Hat{M}$. Given $x\in\Hat{M}$ and $y\in\Hat{W}^{u}(x)$ we define
\begin{align}
\pi_{x,y}^{u}:\Hat{W}^{cs}(x)\to\Hat{W}^{cs}(y),\quad \{\pi_{x,y}^{u}(z)\} = \Hat{W}^{u}(z)\cap\Hat{W}^{cs}(y).
\end{align}
\begin{figure}[h!]
\begin{tikzpicture}
\draw [color = red, very thick] plot [smooth, tension = 0.6] coordinates {(-4,-1.9) (-1.5,-1.7) (1.5,-1.7) (2.5,-2) (3.5,-2)} node[anchor=west] {\scriptsize$W^{u}(q)$};
\draw [color = black!60!green, very thick] plot [smooth cycle, tension = 0.7] coordinates {(-3,-3) (-2.5,0) (-2,1.5) (-2.5,3) (-3,1.5) (-3.2,0) (-3.4,-3)} node[anchor=east] {\scriptsize$W^{c}(x)$};
\draw [color = black!60!green, very thick] plot [smooth cycle, tension = 0.7] coordinates {(3.1,-3) (3.4,0) (3.9,1.5) (3.5,3) (3.1,1.5) (2.7,0) (2.6,-3)} node[anchor=east] {\scriptsize$W^{c}(y)$};
\draw [color = red, very thick] plot [smooth, tension = 0.6] coordinates {(-2.9,0) (-1.5,0.5) (1.5,0.2) (2.5,-0.3) (5,0)} node[anchor=west] {\scriptsize$W^{u}(p)$};
\draw [color = red, very thick] plot [smooth, tension = 0.6] coordinates {(-2.9,1) (-1.5,1.2) (1.5,0.7) (2.5,1.2) (5,1)} node[anchor=west] {\scriptsize$W^{u}(r)$};
\filldraw [purple] (-2.44,0.24) circle (2pt) node[anchor=north west] {\Tiny$p$};
\filldraw [purple] (-3.37,-1.85) circle (2pt) node[anchor=north west] {\Tiny$q$};
\filldraw [purple] (-2.05,1.17) circle (2pt) node[anchor=north west] {\Tiny$r$};
\filldraw [purple] (3.36,-0.25) circle (2pt) node[anchor=north west] {\Tiny$\pi_{x,y}^{u}(p)$};
\filldraw [purple] (3.83,1.13) circle (2pt) node[anchor=south west] {\Tiny$\pi_{x,y}^{u}(r)$};
\filldraw [purple] (2.57,-2.02) circle (2pt) node[anchor=north east] {\Tiny$\pi_{x,y}^{u}(q)$};
\end{tikzpicture}
\caption{Unstable holonomy between $W^{c}(x)$ and $W^{c}(y)$, $y\in W^{u}(x)$.}\label{Fig:DefinitionOfHolonomies1}
\end{figure}
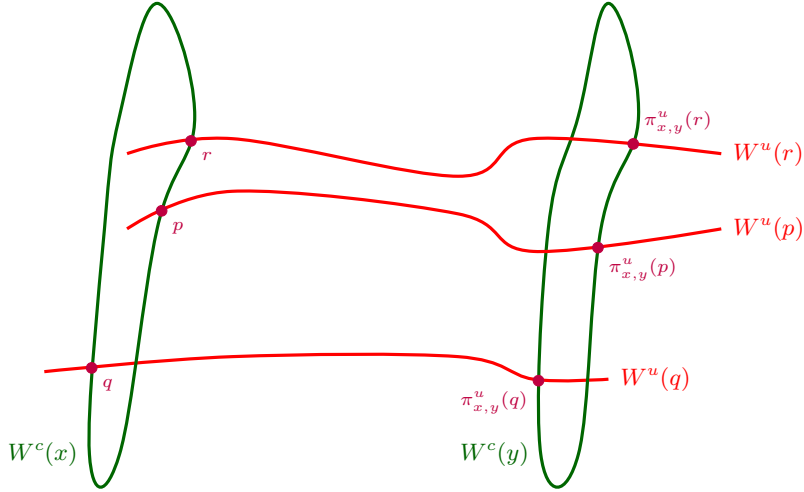
Since $\Hat{W}^{u}$ and $\Hat{W}^{c}$ subfoliate $\Hat{W}^{cu}$ the holonomy maps $\pi_{x,y}^{u}$ restricts to maps $\pi_{x,y}^{u}:\Hat{W}^{c}(x)\to\Hat{W}^{c}(y)$. The holonomy maps between center leaves $\pi_{x,y}^{u}|_{\Hat{W}^{c}(x)}$ descend to holonomy maps between center leaves on $M$ (Figure \ref{Fig:DefinitionOfHolonomies1}). Similarly we define stable holonomies $\pi_{x,y}^{s}:W^{c}(x)\to W^{c}(y)$ when $y\in W^{s}(x)$. When considering holonomies between center leaves then the holonomy maps are $C^{1+\alpha}$ \cite{Brown2022}, and if $f$ is $r-$bunching then the holonomies $\pi_{x,y}^{\sigma}:W^{c}(x)\to W^{c}(y)$, $\sigma = s,u$, are $C^{r}$ \cite{PughShubWilkinson1997}.

We say that a path $\gamma:[0,1]\to M$ is an $su-$path if $[0,1]$ has a subdivision $0 = t_{0} < t_{1} < ... < t_{N-1} < t_{N} = 1$ such that ${\rm Im}\left(\gamma|_{[t_{j},t_{j+1}]}\right)$ is entirely contained in either an $W^{s}-$leaf or a $W^{u}-$leaf. If any two points $x,y\in M$ are connected by an $su-$path, then we say that $f$ is \textbf{accessible}. A set $E\subset M$ is $\sigma-$saturated, $\sigma = s,u$, if $x\in E$ implies $W^{\sigma}(x)\subset E$, and $su-$saturated if it is $s$ and $u-$saturated. Equivalently $f$ is accessible if the only $su-$saturated sets are $M$ and $\emptyset$. For $x\in M$ (or $\Hat{M}$) we define the \textbf{accessibility class} of $x$
\begin{align}
{\rm AC}(x) = \{y\in M\text{ : there is an }su-\text{path connecting }x\text{ and }y\}.
\end{align}
Define a closed set $\Lambda(f)$ by $x\in\Lambda(f)$ if ${\rm AC}(x)$ is not open. That is
\begin{align}
\Lambda(f) = \left[\bigcup_{{\rm AC}(x)\text{ is open}}{\rm AC}(x)\right]^{c} = \bigcap_{{\rm AC}(x)\text{ is open}}{\rm AC}(x)^{c}
\end{align}
If $f$ has $1-$dimensional center direction, then $\Lambda(f)$ is laminated by accessibility classes \cite[Proposition A.3]{RodriguezHertz2006AccessibilityAS}. In particular, if $f$ has $1-$dimensional center and $\Lambda(f) = M$ then $E^{s}\oplus E^{u}$ is jointly integrable to some continuous foliation $W^{su}$ with smooth leaves (in fact, the foliation $W^{su}$ will be a $C^{r}-$foliation if $f$ is $r-$bunching). In the other extreme, $f$ is accessible if and only if $\Lambda(f) = \emptyset$, in this case $f$ has a unique accessibility class.

\subsection{Nilmanifolds and higher rank actions}\label{Sec:BackgroundNilmanifolds}

Let $G$ be a (simply connected) Lie group with Lie algebra $\mathfrak{g}$. We define the lower central series of $\mathfrak{g}$ inductively as 
\begin{align}
\mathfrak{g}^{(1)} = \mathfrak{g},\quad\mathfrak{g}^{(j+1)} = [\mathfrak{g}^{(j)},\mathfrak{g}].
\end{align}
If there is $\ell$ such that $\mathfrak{g}^{(\ell+1)} = 0$ then we say that $\mathfrak{g}$ is nilpotent and the minimal $\ell$ satisfying $\mathfrak{g}^{(\ell+1)} = 0$ is the \textbf{step} of $\mathfrak{g}$. We say that $G$ is a $\ell-$\textbf{step nilpotent Lie group} if $\mathfrak{g}$ is $\ell-$step nilpotent. Given a discrete subgroup $\Gamma\leq G$, we say that $\Gamma$ is a lattice if the quotient space $\Gamma\setminus G$ carries a finite Haar measure $\mu_{\Gamma}$. Equivalently, for nilpotent groups \cite[Corollary 5.4.6]{CorwinGreenleaf1990}, a discrete subgroup $\Gamma\leq G$ is a lattice if the quotient $\Gamma\setminus G$ is compact.

If $G$ is simply connected, nilpotent, and $\Gamma\leq G$ is a lattice then we define the associated \textbf{compact nilmanifold} by 
\begin{align}
X_{\Gamma} = \Gamma\setminus G.
\end{align}
Denote by $\mu_{\Gamma}$ the normalized Haar measure on $X_{\Gamma}$ and write 
\begin{align}
p_{\Gamma}:G\to X_{\Gamma}
\end{align}
for the natural projection map. If $G$ is $\ell-$step then there is a sequence of (normal) subgroups 
\begin{align}
G = G^{(0)}\unrhd G^{(1)}\unrhd...\unrhd G^{(\ell-1)}\unrhd G^{(\ell)} = e,\quad G^{(j)} = \exp\left(\mathfrak{g}^{(j)}\right).
\end{align}
The intersection $\Gamma^{(j)} = G^{(j)}\cap\Gamma$ defines a lattice in $G^{(j)}$ \cite[Theorem 5.2.3]{CorwinGreenleaf1990}, define 
\begin{align}
X_{\Gamma}^{(j)} = X_{\Gamma}/G^{(j)},\quad\pi^{(j)}:X_{\Gamma}\to X_{\Gamma}^{(j)},\quad0\leq j\leq\ell
\end{align}
where $X_{\Gamma}^{(j)}$ is a compact nilmanifold and $\pi^{(j)}:X_{\Gamma}\to X_{\Gamma}^{(j)}$ is a fiber bundle. In particular, if $G$ is $2-$step then we get one (non-trivial) projection map $\pi^{(2)}:X_{\Gamma}\to X_{\Gamma}^{(2)}$. Since $G^{(2)} = [G,G]$, $G/G^{(2)} = G/[G,G]$ is abelian, so $X_{\Gamma}^{(2)}$ is a torus. In the case of $2-$step nilmanifolds, write 
\begin{align}
\pi:X_{\Gamma}\to\mathbb{T}^{d}
\end{align}
for the projection, with $\mathbb{T}^{d}$ the \textbf{base} of $X_{\Gamma}$. The fibers of $\pi$ are $G^{(2)}/\Gamma^{(2)}$ which is also a torus. So $\pi:X_{\Gamma}\to\mathbb{T}^{d}$ is a fiber bundle with base and fibers both tori, but $X_{\Gamma}$ is not a torus.

We define the \textbf{automorphism} and \textbf{affine group} of $X_{\Gamma}$ by
\begin{align}
& {\rm Aut}(X_{\Gamma}) = \{L\in{\rm Aut(G)}\text{ : }L\Gamma = \Gamma\}, \\
& {\rm Aff}(X_{\Gamma}) = \{f_{0}(x) = L(x)g^{-1}\text{ : }L\in{\rm Aut}(X_{\Gamma}),\text{ }g\in G\}.
\end{align}
There is a natural map ${\rm Aff}(X_{\Gamma})\to{\rm Aut}(X_{\Gamma})$ defined by mapping $f_{0}(x) = L(x)g^{-1}$ to the automorphism $L$, and each projection $\pi^{(j)}$ induce a map ${\rm Aut}(X_{\Gamma})\to{\rm Aut}(X_{\Gamma}^{(j)})$.

Fix a homomorphism
\begin{align}
\rho:\mathbb{Z}^{k}\to{\rm Aut}(X_{\Gamma}),\quad\rho^{\mathbf{n}}:X_{\Gamma}\to X_{\Gamma},\text{ }\mathbf{n}\in\mathbb{Z}^{k}.
\end{align}
We say that $\rho$ has a rank$-1$ factor if there is a nilpotent group $\Hat{G}$, of positive dimension, a surjective homomorphism $\phi:G\to\Hat{G}$ such that $\Hat{\Gamma} = \phi\Gamma$ is a lattice in $\Hat{G}$ and an automorphism $L\in{\rm Aut}(\Hat{X}_{\Hat{\Gamma}})$ such that for some finite index subgroup $\Lambda\subset\mathbb{Z}^{k}$ we have $n:\Lambda\to\mathbb{Z}$ satisfying $\phi\rho^{\mathbf{n}} = L^{n(\mathbf{n})}\phi$, $\mathbf{n}\in\Lambda$. That is, $\rho$ has a rank$-1$ factor if there is a factor of $X_{\Gamma}$ where the projected action of $\mathbb{Z}^{k}$ is a $\mathbb{Z}^{1}-$action (up to finite index). More generally, if $\alpha_{0}:\mathbb{Z}^{k}\to{\rm Aff}(X_{\Gamma})$ is a homomorphism then $\alpha_{0}$ has a rank$-1$ factor if
\begin{align}
\mathbb{Z}^{k}\xrightarrow{\alpha_{0}}{\rm Aff}(X_{\Gamma})\xrightarrow{}{\rm Aut}(X_{\Gamma})
\end{align}
has a rank$-1$ factor.
\begin{definition}\label{Def:HigherRankActions}
A homomorphism $\alpha_{0}:\mathbb{Z}^{k}\to{\rm Aff}(X_{\Gamma})$ is higher rank if it has no rank$-1$ factor.
\end{definition}
Let $\alpha_{0}:\mathbb{Z}^{k}\to{\rm Aff}(X_{\Gamma})$ be a homomorphism. We say that $\chi:\mathbb{Z}^{k}\to\mathbb{R}$ is a \textbf{Lyapunov exponent} of $\alpha_{0}$ if there is $v\in\mathfrak{g}\setminus0$ such that
\begin{align}\label{Eq:DefLyapExp}
\chi(\mathbf{n}) = \lim_{\ell\to\pm\infty}\frac{1}{\ell}\log\norm{\alpha_{0}^{\ell\mathbf{n}}(v)}.
\end{align}
The Lyapunov exponents $\chi$ are linear and extends uniquely to $\mathbb{R}^{k}$, we will consider Lyapunov exponents as linear maps on $\mathbb{R}^{k}$. The \textbf{Lyapunov space} associated to $\chi$, $E_{0}^{\chi}\leq\mathfrak{g}$, is the subspace where Equation \ref{Eq:DefLyapExp} hold. For the $0-$functional we write $E_{0}^{c}$. Note that
\begin{align}
\mathfrak{g} = E_{0}^{c}\oplus\bigoplus_{\chi\neq0}E_{0}^{\chi}.
\end{align}
Denote the set of (non-zero) Lyapunov exponents for $\alpha_{0}$ by
\begin{align}
{\rm Lyap}(\alpha_{0}) = \{\chi\neq0\text{ : }\chi\text{ is a Lyapunov exponent of }\alpha_{0}\}.
\end{align}
For $\chi\in{\rm Lyap}(\alpha_{0})$, define the associated \textbf{coarse exponent} and \textbf{coarse space} by
\begin{align}
[\chi] = \{\chi'\in{\rm Lyap}(\alpha_{0})\text{ : }\chi' = c\chi,\text{ for some }c > 0\},\quad E_{0}^{[\chi]} = \bigoplus_{\chi'\in[\chi]}E_{0}^{\chi'}.
\end{align}
If $\chi(\mathbf{n}) > 0$ (or $\chi(\mathbf{n}) < 0$) then $\chi'(\mathbf{n}) > 0$ (or $\chi(\mathbf{n}) < 0$) for every $\chi'\in[\chi]$ so we define $[\chi](\mathbf{n})$ as the sign of $\chi(\mathbf{n})$ (or as $0$ if $\chi(\mathbf{n}) = 0$). We also define $\ker[\chi] = \ker\chi$.
\begin{definition}\label{Def:WeylChambers}
Let $\{[\chi_{1}],...,[\chi_{N}]\}$ be the coarse exponents of $\alpha_{0}$ and 
\begin{align}
U = \left(\bigcup_{j = 1}^{N}\ker[\chi_{j}]\right)^{c}.
\end{align}
Each connected component $\mathcal{C}$ of $U$ is a Weyl chamber of $\alpha_{0}$. The kernels $\ker[\chi]$ are Weyl chamber walls. A wall $\ker[\chi]$ is a wall of $\mathcal{C}$ if $\dim\overline{\mathcal{C}}\cap\ker[\chi] = k-1$.
\end{definition}
Two coarse exponents, $[\chi]$ and $[\eta]$, are \textbf{dependent} if $[\chi](\mathbf{n}) = -[\eta](\mathbf{n})$, otherwise the two exponents are \textbf{independent}. Given any two $\chi',\chi''\in[\chi]$ it is immediate
\begin{align}
\left[E_{0}^{\chi'},E_{0}^{\chi''}\right]\subset E_{0}^{\chi'+\chi''}\quad\left(\text{with }E_{0}^{\chi'+\chi''}=0\text{ if }\chi'+\chi''\not\in{\rm Lyap}(\alpha_{0})\right)
\end{align}
so $E_{0}^{[\chi]}$ is a subalgebra of $\mathfrak{g}$. Define the associated group
\begin{align}
G^{[\chi]}\leq G,\quad G^{[\chi]} = \exp\left(E_{0}^{[\chi]}\right).
\end{align}
In Section \ref{Sec:ImprovedRegularity} we will use that every coarse group $G^{[\chi]}$ has a transverse normal subgroup in the stable group $G^{s}$. More precisely, if $\alpha_{0}^{\mathbf{n}}$ has stable space $E_{0}^{s}$, $[\chi](\mathbf{n}) < 0$, and
\begin{align}
E_{0}^{ss} := \bigoplus_{\substack{[\eta]\neq[\chi] \\
[\eta](\mathbf{n}) < 0}}E_{0}^{[\eta]}
\end{align}
then $E_{0}^{ss}\leq E_{0}^{s}$ is an ideal in $E_{0}^{s}$ \cite[Lemma 3.1]{Rodriguez-HertzWang2014}. Equivalently, the subgroup $G^{ss} = \exp(E_{0}^{ss})\leq G^{s} = \exp(E_{0}^{s})$ is a normal subgroup.

The following two lemmas are well-known, we include proofs in Appendix \ref{Sec:Appendix}.
\begin{lemma}\label{L:Background1}
If $\alpha_{0}:\mathbb{Z}^{k}\to{\rm Aff}(X_{\Gamma})$ is higher rank then there are at least two independent coarse exponents.
\end{lemma}
\begin{lemma}\label{L:Background2}
If $\alpha_{0}:\mathbb{Z}^{k}\to{\rm Aff}(X_{\Gamma})$ is higher rank and $[\chi]$ is a coarse Lyapunov exponent then the space
\begin{align*}
V = \bigoplus_{[\eta]\neq\pm[\chi]}E_{0}^{[\chi]}
\end{align*}
defines a minimal translation action on $X_{\Gamma}$ (the translation action by $V$ is the translation action of the exponential of the Lie algebra generated by $V$).
\end{lemma}
Given a homomorphism $\alpha:\mathbb{Z}^{k}\to{\rm Diff}^{\infty}(X_{\Gamma})$, written $\alpha(\mathbf{n}) = \alpha^{\mathbf{n}}$, we obtain a \textbf{linearization} $\rho:\mathbb{Z}^{k}\to{\rm Aut}(\pi_{1}X_{\Gamma})\cong{\rm Aut}(X_{\Gamma})$.
\begin{definition}
A smooth action $\alpha:\mathbb{Z}^{k}\to{\rm Diff}^{\infty}(X_{\Gamma})$ is higher rank if the linearization $\rho:\mathbb{Z}^{k}\to{\rm Aut}(X_{\Gamma})$ is higher rank.
\end{definition}
Fix $n\geq 1$, $d = 2n$ and define $H^{n} := \mathbb{R}^{n}\times\mathbb{R}^{n}\times\mathbb{R}$. With $g = (q,p,z)\in G$ and $h = (q',p',z')\in G$ we define a multiplication
\begin{align}
gh = (q,p,z)(q',p',z') = \left(q + q', p + p', z + z' + q\cdot p'\right).
\end{align}
This makes $H^{n}$ into a group, the $(d+1)-$\textbf{dimensional Heisenberg group}. Denote by $\omega$ the symplectic form on $\mathbb{R}^{d} = \mathbb{R}^{n}\oplus\mathbb{R}^{n}$. The Lie bracket on $\mathfrak{g}$ is
\begin{align}
\left[(X,Z),(X',Z')\right] = (0,\omega(X,X')),\quad X,X'\in\mathbb{R}^{d}\text{ }Z,Z'\in\mathbb{R}.
\end{align}
Let $\Gamma\leq H^{n}$ be a lattice and $X_{\Gamma} = \Gamma\setminus H^{n}$ the associated nilmanifold. For $L\in{\rm Aut}(X_{\Gamma})$ we obtain a map $L_{su}\in{\rm GL}(d,\mathbb{Z})$ by projecting onto the base, this element $L_{su}$ satisfies $L_{su}\in{\rm Sp}(d,\mathbb{Z})$. In particular, if $[\chi]$ is a coarse exponent of $\alpha_{0}:\mathbb{Z}^{k}\to{\rm Aff}(X_{\Gamma})$ then $-[\chi]$ is also a coarse exponent, so the coarse exponents come in negatively proportional pairs.

In the remainder, we will be interested in groups $G$ of the form $G = \mathbb{R}^{\ell}\times H^{n}$ for some $\ell\geq0$ and $n\geq0$. These groups constitute all abelian simply connected nilpotent groups and all $2-$step, simply connected nilpotent Lie groups with $\dim[G,G] = 1$. Recall the \textbf{Baker–Campbell–Hausdorff formula} \cite{CorwinGreenleaf1990}
\begin{align}
e^{X}e^{Y} = e^{X + Y + [X,Y]/2},\quad X,Y\in\mathfrak{g}.
\end{align}
Fix a left invariant metric, $\intd$, on $G$. Using the Baker–Campbell–Hausdorff formula, it is immediate that for $e^{Z} = g_{c}\in[G,G]$, $\intd(e,g_{c})\leq4\sqrt{\norm{Z}}$.

\section{Some properties of partially hyperbolic diffeomorphisms with quasi isometric leaves in the universal cover}\label{Sec:PropertiesOfPHdiffeos}

In this section, we prove Theorem \ref{Thm:PropertiesOfPHdiffeos}. We begin by proving that $\Phi$ from Theorem \ref{Thm:PropertiesOfPHdiffeos} exists in Section \ref{SubSec:ExistenceFranksManningCoordinate}. In section \ref{SubSec:Accessibility} we show that $f$ is accessible. 

Let $G = H^{n}\times\mathbb{R}^{\ell}$ be the product of some Heisenberg group and some abelian group, $\Gamma\leq G$ a lattice and $X_{\Gamma}$ the associated compact nilmanifold. If $n > 0$ then we write $d = 2n + \ell$ and let $\pi:G\to\mathbb{R}^{d}$ be the base projection. If $n = 0$ and we have an automorphism $L\in{\rm Aut}(X_{\Gamma})$ with $1-$dimensional center, then we let $d = \ell - 1$ and $\pi:G\to\mathbb{R}^{d}$ be the projection from $G$ onto $G/E_{0}^{c}$ (note that if $L$ has $1-$dimensional center then the center direction $E_{0}^{c}$ is a rational line). Assume for the remainder of this section that $f\in{\rm Diff}^{\infty}(X_{\Gamma})$ satisfy all the assumptions of Theorem \ref{Thm:PropertiesOfPHdiffeos}. Denote by $L$ the linearization of $f$ and $L_{su}$ the induced map on the base.

\subsection{Existence of Franks-Manning coordinates}\label{SubSec:ExistenceFranksManningCoordinate}

Write $f:X_{\Gamma}\to X_{\Gamma}$ as
\begin{align}\label{Eq:ExistenceFrankMannin1}
fx = L(x)e^{-v(x)},
\end{align}
with $v:X_{\Gamma}\to\mathfrak{g}$. Fix a lift $F:G\to G$, $Fx = L(x)e^{-v(x)}$. For $x\in G$ let
\begin{align}\label{Eq:ExistenceFrankMannin2}
F^{n}x =: x_{n},\quad n\in\mathbb{Z}.
\end{align}
Let $E_{0}^{s}$, $E_{0}^{c}$ and $E_{0}^{u}$ be the stable, center and unstable direction of $L$ respectively. We decompose any $v\in\mathfrak{g}$ with respect to the splitting $\mathfrak{g} = E_{0}^{s}\oplus E_{0}^{c}\oplus E^{u}$ as $v = v_{s} + v_{c} + v_{u}$. Denote by $\Tilde{\pi}:G\to G/G^{c}\cong\mathbb{R}^{d'}$ the projection, where $G^{c} = \exp(E_{0}^{c})$ is the center of $L$ (we do not know, a priori, that $E_{0}^{c}$ has dimension $1$). Write $A:\mathbb{R}^{d'}\to\mathbb{R}^{d'}$ for the map induced by $L$, then $A$ is hyperbolic (if $\dim E_{0}^{c}\neq 1$ then $A\neq L_{su}$). Recall the following well-known lemma.
\begin{lemma}\label{L:ExistenceFranksManning}
There exists a unique Hölder map $\Phi:G\to G/G^{c}\cong\mathbb{R}^{d'}$
\begin{align}\label{Eq:ExistenceFrankMannin3}
\Phi(x) = \Tilde{\pi}(x) + \varphi(x),\quad\varphi(\gamma x) = \varphi(x),\text{ }\gamma\in\Gamma
\end{align}
such that $\Phi(Fx) = A\Phi(x)$. If $\dim(E_{0}^{c}) = 1$ then $d' = d$, $\Phi:G\to\mathbb{R}^{d}$ descends to a map $\Phi:X_{\Gamma}\to\mathbb{T}^{d}$ homotopic to $\pi$ and $A = L_{su}$.
\end{lemma}
\begin{proof}
The lemma follows from a calculation showing that $\varphi$ satisfy $v_{su}(x) = \varphi(fx) - A(\varphi(x))$, which has a unique solution since $A$ is hyperbolic. If $\dim(E_{0}^{c}) = 1$ then $E_{0}^{c} = [\mathfrak{g},\mathfrak{g}]$ if $G$ is non-abelian (since $[\mathfrak{g},\mathfrak{g}]$ lie in the center of any automorphism) so $G/G^{c}$ is the natural quotient by $[G,G]$. The lemma follows since $\varphi$ is $\Gamma-$invariant. If $G$ is abelian, then $E_{0}^{c}$ is some $1-$dimensional rational line (since we assume that $L_{su}$ have at least one rational eigenvalue) and the last conclusion follows.
\end{proof}
\begin{lemma}\label{L:FrankManningCoordinates1}
If $y\in\Hat{W}^{\sigma}(x)$, $\sigma = s,u$, then $\Phi(y) = \Phi(x)$ if and only if $x = y$. That is $\Phi:\Hat{W}^{\sigma}(x)\to \Phi(x) + E_{0}^{\sigma}$ is injective. Moreover, $\dim E^{\sigma} = \dim E_{0}^{\sigma}$ for $s,c,u$, so $\Phi$ descends and $A = L_{su}$ is hyperbolic.
\end{lemma}
\begin{remark}
We prove the lemma when $G$ is non-abelian. The proof when $G$ is abelian simplifies since all terms from brackets vanish.
\end{remark}
\begin{proof}
Assume that $\sigma = u$, for the other case we reverse time. Let $y\in\hat{W}^{u}(x)$. Write $y_{n} = x_{n}e^{\gamma^{n}}$ where $\gamma^{n}\in\mathfrak{g}$. With respect to the decomposition $\mathfrak{g} = E_{0}^{s}\oplus E_{0}^{c}\oplus E_{0}^{u}$, decompose $\gamma^{n} = \gamma_{s}^{n} + \gamma_{c}^{n} + \gamma_{u}^{n}$. If $\Phi(x) = \Phi(y)$ then $\Phi(x_{n}) = \Phi(y_{n})$ for all $n\geq 0$, so $\norm{\pi(x_{n}) - \pi(y_{n})}\leq2\norm{\varphi}_{C^{0}} =: C$ independently of $n$. On the other hand, $\norm{\pi(x_{n}) - \pi(y_{n})}^{2} = \norm{\gamma_{s}^{n}}^{2} + \norm{\gamma_{u}^{n}}^{2}$. So to show that $\Phi$ is injective on $\Hat{W}^{u}(x)$, it suffices to show that $\gamma_{u}^{n}\to\infty$ as $n\to\infty$ if $x\neq y$.

Suppose for contradiction that $x\neq y$ and $\norm{\gamma_{u}^{n}}\leq K$ uniformly in $n$. From our definitions
\begin{align*}
x_{n+1}e^{\gamma^{n+1}} = & y_{n+1} = F(y_{n}) = L(y_{n})e^{-v(y_{n})} = L(x_{n})e^{L(\gamma^{n})}e^{-v(y_{n})} = \\ &
L(x_{n})e^{L(\gamma^{n}) - v(y_{n}) - [L(\gamma^{n}),v(y_{n})]/2} = \\ &
F(x_{n})e^{v(x_{n}) + L(\gamma^{n}) - v(y_{n}) - \left([L(\gamma^{n}),v(y_{n})] - [v(x_{n}),L(\gamma^{n})] + [v(x_{n}),v(y_{n})]\right)/2}
\end{align*}
or if we take logarithms
\begin{align}\label{Eq:ExistenceFrankMannin4}
\gamma^{n+1} & = v(x_{n}) + L(\gamma^{n}) - v(y_{n}) - \\ &
\label{Eq:ExistenceFrankMannin5}
\frac{[L(\gamma^{n}),v(y_{n})] - [v(x_{n}),L(\gamma^{n})] + [v(x_{n}),v(y_{n})]}{2}.
\end{align}
Using \ref{Eq:ExistenceFrankMannin5} we estimate
\begin{align}\label{Eq:ExistenceFrankMannin6}
& \norm{\gamma_{s}^{n+1}}\leq\norm{L\gamma_{s}^{n}} + C, \\
\label{Eq:ExistenceFrankMannin7}
& \norm{\gamma_{c}^{n+1}}\leq\norm{\gamma_{c}^{n}} + C + k\left(\norm{\gamma_{s}^{n}} + \norm{\gamma_{u}^{n}}\right),
\end{align}
with constants $C,K,k$ that only depend on $\norm{v}_{C^{0}}$. Since $L$ is contracting on $E_{0}^{s}$ there is $\tau\in(0,1)$ such that $\norm{\gamma_{s}^{n+1}}\leq\tau\norm{\gamma_{s}^{n}} + C$, or
\begin{align}\label{Eq:ExistenceFrankMannin8}
\norm{\gamma_{s}^{n}}\leq\frac{C}{1 - \tau}
\end{align}
uniformly in $n$. We have $\norm{\gamma_{u}^{n}}\leq K$ by assumption, so for some possibly larger $C$ we obtain
\begin{align}\label{Eq:ExistenceFrankMannin9}
\norm{\gamma_{c}^{n+1}}\leq\norm{\gamma_{c}^{n}} + C,\quad\norm{\gamma_{c}^{n}}\leq Cn + \norm{\gamma_{c}^{0}}.
\end{align}
After possibly enlarging $C$ again, we have
\begin{align}\label{Eq:ExistenceFrankMannin10}
\intd(x_{n},y_{n}) = \intd(e,e^{\gamma^{n}})\leq\intd(e,e^{\gamma_{s}^{n} + \gamma_{u}^{n}}) + \intd(e,e^{\gamma_{c}^{n}})\leq C(\sqrt{n} + 1).
\end{align}
On the other hand, the assumption that $\Hat{W}^{u}(x_{n})$ is quasi-isometric implies that there is some $\lambda > 1$ and $Q\geq1$ such that
\begin{align}\label{Eq:ExistenceFrankMannin11}
\intd(x_{n},y_{n})\geq\frac{1}{Q}\intd_{u}(F^{n}x,F^{n}y)\geq\frac{1}{Q}\lambda^{n}\intd_{u}(x,y).
\end{align}
If $\intd_{u}(x,y)\neq 0$ then \ref{Eq:ExistenceFrankMannin10} and \ref{Eq:ExistenceFrankMannin11} gives a contradiction for $n$ sufficiently large, so $x = y$.

Since $\Phi:\Hat{W}^{\sigma}(x)\to\Phi(x) + E_{0}^{\sigma}$, $\sigma = s,u$, is injective it follows by invariance of domain that $\dim(E^{\sigma})\leq\dim(E_{0}^{\sigma})$. On the other hand, we have
\begin{align}\label{Eq:ExistenceFrankMannin12}
\dim(E^{s}) + \dim(E^{u}) + 1 = d + 1 = \dim(E_{0}^{s}) + \dim(E_{0}^{c}) + \dim(E_{0}^{u})
\end{align}
or $1\geq\dim(E_{0}^{c})$, so $\dim(E_{0}^{c}) = 1$. This implies $\dim(E^{s}) + \dim(E^{u}) = \dim(E_{0}^{u}) + \dim(E_{0}^{s})$, which only hold if $\dim(E^{\sigma}) = \dim(E_{0}^{\sigma})$, $\sigma = s,u$, since $\dim(E^{\sigma})\leq\dim(E_{0}^{\sigma})$. That $A = L_{su}$ is hyperbolic follows by Lemma \ref{L:ExistenceFranksManning}.
\end{proof}
\begin{remark}
We will make no notational distinction between the map $\Phi:G\to\mathbb{R}^{d}$ and the induced map $\Phi:X_{\Gamma}\to\mathbb{T}^{d}$.
\end{remark}
\begin{lemma}\label{L:FrankManningCoordinates2}
For every $x\in G$ the map $\Phi_{\sigma,x}:\Hat{W}^{\sigma}(x)\to E_{0}^{\sigma}$, defined by
\begin{align}\label{Eq:ExistenceFrankMannin13}
\Phi_{\sigma,x}(y) = \Phi(y) - \Phi(x),
\end{align}
is a homeomorphism. For any $y\in\Hat{W}^{\sigma}(x)$, the map $\Phi_{\sigma,x}:\Hat{W}_{R}^{\sigma}(y)\to E_{0}^{\sigma}$ is uniformly bi-Hölder for fixed $R$.
\end{lemma}
\begin{proof}
By Lemma \ref{L:FrankManningCoordinates1} $\Phi:\Hat{W}^{\sigma}(x)\to\Phi(x) + E_{0}^{\sigma}$ is injective so $\Phi_{\sigma,x}$ is injective. Since $\Hat{W}^{\sigma}(x)$ and $E_{0}^{\sigma}$ have the same dimension, it follows by invariance of domain that $\Phi_{\sigma,x}$ has an open image and is homeomorphic onto its image. In particular, the image of $\Phi_{\sigma,x}$ contain some ball $B_{r_{x}}^{\sigma}(0)$ around $0$ in $E_{0}^{\sigma}$ (to make $r_{x}$ well-defined we take the maximal possible $r_{x}$). Given $\gamma\in\Gamma$ we have 
\begin{align}
\Phi_{\sigma,\gamma x}(\gamma y) = \Phi(\gamma y) - \Phi(\gamma x) = \Phi(y) + \pi(\gamma) - \Phi(x) - \pi(\gamma) = \Phi_{\sigma,x}(y)
\end{align}
from which it follows that $x\mapsto r_{x}$ is $\Gamma-$invariant. Moreover, $\Phi_{\sigma,x}$ and $\Hat{W}^{\sigma}(x)$ vary continuously in $x$, so $r_{x} > r_{0}$ is open. Combined with $\Gamma-$invariance and the fact that $X_{\Gamma}$ is compact, we find $r_{0} > 0$ such that $r_{x}\geq r_{0}$ for all $x\in G$. Assume now $\sigma = u$, the other case follows by reversing time. We have
\begin{align*}
\Phi_{u,x}(\Hat{W}^{u}(x)) = & \Phi_{u,x}(F^{n}\Hat{W}^{u}(x_{-n})) = L^{n}\Phi_{u,x_{-n}}(\Hat{W}^{u}(x_{-n}))\supset L^{n}B_{r_{0}}^{u}(0)
\end{align*}
and letting $n\to\infty$, using that $L$ expand $E_{0}^{u}$, we obtain $\Phi_{u,x}(\Hat{W}^{u}(x)) = E_{0}^{u}$.

So $\Phi_{u,x}$ is a homeomorphism. Since $\Phi$ is Hölder, $(x,y)\mapsto\Phi_{u,x}(y)$ is Hölder in $x$ and $y$. The set $\Phi_{u,x}^{-1}\left(\overline{B_{1}^{u}(0)}\right)$ is compact in $\Hat{W}^{u}(x)$, so we define $K_{x}$ as the minimal radius such that the closure of $\Hat{W}_{K_{x}}^{u}(x)$ contain $\Phi_{u,x}^{-1}\left(\overline{B_{1}^{u}(0)}\right)$. Since $x\mapsto\Phi_{u,x}$ vary continuously the map $x\mapsto K_{x}$ also vary continuosly in $x$. If $K = \sup_{x}K_{x}$, then $K$ is such that if $z,w\in \Hat{W}^{u}(x)$ satisfy $\intd_{u}(z,w)\geq K$ then $\norm{\Phi(z) - \Phi(w)}\geq 1$. Let $\mu = \norm{L|_{E_{0}^{u}}} > 1$ and $\lambda > 1$ such that $\intd_{u}(F^{n}w,F^{n}z)\geq c\lambda^{n}\intd_{u}(z,w)$. Let $z,w\in\Hat{W}^{u}(x)$ satisfy $\intd_{u}(z,w)\leq 1$. There is $n = n(z,w)\geq 0$ such that $n(z,w)\leq-\kappa\log\intd_{u}(z,w) + C$, $\intd_{u}(z_{n},w_{n})\geq K$ and $\kappa$ only depends on $\lambda$. We have
\begin{align*}
\norm{\Phi_{u,x}(z) - \Phi_{u,x}(w)} = & \norm{\Phi(z) - \Phi(w)} = \norm{L^{-n}\left(\Phi(z_{n}) - \Phi(w_{n})\right)}\geq \\ & 
\mu^{-n}\norm{\Phi(z_{n}) - \Phi(w_{n})}\geq\mu^{-n}\geq\mu^{\kappa\log(\intd_{u}(z,w)) - C} = \\ &
\mu^{-C}\intd_{u}(z,w)^{\kappa\log(\mu)}.
\end{align*}
Fix $r > 0$ such that $\Phi_{u,x}^{-1}(B_{r}^{u}(0))\subset\Hat{W}_{1}^{u}(x)$ for all $x$. For $v_{1},v_{2}\in B_{r}^{u}(0)$ we obtain
\begin{align}
\intd_{u}(\Phi_{u,x}^{-1}(v_{1}),\Phi_{u,x}^{-1}(v_{2}))\leq\mu^{C/\kappa\log(\mu)}\norm{v_{1} - v_{2}}^{1/\kappa\log(\mu)}
\end{align}
so $\Phi_{u,x}^{-1}:B_{r}^{u}(0)\to \Hat{W}^{u}(x)$ is uniformly Hölder. Given $y\in \Hat{W}^{u}(x)$ and $v_{0}\in E_{0}^{u}$
\begin{align}
\Phi_{u,y}(\Phi_{u,x}^{-1}(v + v_{0})) = \Phi(\Phi_{u,x}^{-1}(v + v_{0})) - \Phi(y) = v + v_{0} + \Phi(x) - \Phi(y)
\end{align}
or if we choose $\Phi(y) - \Phi(x) = v_{0}$, then
\begin{align}
\Phi_{u,x}^{-1}(v + v_{0}) = \Phi_{u,y}^{-1}(v)
\end{align}
which shows that $\Phi_{u,x}^{-1}$ is uniformly Hölder on any ball of radius $r$. By covering any ball of radius $R$ with balls of radius $r$, the lemma follows.
\end{proof}
\begin{lemma}\label{L:FrankManningCoordinates3}
If $\gamma:[0,1]\to G$ is a $C^{1}-$curve that is tangent to $E^{cs}$ (or $E^{cu}$) then $\Phi(\gamma(1))\in\Phi(\gamma(0)) + E_{0}^{s}$ (or $\Phi(\gamma(1))\in\Phi(\gamma(0)) + E_{0}^{u}$).
\end{lemma}
\begin{proof}
Let $\gamma$ be tangent to $E^{cs}$. Write $H_{x} = \Phi^{-1}(\Phi(x))$. For $y\in H_{x}$ 
\begin{align}
\norm{\pi(y - x)}\leq\norm{\varphi(y) - \varphi(x)}\leq2\norm{\varphi}_{C^{0}} = C
\end{align}
so $H_{x}\subset\pi^{-1}B_{C}(\pi(x))$. In particular, if $y = xe^{\eta(y)}$ for $y\in H_{x}$ then
\begin{align}
\norm{\eta_{s}(y)}^{2} + \norm{\eta_{u}(y)}^{2}\leq C^{2}.
\end{align}
Let $\eta^{n}:\Hat{W}^{u}(x)\to\mathfrak{g}$ be defined by $y_{n} = x_{n}e^{\eta^{n}(y)}$. Since 
\begin{align*}
F^{n}H_{x} = F^{n}\Phi^{-1}(\Phi(x)) = \Phi^{-1}(\Phi(F^{n}x)) = H_{x_{n}}
\end{align*}
we have $\norm{\eta_{s}^{n}(y)},\norm{\eta_{u}^{n}(y)}\leq C$ uniformly in $n\in\mathbb{Z}$. As in the proof of Lemma \ref{L:FrankManningCoordinates1} it follows that we also have $\norm{\eta_{c}^{n}(y)}\leq C|n|$ for some (possibly larger) constant $C$. Let $\gamma:[0,1]\to G$ be a $C^{1}-$curve tangent to $E^{cs}$ and denote the end-points of $\gamma$ by $x = \gamma(0)$ and $y = \gamma(1)$. We find $z\in\Hat{W}^{s}(y)$ and $w\in\Hat{W}^{u}(z)$ such that $w\in H_{x}$ (this corresponds to choosing a two-legged $su-$path from $\Phi(y)$ to $\Phi(x)$ in $\mathbb{R}^{d}$). By the reverse triangle inequality
\begin{align}\label{Eq:TriangleInequality}
\intd(x_{n},w_{n})\geq\intd(w_{n},z_{n}) - \intd(x_{n},y_{n}) - \intd(y_{n},z_{n}).
\end{align}
Since $w$ and $z$ lie in the same unstable leaf and $\Hat{W}^{u}$ have quasi-isometric leaves, $\intd(w_{n},z_{n})\geq c\lambda^{n}\intd(w,z)$ for some $c > 0$ and $\lambda > 1$. But $\gamma$ is a $C^{1}-$curve along $E^{cs}$ from $x$ to $y$, so
\begin{align}
\intd(x_{n},y_{n})\leq\int_{0}^{1}\norm{D_{\gamma(t)}F^{n}(\Dot{\gamma}(t))}\intd t\leq C\Hat{\lambda}^{n}\cdot|\gamma|
\end{align}
where $|\gamma|$ is the length of $\gamma$ and $\Hat{\lambda} < \lambda$. Finally, $\intd(y_{n},z_{n})\leq C$ uniformly for $n\geq 0$ since $y$ and $z$ lie in the same stable leaf. Equation \ref{Eq:TriangleInequality} implies
\begin{align}
c\lambda^{n}\intd(w,z) - C\Hat{\lambda}^{n}|\gamma| - C\leq\intd(x_{n},w_{n}) = \intd(e,e^{\eta^{n}(w)})\leq C\sqrt{n},\quad n\geq0
\end{align}
for some constant $C$. If $\intd(w,z) > 0$ then we obtain a contradiction for $n$ large enough. We conclude that $\intd(w,z) = 0$. That is, $w = z$ so $\Hat{W}^{s}(y)\cap H_{x} = \{w\}\neq\emptyset$, which implies  
\begin{align}
\Phi(x) - \Phi(y) = \Phi(w) - \Phi(y) = \Phi_{s,y}(w)\in E_{0}^{s}
\end{align}
proving the lemma for $sc$. The lemma follows for $cu$ by reversing time.
\end{proof}
We can now prove the first five points of Theorem \ref{Thm:PropertiesOfPHdiffeos}
\begin{proof}[Proof of Theorem \ref{Thm:PropertiesOfPHdiffeos}]
Dynamical coherence of $f$ follows from \cite{Brin2003}. Any curve $\gamma$ tangent to $E^{cs}$ satisfies $\Phi(\gamma(1))\in\Phi(\gamma(0)) + E_{0}^{s}$ by Lemma \ref{L:FrankManningCoordinates3}, so $\Hat{W}^{cs}(x)\subset\Phi^{-1}(\Phi(x) + E_{0}^{s})$ and $\Phi^{-1}(\Phi(x) + E_{0}^{s})$ is a union of $\Hat{W}^{cs}-$leaves. We claim that $\Phi^{-1}(\Phi(x) + E_{0}^{s})$ is path-connected, which proves that $\Phi^{-1}(\Phi(x) + E_{0}^{s}) = \Hat{W}^{cs}(x)$. Given any $y\in G$ there is a unique intersection $\Hat{W}^{u}(y)\cap\Phi^{-1}(\Phi(x) + E_{0}^{s}) = \{w\}$ (since $\Hat{W}^{u}(y)$ maps homeomorphically onto $\Phi(y) + E_{0}^{u}$ under $\Phi$). Since $\Phi^{-1}(\Phi(x) + E_{0}^{s})$ is tangent to $E^{cs}$ and $\Hat{W}^{u}$ is tangent to $E^{u}$ the map 
\begin{align}
G\to\Phi^{-1}(\Phi(x) + E_{0}^{s}),\quad y\mapsto \Hat{W}^{u}(y)\cap\Phi^{-1}(\Phi(x) + E_{0}^{s})
\end{align}
is continuous. Since $G\to\Phi^{-1}(\Phi(x) + E_{0}^{s})$ is surjective and $G$ is path-connected it follows that $\Phi^{-1}(\Phi(x) + E_{0}^{s})$ is path-connected. Properties $(i)$ and $(ii)$ follows. Using 
\begin{align*}
& \Hat{W}^{cs}(x) = \Phi^{-1}(\Phi(x) + E_{0}^{s}), \\
& \Hat{W}^{cu}(x) = \Phi^{-1}(\Phi(x) + E_{0}^{u}),\text{ and} \\
& \Hat{W}^{c}(x) = \Hat{W}^{cs}(x)\cap\Hat{W}^{cu}(x)
\end{align*}
we obtain $\Hat{W}^{c}(x) = \Phi^{-1}(\Phi(x))$. Since $\Phi$ descends to a map $\Phi:X_{\Gamma}\to\mathbb{T}^{d}$ the fibers $W^{c}(x) = \Phi^{-1}(\Phi(x))$ are compact. We define
\begin{align}\label{Eq:HomotopyFixedFundamentalGroup}
\Gamma^{c} = \{\gamma\in\Gamma\text{ : }L\gamma = \gamma\}
\end{align}
where $L\in{\rm Aut}(X_{\Gamma})$ is the linearization of $f$. For $\gamma^{c}\in\Gamma^{c}$ we have 
\begin{align}
\gamma^{c}\Hat{W}^{c}(x) = \Hat{W}^{c}(\gamma^{c}x) = \Phi^{-1}(\Phi(\gamma^{c}x)) = \Phi^{-1}(\Phi(x)) = \Hat{W}^{c}(x)
\end{align}
where we have used $\Phi(\gamma x) = \pi(\gamma) + \Phi(x)$ and $\ker\pi = \Gamma^{c}$. So, if we fix a generator $\gamma_{0}^{c}$ of $\Gamma^{c}$ then we can orient $\Hat{W}^{c}(x)$ by letting $\gamma_{0}^{c}x > x$. This is a well-defined orientation of $\Hat{W}^{c}(x)$ since $x\mapsto\gamma_{0}^{c}x$ have no fixed points. That $W^{c}$ are circles follows since they are compact $1-$dimensional manifolds.

Let $\varepsilon > 0$ be small and fix $x\in X_{\Gamma}$. We denote by $D = W_{\varepsilon}^{s}(W_{\varepsilon}^{u}(x))$, $U = W^{c}(D)$ and $D' = \Phi(D) = \Phi(U)$. With $\varepsilon$ small enough we define 
\begin{align}
\pi^{su}:U\to W^{c}(x)
\end{align}
by the unique holonomy first along $W^{s}$ then $W^{u}$ in $U$. We obtain a map
\begin{align}
U\to D'\times W^{c}(x)\quad y\mapsto(\Phi(y),\pi^{su}(y)),
\end{align}
this map is smooth along $W^{c}$ since the holonomies are $C^{r}$ (where $r$ depends on the bunching). That $\Phi$ semiconjugates (a finite index subgroup of) the centralizer $Z^{\infty}(f)$ onto its linearization is immediate since $\Phi$ is (essentially) unique homotopic to $\pi$ (note that the uniqueness in Lemma \ref{L:ExistenceFranksManning} implies that $\Phi:X_{\Gamma}\to\mathbb{T}^{d}$ is unique modulo the fact that we could change $\Phi(x)$ to $\Phi(x) + p_{0}$ where $L_{su}p_{0} = p_{0}$). Indeed, for any $g\in Z^{\infty}(g)$ we let $B\in{\rm GL}(d,\mathbb{Z})$ be the induced map on $\mathbb{T}^{d}$ by its linearization and let $\Tilde{\Phi}(x) = B^{-1}\Phi(gx)$. Then $\Tilde{\Phi}$ is homotopic to $\pi$ and we still have
\begin{align}
\Tilde{\Phi}(fx) = B^{-1}\Phi(gfx) = B^{-1}\Phi(fgx) = L_{su}B^{-1}\Phi(gx) = L_{su}\Tilde{\Phi}(x)
\end{align}
so $\Tilde{\Phi}(x) = \Phi(x) + p_{0}$ where $p_{0}$ is a fixed point for $L_{su}$. Define
\begin{align}
Z_{\rm fix}^{\infty}(f) = \{g\in Z^{\infty}(f)\text{ : }gW^{c}(x) = W^{c}(x)\text{ if }fW^{c}(x) = W^{c}(x)\}.
\end{align}
Since $\Phi(fx) = L_{su}\Phi(x)$ and $L_{su}$ has finitely many fixed points $Z_{\rm fix}^{\infty}(f)$ has finite index in $Z^{\infty}(f)$. For $g\in Z_{\rm fix}^{\infty}(f)$ and $\Phi(x_{0}) = 0$ we have $fW^{c}(x_{0}) = W^{c}(x_{0})$ so $gW^{c}(x_{0}) = W^{c}(x_{0})$. It follows that
\begin{align}
p_{0} = \Tilde{\Phi}(x_{0}) = B^{-1}\Phi(gx_{0}) = B^{-1}\Phi(x_{0}) = 0
\end{align}
or $B^{-1}\Phi(gx) = \Phi(x)$.
\end{proof}

\subsection{Proof of accessibility}\label{SubSec:Accessibility}

In this section, we show the last point of Theorem \ref{Thm:PropertiesOfPHdiffeos}. For $\gamma\in\Gamma$, define $T_{\gamma}:G\to G$ by
\begin{align}
T_{\gamma}(x) = \Hat{W}^{u}\left[\Hat{W}^{s}(\gamma x)\cap\Hat{W}^{cu}(x)\right]\cap\Hat{W}^{c}(x) = \Hat{W}^{u}(\Hat{W}^{s}(\gamma x))\cap \Hat{W}^{c}(x)
\end{align}
which is well-defined by point $(i)$ of Theorem \ref{Thm:PropertiesOfPHdiffeos}. A calculation shows
\begin{align*}
FT_{\gamma}(x) = T_{L\gamma}(Fx).
\end{align*}
Recall that we denote by $\Lambda(F)\subset G$ and $\Lambda(f)\subset X_{\Gamma}$ the complement of open accessibility classes. The set $\Lambda(F)$ is closed and $su-$saturated. Since $E^{c}$ is $1-$dimensional the set $\Lambda(F)$ is laminated by accessibility classes, denoted ${\rm AC}(x) = \Hat{W}^{su}(x)$, \cite{RodriguezHertz2006AccessibilityAS}. The union of open accessibility classes is $\Gamma-$invariant, so $\Lambda(F)$ is $\Gamma-$invariant. Moreover, $F$ maps accessibility classes to accessibility classes so $F\Lambda(F) = \Lambda(F)$. Given $\gamma,\gamma'\in\Gamma$ and $x\in\Lambda(F)$ we have
\begin{align*}
T_{\gamma}(T_{\gamma'}(x)) = & \Hat{W}^{su}\left(\gamma\left[\Hat{W}^{su}(\gamma'x)\cap \Hat{W}^{c}(x)\right]\right))\cap\Hat{W}^{c}\left(\left[\Hat{W}^{su}(\gamma'x)\cap \Hat{W}^{c}(x)\right]\right) = \\ &
\Hat{W}^{su}(\Hat{W}^{su}(\gamma\gamma'x)\cap \Hat{W}^{c}(\gamma x))\cap \Hat{W}^{c}(x) = \\ &
\Hat{W}^{su}(\gamma\gamma'x)\cap\Hat{W}^{c}(x) = T_{\gamma\gamma'}(x)
\end{align*}
so restricted to $\Lambda(F)$ the map $(\gamma,x)\mapsto T_{\gamma}(x)$ defines a group action of $\Gamma$, see also \cite[Lemma 6.1]{Rodriguez-Hertz2005}. Before starting the proof we will need an elementary, but important, auxiliary lemma on $\mathbb{Z}^{k}-$actions on the circle.
\begin{lemma}\label{L:AbCHasFixedPoints}
Let $f_{1},...,f_{k},g\in{\rm Homeo}_{+}(\mathbb{T})$ be orientation preserving homeomorphisms on the circle and let $K\subset\mathbb{T}$ be a compact subset that is invariant by $g$ and each $f_{j}$. Moreover, assume that $f_{i}f_{j} = f_{j}f_{i}$ on $K$ and that there is some hyperbolic integer matrix $(A_{i}^{j})_{1\leq i,j\leq k}$ such that
\begin{align}\label{Eq:AbelianByCyclicCondition}
gf_{i} = f_{1}^{A_{i}^{1}}...f_{k}^{A_{i}^{k}}g,\quad{\rm on}\text{ }K.
\end{align}
Then the $\mathbb{Z}^{k}-$action generated by $f_{1},...,f_{k}$ on $K$, $\beta:\mathbb{Z}^{k}\times K\to K$, has a periodic point. That is, there is a point $p\in K$ and a finite index subgroup $\Lambda\leq\mathbb{Z}^{k}$ such that $\beta^{\mathbf{n}}p = p$ for $\mathbf{n}\in\Lambda$.
\end{lemma}
\begin{remark}
The condition in Equation \ref{Eq:AbelianByCyclicCondition} says that $\beta:\mathbb{Z}^{k}\times K\to K$ joint with $g$ form an Abelian-by-Cyclic (AbC) action on $K$.
\end{remark}
In the proof, we will use the following two lemmas whose proofs are standard, but we include them for completeness.
\begin{lemma}\label{L:FixedPointsInCptInvSets}
If $f\in{\rm Homeo}_{+}(\mathbb{T})$ is an orientation preserving homeomorphism on the circle with zero rotation number, $\omega(f) = 0$, and $K\subset\mathbb{T}$ is a compact $f-$invariant set, then $f$ have a fixed point in $K$.
\end{lemma}
\begin{proof}
If $f$ has zero rotation number then for any $p\in\mathbb{T}$ the sequence $f^{n}p$ converges to a fixed point of $f$ as $n\to\infty$. For any $x\in K$ the sequence $f^{n}x$ lie in $K$ since $K$ is $f-$invariant. By compactness of $K$ any limit point of $f^{n}x$ also lies in $K$, so $f$ has a fixed point in $K$.
\end{proof}
\begin{lemma}\label{L:AdditivityOfRotationNumber}
Let $f,g\in{\rm Homeo}_{+}(\mathbb{T})$ be orientation preserving homeomorphisms of the circle and assume that there is a $f$ and $g-$invariant probability measure $\nu$. Then the rotation numbers satisfies $\omega(fg) = \omega(f) + \omega(g)$.
\end{lemma}
\begin{proof}
Write $f(x) = x + u(x)$, $g(x) = x + v(x)$ where $u,v:\mathbb{T}\to\mathbb{R}$. Recall that if $\mu$ is a $f-$invariant measure then we obtain the rotation number of $f$ as
\begin{align}
\omega(f) = \int_{\mathbb{T}}u(x)\intd\mu(x) + \mathbb{Z}.
\end{align}
Similarly we obtain the rotation number of $g$. The measure $\nu$ is $fg-$invariant and $f(gx) = x + u(gx) + v(x)$, so we can write the rotation number of $fg$ as
\begin{align*}
\omega(fg) = & \int_{\mathbb{T}}\left[u(gx) + v(x)\right]\intd\nu(x) + \mathbb{Z} = \\ &
\int_{\mathbb{T}}u(gx)\intd\nu(x) + \int_{\mathbb{T}}v(x)\intd\nu(x) + \mathbb{Z} = \\ &
\int_{\mathbb{T}}u(x)\intd\nu(x) + \int_{\mathbb{T}}v(x)\intd\nu(x) + \mathbb{Z} = \\ &
\omega(f) + \omega(g)
\end{align*}
where the second to last equality use that $\nu$ is $g-$invariant.
\end{proof}
\begin{proof}[Proof of Lemma \ref{L:AbCHasFixedPoints}]
Since the $\beta-$action is abelian on $K$ it has an invariant measure on $K$. By Lemma \ref{L:AdditivityOfRotationNumber} the rotation numbers satisfy 
\begin{align}\label{Eq:AdditivityRotationNumber}
\omega(f_{1}^{n_{1}}...f_{k}^{n_{k}}) = n_{1}\omega(f_{1})+...+n_{k}\omega(f_{k})
\end{align}
for all integers $n_{1},...,n_{k}\in\mathbb{Z}$. Since $g$ preserves orientation, conjugacy invariance of rotation number and Equation \ref{Eq:AbelianByCyclicCondition} implies
\begin{align}
\omega(f_{i}) = \omega(f_{1}^{A_{i}^{1}}) + ... + \omega(f_{k}^{A_{i}^{k}}).
\end{align}
Or if we denote the map $f_{1}^{n_{1}}...f_{k}^{n_{k}}$ by $\beta^{\mathbf{n}}$, then we can write $\omega(\beta^{\mathbf{n}}) = \omega(\beta^{A\mathbf{n}})$. Using Equation \ref{Eq:AdditivityRotationNumber} we obtain
\begin{align}
\omega\left(\beta^{(A - I)\mathbf{n}}\right) = 0
\end{align}
for all $\mathbf{n}\in\mathbb{Z}^{k}$. Since $A$ is hyperbolic $A - I$ is invertible over the rationals. So there is a finite index subgroup $\Lambda\subset\mathbb{Z}^{k}$ such that $\omega(\beta^{\mathbf{n}}) = 0$ for all $\mathbf{n}\in\Lambda$. Let $e_{1},...,e_{k}\in\Lambda$ be generators. By Lemma \ref{L:FixedPointsInCptInvSets} the map $\beta^{e_{1}}$ has a fixed point in $K$. Since $\beta^{e_{2}}$ commute with $\beta^{e_{1}}$ within $K$ it follows that $\beta^{e_{2}}$ preserve the compact set ${\rm Fix}(\beta^{e_{1}})\cap K\neq\emptyset$. So, if we apply Lemma \ref{L:FixedPointsInCptInvSets} once more we see that the set
\begin{align}
K\cap{\rm Fix}(\beta^{e_{1}})\cap{\rm Fix}(\beta^{e_{2}})
\end{align}
is non-empty. Proceeding by induction, we find a point $p\in K$ that is fixed by $\beta^{e_{1}},...,\beta^{e_{k}}$, and therefore by $\Lambda$. Since $\Lambda$ has a finite index in $\mathbb{Z}^{k}$ the lemma follows.
\end{proof}
Recall that $L\in{\rm Aut}(X_{\Gamma})$ is the linearization of $f$ and the $L-$fixed part of $\Gamma$ is $\Gamma^{c}$ (Equation \ref{Eq:HomotopyFixedFundamentalGroup}).
\begin{lemma}\label{L:TrivialityOfHolonomyAction1}
Assume that $\Lambda(f)\neq\emptyset$. There is a finite index subgroup $\Gamma'\leq\Gamma$ and $x\in\Lambda(F)$ such that $T_{\gamma}(x)\in\Gamma^{c}x$ for $\gamma\in\Gamma'$.
\end{lemma}
\begin{proof}
Fix generators $\gamma_{1},\gamma_{2},...,\gamma_{d}\in\Gamma$. Let $x_{0}\in G$ be such that $\Phi(x_{0}) = e$. Identify $\Hat{W}^{c}(x_{0})/\Gamma^{c}\cong W^{c}(p_{\Gamma}(x_{0}))$ with $\mathbb{T}$. Since $\Gamma^{c}$ is central in $\Gamma$ we can identify $T_{\gamma_{j}}:\Hat{W}^{c}(x_{0})/\Gamma^{c}\to\Hat{W}^{c}(x_{0})/\Gamma^{c}$ with circle diffeomorphisms. We also identify $F$ with a circle diffeomorphism (by our choice of $x_{0}$ and the fact that $\Gamma^{c}$ is $L-$fixed we have $F\Hat{W}^{c}(x_{0})/\Gamma^{c} = \Hat{W}^{c}(x_{0})/\Gamma^{c}$). Let $K = (\Lambda(F)\cap\Hat{W}^{c}(x_{0}))/\Gamma^{c}$, which is compact, $T_{\gamma_{j}}-$invariant, $F-$invariant, and non-empty. Any $\gamma^{c}\in\Gamma^{c}$ act trivially on $\mathbb{T}$ (and therefore $K$) under $T_{\gamma^{c}}$, so the action $T_{\gamma}$ on $K$ factor through $\Gamma/\Gamma^{c}\cong\mathbb{Z}^{d}$. Moreover, $F$ satisfy $FT_{\gamma} = T_{L\gamma}F$ so the assumptions of Lemma \ref{L:AbCHasFixedPoints} are satisfied with $T_{\gamma_{j}} = f_{j}$ and $g = F$. Therefore there is a finite index subgroup of $\Gamma/\Gamma^{c}$ that admits a fixed point on $\mathbb{T}$ which implies the lemma.
\end{proof}
An immediate corollary of Lemma \ref{L:TrivialityOfHolonomyAction1} is that there is a compact $su-$leaf.
\begin{lemma}
If $\Lambda(f)\neq\emptyset$, or equivalently if $f$ is not accessible, then there is a compact $su-$leaf intersecting each center leaf $q < \infty$ times.
\end{lemma}
\begin{proof}
For any $x\in\Lambda(F)$ the map
\begin{align}
\Hat{W}^{su}(x)\xrightarrow{\Phi}\mathbb{R}^{d}
\end{align}
is a homeomorphism by points $(i)$ and $(iv)$ in Theorem \ref{Thm:PropertiesOfPHdiffeos}. Choose $x_{0}$ as in Lemma \ref{L:TrivialityOfHolonomyAction1} and let $\Gamma'\leq\Gamma$ be of finite index such that $\gamma\Hat{W}^{su}(x_{0})\subset\Gamma^{c}\Hat{W}^{su}(x_{0})$ for all $\gamma\in\Gamma'$. Note that $\Phi_{*}\Gamma' = \pi_{*}\Gamma' = \Lambda\subset\mathbb{Z}^{d}$ has finite index in $\mathbb{Z}^{d}$. Define
\begin{align}\label{Eq:TransverseToDerivedSubgroup}
\Lambda_{*} = \{\gamma\in\Gamma\text{ : }\gamma\Hat{W}^{su}(x_{0}) = \Hat{W}^{su}(x_{0})\}\subset\Gamma'.
\end{align}
Since $\gamma\Hat{W}^{su}(x_{0})\subset\Gamma^{c}\Hat{W}^{su}(x_{0})$ for all $\gamma\in\Gamma'$ there is for each $\gamma\in\Gamma'$ some $\gamma^{c}\in\Gamma^{c}$ such that $\gamma\gamma^{c}\in\Lambda_{*}$. In particular, $\Phi_{*}\Lambda_{*} = \Lambda$. Since $\Lambda$ is a lattice in $\mathbb{R}^{d}$ and we obtain a homeomorphism $\Hat{W}^{su}(x_{0})/\Lambda_{*}\cong\mathbb{R}^{d}/\Lambda$ the set $\Hat{W}^{su}(x_{0})/\Lambda_{*}$ is compact. Since $p_{\Gamma}\Hat{W}^{su}(x_{0}) = W^{su}(p_{\Gamma}(x_{0}))$ is homeomorphic to $\Hat{W}^{su}(x_{0})/\Lambda_{*}$ the lemma follows.
\end{proof}
We can now prove the last claim of Theorem \ref{Thm:PropertiesOfPHdiffeos}
\begin{proof}[Proof of $(vi)$ in Theorem \ref{Thm:PropertiesOfPHdiffeos}]
If $f$ is not accessible then we construct $\Lambda_{*}\subset\Gamma$ as in Equation \ref{Eq:TransverseToDerivedSubgroup}. Since $\pi_{*}\Lambda_{*}$ has finite index in $\mathbb{Z}^{d}$ the group $\Lambda_{*}\times\ker\pi_{*}$ has finite index in $\Gamma$, so $\Gamma$ is virtually abelian. This is a contradiction if $G$ is a non-abelian nilpotent Lie group.
\end{proof}

\section{Action of the su-path group}\label{Sec:suPathGroup}

In this section we introduce and prove basic properties of the $su-$path group $\mathcal{P}$. The $su-$path group naturally acts on $X_{\Gamma}$ (Definition \ref{Def:suPathGroupAction}). The group $\mathcal{P}$, its various subgroups, and its action on $X_{\Gamma}$ will be the key object in the following sections.

\subsection{The su-path group}

Let $\alpha:\mathbb{Z}^{k}\times X_{\Gamma}\to X_{\Gamma}$ be a smooth action satisfying the assumptions of Theorem \ref{Thm:MainTheorem1}. Moreover, let $\rho:\mathbb{Z}^{k}\to{\rm GL}(d,\mathbb{Z})$ the action defined by $\rho^{\mathbf{n}}\Phi = \Phi\alpha^{\mathbf{n}}$ from Theorem \ref{Thm:PropertiesOfPHdiffeos}, $f = \alpha(\gamma_{0})$ the partially hyperbolic element, and let $L$ be the linearization of $f$. Denote by $E_{0}^{s}$ and $E_{0}^{u}$ the stable and unstable distributions for $L$.
\begin{definition}\label{Def:suPathGroupAction}
We define the $su-$path group $\mathcal{P}$ as the free product
\begin{align}
\mathcal{P} = E_{0}^{s}*E_{0}^{u}.
\end{align}
If $w\in\mathcal{P}$ is a word in $\mathcal{P}$ then we define $\Pi(w) = w_{1}+...+w_{N}$ to be the sum of all factors in $w$. We define the normal subgroup $\mathcal{P}^{c} = \Pi^{-1}(0)$.
\end{definition}
Given any pair of negatively proportional course Lyapunov spaces, $E_{0}^{-[\chi]}$, $E_{0}^{[\chi]}$, of $\rho$, we define the $[\chi]-$\textbf{path group}
\begin{align}
\mathcal{P}_{[\chi]} = E_{0}^{-[\chi]}*E_{0}^{[\chi]},\quad\mathcal{P}_{[\chi]}^{c} = \mathcal{P}_{[\chi]}\cap\Pi^{-1}(0)
\end{align}
with $E_{0}^{-[\chi]} = 0$ if $-[\chi]$ is not a coarse exponent. We also define the \textbf{complementary $[\chi]-$path group}
\begin{align}
\mathcal{Q}_{[\chi]} = \left(\bigoplus_{\substack{[\eta]\neq\pm[\chi] \\
[\eta](\mathbf{n}_{0}) > 0}}E^{-[\eta]}\right)*\left(\bigoplus_{\substack{[\eta]\neq\pm[\chi] \\
[\eta](\mathbf{n}_{0}) > 0}}E^{[\eta]}\right),\quad\mathcal{Q}_{[\chi]}^{c} = \mathcal{Q}_{[\chi]}\cap\mathcal{P}^{c}.
\end{align}
It is immediate that $\mathcal{P}_{[\chi]},\mathcal{Q}_{[\chi]}\subset\mathcal{P}$, $\mathcal{P}_{[\chi]}^{c},\mathcal{Q}_{[\chi]}^{c}\subset\mathcal{P}^{c}$. The following well-known lemma on free products will be useful.
\begin{lemma}\label{L:CommutatorOfFreeProduct}
Let $V,W$ be vector spaces, $\mathcal{G} = V*U$ and $\Pi:G\to V\oplus W$ the map defined by 
\begin{align}
\Pi(v_{1}u_{1}...v_{n}u_{n}) = v_{1}+...+v_{n}+u_{1}+...+u_{n}.
\end{align}
Then $\Pi^{-1}(0) = \ker\Pi = [\mathcal{G},\mathcal{G}]$ and any $w\in\mathcal{G}$ can be written $w = \Tilde{w}vu$ with $v\in V$, $u\in U$ and $\Tilde{w}\in[\mathcal{G},\mathcal{G}]$.
\end{lemma}
\begin{proof}
For $w = v_{1}u_{1}...v_{n}u_{n}\in\mathcal{G}$ we have
\begin{align*}
v_{1}u_{1}...v_{n}u_{n} = & \left[v_{1}u_{1}(-v_{1})(-u_{1})\right]u_{1}(v_{1}+v_{2})u_{2}...v_{n}u_{n} = \\ &
\left[v_{1}u_{1}(-v_{1})(-u_{1})\right]\cdot\left[u_{1}(v_{1}+v_{2})(-u_{1})(-v_{1}-v_{2})\right]\cdot \\ & 
(v_{1}+v_{2})(u_{1}+u_{2})v_{3}u_{3}...v_{n}u_{n}.
\end{align*}
Note that $v_{1}u_{1}(-v_{1})(-u_{1}),u_{1}(v_{1}+v_{2})(-u_{1})(-v_{1}-v_{2})\in[\mathcal{G},\mathcal{G}]$ so $(v_{1}+v_{2})(u_{1}+u_{2})v_{3}u_{3}...v_{n}u_{n}$ consists of $n-1$ pairs. By induction we find $w_{1},...,w_{\ell}\in[\mathcal{G},\mathcal{G}]$ such that
\begin{align}
w = w_{1}...w_{\ell}\cdot\left(v_{1}+...+v_{n}\right)\left(u_{1}+...+u_{n}\right).
\end{align}
This proves the last part of the lemma. Since $\Pi(w_{j}) = 0$ for each $j$, we have $\Pi(w) = v_{1}+...+v_{n} + u_{1}+...+u_{n}$. In particular, if $w\in\ker\Pi$ then $v_{1}+...+v_{n} = 0$ and $u_{1}+...+u_{n} = 0$, so $w = w_{1}...w_{\ell}\in[\mathcal{G},\mathcal{G}]$.
\end{proof}
\begin{definition}\label{Def:LiftedTranslations}
For $t\in E_{0}^{\sigma}$ we define $\eta_{\sigma}^{t}:G\to G$ (or $\eta_{\sigma}^{t}:X_{\Gamma}\to X_{\Gamma}$) by $\eta_{\sigma}^{t}x = \Phi_{\sigma,x}^{-1}(t)$, where $\Phi_{\sigma,x}$ is defined in Lemma \ref{L:FrankManningCoordinates2}.
\end{definition}
\begin{remark}
By Lemma \ref{L:FrankManningCoordinates2} the map $\eta_{\sigma}^{t}$ is well-defined.
\end{remark}
\begin{lemma}\label{L:BasicPropertiesLiftedTranslations}
The map $\eta_{\sigma}:E_{0}^{\sigma}\times G\to G$, $(t,x)\mapsto\eta_{\sigma}^{t}x$, is a Hölder $E_{0}^{\sigma}-$action that satisfies $F\eta_{\sigma}^{t}x = \eta_{\sigma}^{L_{su}t}Fx$. The action $\eta_{\sigma}^{t}$ naturally descends to $X_{\Gamma}$ and if $g\in Z^{\infty}(f)$, $B$ is the automorphism defined by $\Phi(gx) = B\Phi(x)$, then $g\eta_{\sigma}^{t}x = \eta_{\sigma}^{Bt}gx$. Finally, $\Phi(\eta_{\sigma}^{t}x) = t + \Phi(x)$, so $\Phi$ semi-conjugates $\eta_{\sigma}^{t}$ to the standard translation action along $E_{0}^{\sigma}$ on the base.
\end{lemma}
\begin{proof}
That $\eta_{\sigma}^{t}$ defines an action is immediate from the definition. Indeed, for $t,s\in E_{0}^{\sigma}$
\begin{align*}
\Phi(\eta_{\sigma}^{s}\eta_{\sigma}^{t}x) = & \Phi(\Phi_{\sigma,\eta_{\sigma}^{t}x}^{-1}(s)) = s + \Phi(\eta_{\sigma}^{t}x) = s + \Phi_{\sigma,x}(\eta_{\sigma}^{t}x) + \Phi(x) = \\ &
s + t + \Phi(x)
\end{align*}
so if we subtract $\Phi(x)$ then $\Phi_{\sigma,x}(\eta_{\sigma}^{s}\eta_{\sigma}^{t}x) = s + t$. Applying $\Phi_{\sigma,x}^{-1}$ on both sides of the equality yields $\eta_{\sigma}^{s}\eta_{\sigma}^{t}x = \eta_{\sigma}^{s+t}x$. This also shows that $\Phi(\eta_{u}^{t}x) = \Phi(x) + t$. Given $\gamma\in\Gamma$ we have $\Phi_{\sigma,\gamma x}(\gamma y) = \Phi(\gamma y) - \Phi(\gamma x) = \Phi(y) - \Phi(x) = \Phi_{\sigma,x}(y)$ which implies $\eta_{\sigma}^{t}(\gamma x) = \gamma\eta_{\sigma}^{t}x$, so $\eta_{\sigma}^{t}$ descend to $X_{\Gamma}$. For $g\in Z^{\infty}(f)$ we have $\Phi(gx) = B\Phi(x) + p_{0}$ for some $p_{0}\in\mathbb{T}^{d}$ (that is fixed by $L_{su}$). It follows that
\begin{align}
\Phi_{\sigma,gx}(gy) = \Phi(gy) - \Phi(gx) = B\left(\Phi(x) - \Phi(y)\right) = B\Phi_{\sigma,x}(y)
\end{align}
or with $\Phi_{\sigma,x}(y) = t$
\begin{align}
\eta_{\sigma}^{Bt}(gx) = gy = g\eta_{\sigma}^{t}(x).
\end{align}
Next, we show that $\eta_{\sigma}$ is Hölder. Let $u = \sigma$, the other case is similar. Since $\Phi_{u,x}$ is a bi-Hölder homeomorphism (Lemma \ref{L:FrankManningCoordinates2}) it is immediate that $(t,x)\mapsto\eta_{u}^{t}x$ is Hölder in $t$. The foliations $\Hat{W}^{u}$ and $\Hat{W}^{cs}$ are uniformly transverse, so we find $\varepsilon_{0} > 0$ and $K$ such that for $x,y\in X_{\Gamma}$ with $\intd(x,y)\leq\varepsilon_{0}$ we have $\Hat{W}_{\rm loc}^{u}(x)\cap\Hat{W}_{\rm loc}^{cs}(y) = \{z\}$,
\begin{align}\label{Eq:LocalProductStructureEstimate}
\intd_{u}(x,z)\leq K\intd(x,y),\quad\intd_{cs}(z,y)\leq K\intd(x,y).
\end{align}
If $y\in\Hat{W}^{cs}(x)$ then Lemma \ref{L:FrankManningCoordinates3} shows that $\Phi(y)\in\Phi(x) + E_{0}^{s}$, so
\begin{align*}
\Phi(\eta_{u}^{t}y) = & \Phi(y) + t = \Phi(x) + t + [\Phi(y) - \Phi(x)] = \\ & \Phi(\eta_{u}^{t}x) + [\Phi(y) - \Phi(x)]\in\Phi(\eta_{u}^{t}x) + E_{0}^{s},
\end{align*}
or $\eta_{u}^{t}y\in\Hat{W}^{cs}(\Phi(\eta_{u}^{t}x))$. That is, $\eta_{u}^{t}$ preserve the foliation $\Hat{W}^{cs}$. Given $t\in E_{0}^{u}$
\begin{align}
\intd(\eta_{u}^{t}x,\eta_{u}^{t}y)\leq\intd_{u}(\eta_{u}^{t}x,\eta_{u}^{t}z) + \intd_{cs}(\eta_{u}^{t}z,\eta_{u}^{t}y)
\end{align}
with $x,y$ and $z$ as in Equation \ref{Eq:LocalProductStructureEstimate}. So, it suffices to show that $\eta_{u}^{t}$ is Hölder along $\Hat{W}^{u}$ and $\Hat{W}^{cs}$. For $y\in\Hat{W}^{u}(x)$
\begin{align*}
\Phi_{u,x}(\eta_{u}^{t}y) = & \Phi(\eta_{u}^{t}y) - \Phi(x) = \Phi_{u,y}(\eta_{u}^{t}y) + \Phi(y) - \Phi(x) = \\ &
t + \Phi(y) - \Phi(x)
\end{align*}
or $\eta_{u}^{t}y = \Phi_{u,x}^{-1}(t + \Phi(y) - \Phi(x))$. Since $\Phi$ is Hölder and $\Phi_{u,x}^{-1}$ is uniformly Hölder
\begin{align}
\intd(\eta_{u}^{t}x,\eta_{u}^{t}y)\leq C\norm{\Phi(y) - \Phi(x)}^{\theta}\leq C'\intd(x,y)^{\theta'}
\end{align}
so $\eta_{u}^{t}$ is Hölder along $\Hat{W}^{u}$. Given $y\in\Hat{W}^{cs}(x)$ we have $\eta_{u}^{t}y\in\Hat{W}^{cs}(\eta_{u}^{t}x)$ since $\eta_{u}^{t}$ preserve $\Hat{W}^{cs}$. On the other hand, $\eta_{u}^{t}y\in\Hat{W}^{u}(y)$ by the definition of $\Phi_{u,y}$. So $\eta_{u}^{t}y = \pi_{x,\eta_{u}^{t}x}^{u}(y)$. The unstable Holonomy is (uniformly) Hölder \cite{PughShubWilkinson1997}, so $\eta_{u}^{t}$ is Hölder along $\Hat{W}^{cs}$.
\end{proof}
\begin{definition}\label{Def:ActionOfsuPathGroup}
We define an action of $\mathcal{P}$ (or $\mathcal{P}^{c},\mathcal{P}_{[\chi]},\mathcal{Q}_{[\chi]},\mathcal{P}_{[\chi]}^{c},\mathcal{Q}_{[\chi]}^{c}$) on $X_{\Gamma}$ (and on $G$) by
\begin{align}
\eta_{\mathcal{P}}^{w}x = \eta_{\mathcal{P}}^{w_{1}^{s}w_{1}^{u}...w_{N}^{s}w_{N}^{u}}x = \eta_{s}^{w_{1}^{s}}\eta_{u}^{w_{1}^{u}}...\eta_{s}^{w_{N}^{s}}\eta_{u}^{w_{N}^{u}}x
\end{align}
where $\eta_{\sigma}$ is defined in Definition \ref{Def:LiftedTranslations} and shown to be an action in Lemma \ref{L:BasicPropertiesLiftedTranslations}.
\end{definition}
\begin{remark}
We make no notational distinction between the action on $X_{\Gamma}$ and $G$. It is clear that the action on $G$ covers the action on $X_{\Gamma}$ in the sense that the projection $p_{\Gamma}:G\to X_{\Gamma}$ intertwines the two actions.
\end{remark}
\begin{lemma}\label{L:PropertiesSUaction1}
We have $\Phi(\eta_{\mathcal{P}}^{w}x) = \Phi(x) + \Pi(w)$. That is $\Phi$ semiconjugates the $\mathcal{P}-$action onto the translation action on $\mathbb{T}^{d}$ (and $\mathbb{R}^{d}$).
\end{lemma}
\begin{proof}
By induction, it suffices to consider $v\in E_{0}^{\sigma}$. The Lemma follows from Lemma \ref{L:BasicPropertiesLiftedTranslations}.
\end{proof}
\begin{lemma}\label{L:PriopertiesSUactionHolonomy}
For any $w\in\mathcal{P}$ the homeomorphism $\eta_{\mathcal{P}}^{w}$ preserve the center foliation, $W^{c}$. Moreover, if $w = v_{N}^{s}v_{N}^{u}...v_{1}^{s}v_{1}^{u}$ with $v_{j}^{\sigma}\in E_{0}^{\sigma}$, $x_{0} = x$ and
\begin{align*}
& x_{1}^{u} = \eta_{u}^{v_{1}^{u}}x_{0},\quad x_{1}^{s} = \eta_{s}^{v_{1}^{s}}x_{1}^{u}, \\
\vdots \\
& x_{N}^{u} = \eta_{u}^{v_{N}^{u}}x_{N-1}^{s},\quad x_{N}^{s} = \eta_{s}^{v_{N}^{s}}x_{N}^{u}
\end{align*}
then $\eta_{\mathcal{P}}^{w}x = x_{N}^{s}$ and the map $\eta_{\mathcal{P}}^{w}:W^{c}(x)\to W^{c}(\eta_{\mathcal{P}}^{w}x)$ coincide with the composition
\begin{align}
W^{c}(x_{0})\xrightarrow{\pi_{x,x_{1}^{u}}^{u}}W^{c}(x_{1}^{u})\xrightarrow{\pi_{x_{1}^{u},x_{1}^{s}}^{s}}W^{c}(x_{1}^{s})\xrightarrow{\pi_{x_{1}^{s},x_{2}^{u}}^{u}}...\xrightarrow{\pi^{s}_{x_{N}^{u},x_{N}^{s}}}W^{c}(x_{N}^{s}).
\end{align}
That is, we have
\begin{align}
\eta_{\mathcal{P}}^{w}|_{W^{c}(x)} = \pi_{x_{N}^{u},x_{N}^{s}}^{s}\circ\pi_{x_{N-1}^{s},x_{N}^{u}}^{u}\circ\pi_{x_{N-1}^{s},x_{N-1}^{u}}^{s}\circ...\circ\pi_{x_{1}^{u},x_{1}^{s}}^{s}\circ\pi_{x_{0},x_{1}^{u}}^{u}.
\end{align}
If $f$ is $r-$bunching then holonomies between center manifolds are $C^{r}-$smooth, so $\eta_{\mathcal{P}}^{w}$ is $C^{r}$ along $W^{c}$ for all $w\in\mathcal{P}$.
\end{lemma}
\begin{proof}
By induction it suffices to consider $w = v\in E_{0}^{\sigma}$ for $\sigma = s,u$. This was shown in the proof of Lemma \ref{L:BasicPropertiesLiftedTranslations}. The regularity follows from \cite{PughShubWilkinson1997}.
\end{proof}
\begin{lemma}\label{L:PropertiesSUaction2}
We have $w\in\mathcal{P}^{c}$ if and only if $\eta_{\mathcal{P}}^{w}x\in\Hat{W}^{c}(x)$ for every $x\in G$.
\end{lemma}
\begin{remark}
When $w\in\mathcal{P}^{c}$ then Lemma \ref{L:PropertiesSUaction2} shows that $\eta_{\mathcal{P}}^{w}W^{c}(x) = W^{c}(x)$ for every $x$, so $\mathcal{P}^{c}$ is the (homotopically trivial) \textbf{center fixing} part of $\mathcal{P}$.
\end{remark}
\begin{proof}
The lemma is immediate from $\Phi(\eta_{\mathcal{P}}^{w}x) = \Phi(x) + \Pi(w)$ (Lemma \ref{L:PropertiesSUaction1}).
\end{proof}
\begin{lemma}\label{L:PropertiesSUaction3}
If $f$ is accessible, then the $\mathcal{P}^{c}-$action is transitive on $\Hat{W}^{c}(x)$ (and $W^{c}(x)$) for all $x\in G$ (and $x\in X_{\Gamma}$).
\end{lemma}
\begin{proof}
Since $f$ is accessible $\eta_{\mathcal{P}}^{\mathcal{P}}x\supset W^{c}(x)$ for every $x\in X_{\Gamma}$, and $\eta_{\mathcal{P}}^{w}x\in W^{c}(x)$ if and only if $\Pi(w)\in\mathbb{Z}^{d}$ (Lemma \ref{L:PropertiesSUaction1}). Given $\mathbf{n}\in\mathbb{Z}^{d}$ write $\mathcal{P}_{\mathbf{n}} = \{w\in\mathcal{P}\text{ : }\Pi(w) = \mathbf{n}\}$. Then $\mathcal{P}^{c} = \mathcal{P}_{0}$. It follows that
\begin{align*}
W^{c}(x) = \bigcup_{\mathbf{n}\in\mathbb{Z}^{d}}\eta_{\mathcal{P}}^{\mathcal{P}_{\mathbf{n}}}x
\end{align*}
and since $W^{c}(x)$ is uncountable and $\mathbb{Z}^{d}$ is countable there is at least one $\mathbf{n}_{0}$ such that $\#\eta_{\mathcal{P}}(\mathcal{P}_{\mathbf{n}_{0}})x > 1$. If we fix some $w\in\mathcal{P}_{-\mathbf{n}_{0}}$ then 
\begin{align*}
\#\eta_{\mathcal{P}}^{w}\eta_{\mathcal{P}}^{\mathcal{P}_{\mathbf{n}_{0}}}x = \#\eta_{\mathcal{P}}^{w\mathcal{P}_{\mathbf{n}_{0}}}x > 1.
\end{align*}
For any $w'\in\mathcal{P}_{\mathbf{n}_{0}}$ we have $\Pi(ww') = \Pi(w) + \Pi(w') = 0$ so $w\mathcal{P}_{\mathbf{n}_{0}}\subset\mathcal{P}^{c}$, or $\#\eta_{\mathcal{P}}^{\mathcal{P}^{c}}x > 1$. Given any $w = w_{1}^{s}w_{1}^{u}...w_{N}^{s}w_{N}^{u}\in\mathcal{P}^{c}$ we define a path $w_{t} = (tw_{1}^{s})(tw_{1}^{u})...(tw_{N}^{s})(tw_{N}^{u})\in\mathcal{P}^{c}$, $t\in[0,1]$, from $0$ to $w$, so $\eta_{\mathcal{P}}^{\mathcal{P}^{c}}x$ is path connected. The image $I = \eta_{\mathcal{P}}^{\mathcal{P}^{c}}x$ is an interval in $W^{c}(x)$ since it contains at least $2$ distinct points and is path connected. We claim that $x$ is an interior point in this interval. If $x$ is not an interior point, then $I = [x,y),[x,y]$ or $(y,x],[y,x]$ for some $y\in W^{c}(x)$. We will assume that one of the first two cases holds the other two cases are similar. Let $x\neq z\in I$ and let $\eta_{\mathcal{P}}^{w}x = z$, then $[x,z]\subset I$. Since $w^{-1}\in\mathcal{P}^{c}$ and $\eta_{\mathcal{P}}^{w^{-1}}$ preserves orientation we have 
\begin{align}
I\supset\eta_{\mathcal{P}}^{w^{-1}}[x,z] = [\eta_{\mathcal{P}}^{w^{-1}}x,x]
\end{align}
which would imply $\eta_{\mathcal{P}}^{w^{-1}}x = x$ if $x$ is an end point of $I$. After applying $\eta_{\mathcal{P}}^{w}$ we obtain $x = \eta_{\mathcal{P}}^{w}x = z$. This is a contradiction since we assumed $z\neq x$. The point $x$ is interior in $\eta_{\mathcal{P}}^{\mathcal{P}^{c}}x$ and $x$ was arbitrary, so the orbit of $x$ under $\eta_{\mathcal{P}}^{\mathcal{P}_{c}}$ is open. This holds for every $x$ and $W^{c}(x)$ is connected, so $\eta_{\mathcal{P}}^{\mathcal{P}^{c}}x = W^{c}(x)$. The second part of the lemma follows.

The first claim follows from the second part. Indeed the second part implies that for any $x\in G$ the $\mathcal{P}^{c}-$orbit of $x$ is open in $\Hat{W}^{c}(x)$. Connectedness of $\Hat{W}^{c}(x)$ implies the first part of the lemma.
\end{proof}
Recall that $\rho:\mathbb{Z}^{k}\to{\rm GL}(d,\mathbb{Z})$ is defined by $\Phi\alpha^{\mathbf{n}} = \rho^{\mathbf{n}}\Phi$. Let $\rho^{\mathbf{n}}:\mathcal{P}\to\mathcal{P}$ be the map 
\begin{align*}
\rho^{\mathbf{n}}(w_{1}^{s}w_{1}^{u}...w_{N}^{s}w_{N}^{u}) = (\rho^{\mathbf{n}}w_{1}^{s})(\rho^{\mathbf{n}}w_{1}^{u})...(\rho^{\mathbf{n}}w_{N}^{s})(\rho^{\mathbf{n}}w_{N}^{u})
\end{align*}
The following lemma is immediate from Lemma \ref{L:BasicPropertiesLiftedTranslations}.
\begin{lemma}\label{L:SolvableRelations}
We have $\alpha^{\mathbf{n}}\eta_{\mathcal{P}}^{w} = \eta_{\mathcal{P}}^{\rho^{\mathbf{n}}w}\alpha^{\mathbf{n}}$ and $\rho^{\mathbf{n}}$ preserve $\mathcal{P}^{c}$, $\mathcal{P}_{[\chi]}$, $\mathcal{Q}_{[\chi]}$, $\mathcal{P}_{[\chi]}^{c}$, and $\mathcal{Q}_{[\chi]}^{c}$.
\end{lemma}

\section{An invariance principle for higher rank Anosov actions}\label{Sec:HigherRankInvariancePrinciple}

Let $X_{\Gamma}$ be any nilmanifold and $\rho:\mathbb{Z}^{k}\to{\rm Aut}(X_{\Gamma})$ a higher rank action. We will assume that $\rho$ is the restriction of some map $Q:\mathbb{R}^{k}\to{\rm Aut}(G)$. Let $\Phi_{\mathcal{X}}:\mathcal{X}\to X_{\Gamma}$ be a Hölder fiber bundle over $X_{\Gamma}$ with fibers $\mathcal{X}_{x} = \Phi_{\mathcal{X}}^{-1}(x)$ uniformly $C^{r}$ for some $r > 1$ (we allow $r\in(1,2)$). We will assume throughout this section that $\mathcal{X}$ is compact, and therefore have compact fibers.
\begin{definition}
We say that $F:\mathbb{Z}^{k}\times\mathcal{X}\to\mathcal{X}$ is a cocycle over $\rho$ if $F$ is a $\mathbb{Z}^{k}-$action covering $\rho$. Moreover, $F$ is a $C^{s}-$cocycle if $F(\mathbf{n},\cdot):\mathcal{X}_{x}\to\mathcal{X}_{\rho^{\mathbf{n}}x}$ is uniformly $C^{s}$.
\end{definition}
For a cocycle $F$ over $\rho$ we write $F^{\mathbf{n}}:\mathcal{X}_{x}\to\mathcal{X}_{\rho^{\mathbf{n}}x}$ for the map $F(\mathbf{n},\cdot)$. For each coarse exponent $[\chi]$ of $\rho$ we have a translation action $T_{[\chi]}:G^{[\chi]}\times X_{\Gamma}\to X_{\Gamma}$, $T_{[\chi]}^{g}x = xg^{-1}$ for $g\in G^{[\chi]}$.
\begin{definition}
A cocycle $F$ has $[\chi]-$translations if there is a Hölder action $\eta_{[\chi]}:G^{[\chi]}\times\mathcal{X}\to\mathcal{X}$ covering $T_{[\chi]}$ such that
\begin{align}
F^{\mathbf{n}}\eta_{[\chi]}^{g} = \eta_{[\chi]}^{\rho^{\mathbf{n}}g}F^{\mathbf{n}},
\end{align}
for any $\mathbf{n}\in\mathbb{Z}^{k}$.
\end{definition}
Our interest in cocycles over algebraic actions comes from the following lemma.
\begin{lemma}\label{L:HigherRankActionIsCocycle}
If $\alpha$ is as in Theorem \ref{Thm:MainTheorem1}, with $\Phi\alpha^{\mathbf{n}} = \rho^{\mathbf{n}}\Phi$, then $\alpha$ is a $C^{1+\alpha}-$cocycle over $\rho$ that admit $[\chi]-$translations for every coarse $[\chi]$ and $\eta_{[\chi]} = \eta_{\mathcal{P}}|_{E_{0}^{[\chi]}}$.
\end{lemma}
\begin{proof}
The lemma is immediate from Theorem \ref{Thm:PropertiesOfPHdiffeos} and Lemma \ref{L:PropertiesSUaction1}.
\end{proof}
The main result of this section is a sufficient condition for the translation action $\eta_{[\chi]}$ to preserve a $F-$invariant measure.
\begin{theorem}\label{Thm:InvariancePrincipleForHiherRankAction}
Let $F$ be a $C^{s}-$cocycle, $s > 1$, over $\rho$. Let $\nu$ be a $F-$invariant probability measure projecting onto $\mu_{\Gamma}$, $(\Phi_{\mathcal{X}})_{*}\nu = \mu_{\Gamma}$, and $\lambda_{F,\nu}^{1},...,\lambda_{F,\nu}^{N}:\mathbb{Z}^{k}\to\mathbb{R}$ the $\nu-$Lyapunov exponents of $F$ along the fibers of $\mathcal{X}$. If
\begin{align}
\bigcap_{i = 1}^{N}\ker\lambda_{F,\nu}^{i}\not\subset\ker[\chi]
\end{align}
then $\nu$ is $\eta_{[\chi]}-$invariant.
\end{theorem}
\begin{remark}
In Theorem \ref{Thm:InvariancePrincipleForHiherRankAction} we assume that $\rho$ is the restriction of some homomorphism $Q:\mathbb{R}^{k}\to{\rm Aut}(G)$. This is without loss of generality after possibly dropping to a finite index subgroup.
\end{remark}
To apply results from \cite{AvilaViana2010} it will be convenient to reformulate $\eta_{[\chi]}-$invariance of $\nu$ into essential holonomy invariance. Let $x\in X_{\Gamma}$, $g\in G^{[\chi]}$ and $y = T_{[\chi]}^{g}x$. Since $\eta_{[\chi]}^{g}$ cover $T_{[\chi]}^{g}$, we define
\begin{align}\label{Eq:NonSuspendedHolonomy}
h_{x,y}^{[\chi]}:\mathcal{X}_{x}\to\mathcal{X}_{y},\quad h_{x,y}^{[\chi]}(\xi) = \eta_{[\chi]}^{g}(\xi).
\end{align}
We say that $h_{x,y}^{[\chi]}$ is the $[\chi]-$holonomy between $\mathcal{X}_{x}$ and $\mathcal{X}_{y}$.
\begin{definition}
Let $\nu$ be $F-$invariant such that $(\Phi_{\mathcal{X}})_{*}\nu = \mu_{\Gamma}$ and let $\{\nu_{x}\}_{x\in X_{\Gamma}}$ be the disintegration of $\nu$ over $\Phi_{\mathcal{X}}$. We say that $\nu$, or $\{\nu_{x}\}_{x\in X_{\Gamma}}$, is essentially $[\chi]-$holonomy invariant if there is a $\mu_{\Gamma}-$full measure set $Y\subset X_{\Gamma}$ such that $(h_{x,y}^{[\chi]})_{*}\nu_{x} = \nu_{y}$ for $x,y\in Y$.
\end{definition}
\begin{lemma}\label{L:HolonomyInvarianceGiveInvariance}
Let $\nu$ be $F-$invariant and projecting onto $\mu_{\Gamma}$. If $\nu$ is essentially $[\chi]-$holonomy invariant then $\nu$ is $\eta_{[\chi]}-$invariant.
\end{lemma}
\begin{proof}
Let $Y\subset X_{\Gamma}$ be a full measure subset such that $(h_{x,y}^{[\chi]})_{*}\nu_{x} = \nu_{y}$ for $x,y\in Y$. Let $g\in G^{[\chi]}$ and $\Tilde{Y} = Y\cap T_{[\chi]}^{g^{-1}}Y$ so that $x,T_{[\chi]}^{g}x\in Y$ for $x\in\Tilde{Y}$. If $\varphi\in C^{0}(\mathcal{X})$ then
\begin{align*}
\int_{\mathcal{X}}\varphi(\eta_{[\chi]}^{g}\xi)\intd\nu(\xi) = & \int_{\Tilde{Y}}\left(\int_{\mathcal{X}_{x}}\varphi(\eta_{[\chi]}^{g}\xi)\intd\nu_{x}(\xi)\right)\intd\mu_{\Gamma}(x) = \\ &
\int_{\Tilde{Y}}\left(\int_{\eta_{[\chi]}^{g}\mathcal{X}_{x}}\varphi(\xi)\intd(h_{x,T_{[\chi]}^{g}x}^{[\chi]})_{*}\nu_{x}(\xi)\right)\intd\mu_{\Gamma}(x) = \\ &
\int_{\Tilde{Y}}\left(\int_{\mathcal{X}_{T_{[\chi]}^{g}x}}\varphi(\xi)\intd\nu_{T_{[\chi]}^{g}x}(\xi)\right)\intd\mu_{\Gamma}(x) = \\ &
\int_{\mathcal{X}}\varphi(\xi)\intd\nu(\xi),
\end{align*}
so $\nu$ is $\eta_{[\chi]}^{g}-$invariant.
\end{proof}

\subsection{The suspension construction}

Fix a higher rank action $\rho:\mathbb{Z}^{k}\to{\rm Aut}(X_{\Gamma})$, a cocycle $F^{\mathbf{n}}:\mathcal{X}\to\mathcal{X}$ over $\rho$, and a measure $\nu$ as in Theorem \ref{Thm:InvariancePrincipleForHiherRankAction}. We recall the definition of the suspension of an action $\alpha:\mathbb{Z}^{k}\times M\to M$.
\begin{definition}
Let $\tau:\mathbb{Z}^{k}\times(M\times\mathbb{R}^{k})\to M\times\mathbb{R}^{k}$ be defined by $\tau^{\mathbf{n}}(x,\mathbf{s}) = (\alpha^{\mathbf{n}}x,\mathbf{s}-\mathbf{n})$. We define the suspension $\mathcal{S}$ of $\alpha$ as
\begin{align*}
\mathcal{S} := \left(M\times\mathbb{R}^{k}\right)/\tau.
\end{align*}
Given $(x,\mathbf{s})\in M\times\mathbb{R}^{k}$ we denote by $[x,\mathbf{s}]$ the equivalence class of $(x,\mathbf{s})$ in $\mathcal{S}$. We also define a natural action on $\mathcal{S}$ by $\alpha_{\mathcal{S}}^{\mathbf{t}}(x,\mathbf{s}) = (x,\mathbf{s} + \mathbf{t})$. Since
\begin{align*}
\alpha_{\mathcal{S}}^{\mathbf{t}}\tau^{\mathbf{n}}(x,\mathbf{s}) = & \alpha_{\mathcal{S}}^{\mathbf{t}}(\alpha^{\mathbf{n}}x,\mathbf{s}-\mathbf{n}) = (\alpha^{\mathbf{n}}x,\mathbf{t}+\mathbf{s}-\mathbf{n}) = \tau^{\mathbf{n}}(x,\mathbf{t}+\mathbf{s}) = \\ &
\tau^{\mathbf{n}}\alpha_{\mathcal{S}}^{\mathbf{t}}(x,\mathbf{s})
\end{align*}
the action $\alpha_{\mathcal{S}}$ descends to an action on $\mathcal{S}$. Moreover the map $M\times\mathbb{R}^{k}\ni(x,\mathbf{s})\mapsto\mathbf{s} + \mathbb{Z}^{k}\in\mathbb{T}^{k}$ descends to a map $\pi_{\mathcal{S}}:\mathcal{S}\to\mathbb{T}^{k}$ with fibers $M$. The map $\pi_{\mathcal{S}}$ semi-conjugates $\alpha_{\mathcal{S}}$ to the natural translation action on $\mathbb{T}^{k}$.
\end{definition}
Given any $\alpha-$invariant measure $\mu$ on $M$ we define a measure $\mu_{\mathcal{S}}$ on $\mathcal{S}$ as follows. For each $x\in\mathbb{T}^{k}$ we choose some $\mathbf{s}\in\mathbb{R}^{k}$ such that $\mathbf{s} + \mathbb{Z}^{k} = x$. Let $\iota_{\mathbf{s}}:M\to\pi_{\mathcal{S}}^{-1}(x)\subset\mathcal{S}$ be defined by $\iota_{\mathbf{s}}(y) = [y,\mathbf{s}]$. Define a measure $\mu_{x}$ on $\pi_{\mathcal{S}}^{-1}(x)$ by 
\begin{align}
(\iota_{\mathbf{s}})_{*}\mu = \mu_{x}. 
\end{align}
Given any $\mathbf{n}\in\mathbb{Z}^{k}$ we have $\iota_{\mathbf{s} + \mathbf{n}}(y) = [y,\mathbf{s}+\mathbf{n}] = [\alpha^{\mathbf{n}}y,\mathbf{s}] = \iota_{\mathbf{s}}\alpha^{\mathbf{n}}y$. So $\alpha-$invariance of $\mu$ implies
\begin{align*}
(\iota_{\mathbf{s} + \mathbf{n}})_{*}\mu = (\iota_{\mathbf{s}})_{*}\alpha^{\mathbf{n}}_{*}\mu = (\iota_{\mathbf{s}})_{*}\mu
\end{align*}
showing that $\mu_{x}$ is well-defined. Define a suspended measure $\mu_{\mathcal{S}}$ by
\begin{align}
\mu_{\mathcal{S}} = \int_{\mathbb{T}^{k}}\mu_{x}\intd{\rm vol}_{\mathbb{T}^{k}}(x).
\end{align}
One checks that $\mu_{\mathcal{S}}$ is $\alpha_{\mathcal{S}}-$invariant.

In the remainder of this section we denote by $\mathcal{S}_{0}$ the suspension of $\rho$ with action $\rho_{\mathcal{S}_{0}}:\mathbb{R}^{k}\times\mathcal{S}_{0}\to\mathcal{S}_{0}$ and by $\mathcal{S}$ the suspension of $F$ with action $F_{\mathcal{S}}:\mathbb{R}^{k}\times\mathcal{S}\to\mathcal{S}$. We also denote by $\mu_{\mathcal{S}_{0}}$ the suspension of $\mu_{\Gamma}$ and $\nu_{\mathcal{S}}$ the suspension of $\nu$. Note that $\mu_{\mathcal{S}_{0}}$ is a volume on $\mathcal{S}_{0}$. Let $Q:\mathbb{R}^{k}\to{\rm Aut}(G)$ be a homomorphism such that
\begin{align}
Q|_{\mathbb{Z}^{k}} = \rho. 
\end{align}
We suspend the actions $T_{[\chi]}^{g}$ and $\eta_{[\chi]}^{g}$ as
\begin{align}
\Hat{T}_{[\chi]}^{g}\left([x,\mathbf{s}]\right) = \left[T_{[\chi]}^{Q^{-\mathbf{s}}g}x,\mathbf{s}\right],\quad\Hat{\eta}_{[\chi]}^{g}\left([\xi,\mathbf{s}]\right) = \left[\eta_{[\chi]}^{Q^{-\mathbf{s}}g}\xi,\mathbf{s}\right].
\end{align}
For any $\mathbf{n}\in\mathbb{Z}^{k}$ we have
\begin{align}
\Hat{T}_{[\chi]}^{g}([\rho^{\mathbf{n}}x,\mathbf{s} - \mathbf{n}]) = \left[T_{[\chi]}^{Q^{\mathbf{n}}Q^{-\mathbf{s}}g}\rho^{\mathbf{n}}x,\mathbf{s} - \mathbf{n}\right] = \left[\rho^{\mathbf{n}}T_{[\chi]}^{Q^{-\mathbf{s}}g}x,\mathbf{s}-\mathbf{n}\right]
\end{align}
so $\Hat{T}_{[\chi]}$ is a well-defined action on $\mathcal{S}_{0}$ that acts in the fibers of $\pi_{\mathcal{S}_{0}}:\mathcal{S}_{0}\to\mathbb{T}^{k}$. Similarly the action $\Hat{\eta}_{[\chi]}$ is well-defined on $\mathcal{S}$. Define $\Phi_{\mathcal{S}}:\mathcal{S}\to\mathcal{S}_{0}$ by
\begin{align}
\Phi_{\mathcal{S}}([\xi,\mathbf{s}]) = [\Phi_{\mathcal{X}}(\xi),\mathbf{s}].
\end{align}
Since $\Phi_{\mathcal{X}}$ semi-conjugates $\rho$ onto $F$, $\Phi_{\mathcal{S}}$ is well-defined. The following lemma is immediate from our definitions.
\begin{lemma}
Let $\Phi_{\mathcal{S}}$, $F_{\mathcal{S}}$, $\rho_{\mathcal{S}_{0}}$, $\Hat{T}_{[\chi]}$ and $\Hat{\eta}_{[\chi]}$ be as above. The following holds
\begin{enumerate}[label = (\roman*)]
    \item the map $\Phi_{\mathcal{S}}:\mathcal{S}\to\mathcal{S}_{0}$ is a Hölder fiber bundle with uniformly $C^{r}$ fibers. In fact, for any $\mathbf{s}\in\mathbb{T}^{k}$ the restriction $\Phi_{\mathcal{S}}|_{\pi_{\mathcal{S}}^{-1}(\mathbf{s})}$ coincides with $\Phi_{\mathcal{X}}$ using natural identifications of $\pi_{\mathcal{S}}^{-1}(\mathbf{s})\cong\mathcal{X}$ and $\pi_{\mathcal{S}_{0}}^{-1}(\mathbf{s})\cong X_{\Gamma}$,
    \item the map $F_{\mathcal{S}}$ is a $C^{s}$ cocycle over $\rho_{\mathcal{S}_{0}}$,
    \item the Lyapunov exponents for $\rho_{\mathcal{S}_{0}}$ coincide with the Lyapunov exponents of $\rho$,
    \item the $\nu_{\mathcal{S}}-$Lyapunov exponents along the fibers of $\Phi_{\mathcal{S}}$ for $F_{\mathcal{S}}$ coincide with the fiberwise Lyapunov exponents of $F$,
    \item the map $\Phi_{\mathcal{S}}$ semiconjugates $\Hat{\eta}_{[\chi]}$ to $\Hat{T}_{[\chi]}$.
\end{enumerate}
\end{lemma}
\begin{proof}
Point $(i)$ follows from the analogous properties of $\Phi_{\mathcal{X}}$ since $\Phi_{\mathcal{S}}$ is defined in the fibers of $\pi_{\mathcal{S}}:\mathcal{S}\to\mathbb{T}^{k}$. That $\Phi_{\mathcal{S}}$ conjugates $F_{\mathcal{S}}$ to $\rho_{\mathcal{S}_{0}}$ is immediate from its definition: $\Phi_{\mathcal{S}}[\xi,\mathbf{s}+\mathbf{t}] = [\Phi_{\mathcal{X}}(\xi),\mathbf{s} + \mathbf{t}]$. That $F_{\mathcal{S}}$ is $C^{s}$ along the fibers of $\Phi_{\mathcal{S}}$ is immediate since $F$ is $C^{s}$ along the fibers of $\Phi_{\mathcal{X}}$ (note that the identifications $\iota_{\mathbf{s}}(\xi) = [\xi,\mathbf{s}]$, $\mathcal{X}\to\pi_{\mathcal{S}}^{-1}(\mathbf{s} + \mathbb{Z}^{k})$, defines a smooth structure on the fibers of $\Phi_{\mathcal{S}}$ in which $F_{\mathcal{S}}$ is uniformly $C^{s}$). Points $(iii)$ and $(iv)$ holds for $\mathbb{Z}^{k}\subset\mathbb{R}^{k}$, and any functional is determined by its values on a lattice, proving $(iii)$ and $(iv)$. Point $(v)$ is immediate from the definitions and the fact that $\Phi_{\mathcal{X}}$ conjugate $\eta_{[\chi]}$ to $T_{[\chi]}$.
\end{proof}
Define holonomies along the orbits of $\Hat{T}_{[\chi]}$ as in Equation \ref{Eq:NonSuspendedHolonomy}. That is, if $x\in\mathcal{S}_{0}$, $y = \Hat{T}_{[\chi]}^{g}x$ then
\begin{align}
\Hat{h}_{x,y}^{[\chi]}(\xi) = \Hat{\eta}_{[\chi]}^{g}(\xi).
\end{align}
We say that $\nu_{\mathcal{S}}$ (or the disintegration $\{\nu_{\mathcal{S},x}\}_{x\in\mathcal{S}_{0}}$) is essentially $[\chi]-$holonomy invariant if there is a $\mu_{\mathcal{S}_{0}}$ full measure set $Y\subset\mathcal{S}_{0}$ such that
\begin{align}
(\Hat{h}_{x,y}^{[\chi]})_{*}\nu_{\mathcal{S},x} = \nu_{\mathcal{S},y},\quad x,y\in Y.
\end{align}
The key fact about the suspension, is that holonomy invariance of $\nu_{\mathcal{S}}$ implies holonomy invariance of $\nu$. So Theorem \ref{Thm:InvariancePrincipleForHiherRankAction} follows by $[\chi]-$holonomy invariance of $\nu_{\mathcal{S}}$ (by Lemma \ref{L:HolonomyInvarianceGiveInvariance}).
\begin{lemma}\label{L:HolonmyInvSuspGiveHolonomyInv}
The measure $\nu$ is essentially $[\chi]-$holonomy invariant if and only if $\nu_{\mathcal{S}}$ is essentially $[\chi]-$holonomy invariant.
\end{lemma}
Lemma \ref{L:HolonmyInvSuspGiveHolonomyInv} is immediate from the following lemma.
\begin{lemma}\label{L:FormulaForSuspendedDisintegration}
We have $(\Phi_{\mathcal{S}})_{*}\nu_{\mathcal{S}} = \mu_{\mathcal{S}_{0}}$. The disintegration of $\nu_{\mathcal{S}}$ over $\mu_{\mathcal{S}_{0}}$ is given by
\begin{align}
\nu_{\mathcal{S},[x,\mathbf{s}]} = (\iota_{\mathbf{s}})_{*}\nu_{x},\quad[x,\mathbf{s}]\in\mathcal{S}_{0}
\end{align}
where $\nu$ is the disintegration of $\nu$ over $\mu$.
\end{lemma}
\begin{proof}
By construction we have $\nu_{\mathcal{S}} = (\iota_{\mathbf{s}})_{*}\nu\otimes\intd\mathbf{s}$, so
\begin{align*}
(\Phi_{\mathcal{S}})_{*}\nu_{\mathcal{S}} = & \int_{\mathbb{T}^{k}}(\Phi_{\mathcal{S}})_{*}(\iota_{\mathbf{s}})_{*}\nu\intd\mathbf{s} = \int_{\mathbb{T}^{k}}(\Phi_{\mathcal{S}}\iota_{\mathbf{s}})_{*}\nu\intd\mathbf{s} = \int_{\mathbb{T}^{k}}(\iota_{\mathbf{s}}\Phi_{\mathcal{X}})_{*}\nu\intd\mathbf{s} = \\ &
\int_{\mathbb{T}^{k}}(\iota_{\mathbf{s}})_{*}\mu\intd\mathbf{s} = \mu_{\mathcal{S}_{0}}.
\end{align*}
If we define $\nu_{\mathcal{S},[x,\mathbf{s}]} := (\iota_{\mathbf{s}})_{*}\nu_{x}$ then
\begin{align}
\nu_{\mathcal{S},[\rho^{\mathbf{n}}x,\mathbf{s} - \mathbf{n}]} = (\iota_{\mathbf{s}-\mathbf{n}})_{*}\nu_{\rho^{\mathbf{n}}x} = (\iota_{\mathbf{s}}F^{\mathbf{-n}})_{*}\nu_{\rho^{\mathbf{n}}x} = (\iota_{\mathbf{s}})_{*}\nu_{x} = \nu_{\mathcal{S},[x,\mathbf{s}]}
\end{align}
where we have used $\iota_{\mathbf{s}+\mathbf{n}}(\xi) = [\xi,\mathbf{s} + \mathbf{n}] = [F^{-\mathbf{n}}F^{\mathbf{n}}\xi,\mathbf{s} + \mathbf{n}] = [\rho^{\mathbf{n}}\xi,\mathbf{s}] = \iota_{\mathbf{s}}(\rho^{\mathbf{n}}\xi)$ for $\xi\in\mathcal{X}$. So $\nu_{\mathcal{S},[x,\mathbf{s}]}$ is well-defined. We calculate
\begin{align*}
\int_{\mathcal{S}_{0}}\nu_{\mathcal{S},[x,\mathbf{s}]}\intd\mu_{\mathcal{S}_{0}}([x,\mathbf{s}]) = & \int_{\mathcal{S}_{0}}(\iota_{\mathbf{s}})_{*}\nu_{x}\intd\mu_{\mathcal{S}_{0}}([x,\mathbf{s}]) = \\ & \int_{\mathbb{T}^{k}}\left[\int_{\iota_{\mathbf{s}}X_{\Gamma}}(\iota_{\mathbf{s}})_{*}\nu_{x}\intd(\iota_{\mathbf{s}})_{*}\mu_{\Gamma}([x,\mathbf{s}])\right]\intd\mathbf{s} = \\ &
\int_{\mathbb{T}^{k}}\left[(\iota_{\mathbf{s}})_{*}\int_{X_{\Gamma}}\nu_{x}\intd\mu_{\Gamma}(x)\right]\intd\mathbf{s} = \\ &
\int_{\mathbb{T}^{k}}(\iota_{\mathbf{s}})_{*}\nu\intd\mathbf{s} = \nu_{\mathcal{S}}
\end{align*}
which proves that $\nu_{\mathcal{S},[x,\mathbf{s}]}$ is a disintegration of $\nu_{\mathcal{S}}$ over $\Phi_{\mathcal{S}}$.
\end{proof}
\begin{proof}[Proof of Lemma \ref{L:HolonmyInvSuspGiveHolonomyInv}]
If $\nu$ is essentially $[\chi]-$holonomy invariant, then we find $Y\subset X_{\Gamma}$ such that $(h_{x,y}^{[\chi]})_{*}\nu_{x} = \nu_{y}$ for $x,y\in Y$. Letting $\Tilde{Y}\subset\mathcal{S}_{0}$ be the image of $Y\times\mathbb{R}^{k}$ in $\mathcal{S}_{0}$, one direction in Lemma \ref{L:HolonmyInvSuspGiveHolonomyInv} follows from the formula in Lemma \ref{L:FormulaForSuspendedDisintegration}. For the converse direction, let $Y\subset\mathcal{S}_{0}$ be such that $(\Hat{h}_{x,y}^{[\chi]})_{*}\nu_{\mathcal{S},x} = \nu_{\mathcal{S},y}$ for $x,y\in Y$. Since $\mu_{\mathcal{S}_{0}}(Y) = 1$ we have
\begin{align*}
\mu_{\Gamma}(\iota_{\mathbf{s}}^{-1}(Y\cap\pi_{\mathcal{S}_{0}}^{-1}(\mathbf{s}))) = 1
\end{align*}
for $\intd\mathbf{s}-$almost every $\mathbf{s}\in\mathbb{T}^{k}$. For any $\mathbf{s}\in\mathbb{R}^{k}$, $g\in G^{[\chi]}$ and $x\in X_{\Gamma}$ we have
\begin{align*}
\iota_{\mathbf{s}}T_{[\chi]}^{g}(x) = \left[T_{[\chi]}^{g}(x),\mathbf{s}\right] = \left[T_{[\chi]}^{Q^{-\mathbf{s}}Q^{\mathbf{s}}g}(x),\mathbf{s}\right] = \Hat{T}_{[\chi]}^{Q^{\mathbf{s}}g}\left[x,\mathbf{s}\right] = \Hat{T}_{[\chi]}^{Q^{\mathbf{s}}g}\iota_{\mathbf{s}}(x).
\end{align*}
Similarly, $\iota_{\mathbf{s}}\eta_{[\chi]}^{g} = \Hat{\eta}_{[\chi]}^{Q^{\mathbf{s}}g}\iota_{\mathbf{s}}$. Since $\iota_{\mathbf{s}}$ maps fibers of $\Phi_{\mathcal{X}}$ to fibers of $\Phi_{\mathcal{S}}$, it follows that $\iota_{\mathbf{s}}^{-1}\Hat{h}_{[x,\mathbf{s}],[y,\mathbf{s}]}^{[\chi]}\iota_{\mathbf{s}} = h_{x,y}^{[\chi]}$. With $\Tilde{Y} = \iota_{\mathbf{s}}^{-1}(Y\cap\pi_{\mathcal{S}_{0}}^{-1}(\mathbf{s}))$ and Lemma \ref{L:FormulaForSuspendedDisintegration}
\begin{align*}
(h_{x,y}^{[\chi]})_{*}\nu_{x} = & \left(\iota_{\mathbf{s}}^{-1}\Hat{h}_{[x,\mathbf{s}],[y,\mathbf{s}]}^{[\chi]}\iota_{\mathbf{s}}\right)_{*}\nu_{x} = \left(\iota_{\mathbf{s}}^{-1}\Hat{h}_{[x,\mathbf{s}],[y,\mathbf{s}]}^{[\chi]}\right)_{*}\nu_{\mathcal{S},[x,\mathbf{s}]} = \\ &
(\iota_{\mathbf{s}}^{-1})_{*}\nu_{\mathcal{S},[y,\mathbf{s}]} = \nu_{y},
\end{align*}
for $x,y\in\Tilde{Y}$. Choosing $\mathbf{s}$ such that $\mu_{\Gamma}(\Tilde{Y}) = 1$, we see that $\nu$ is essentially $[\chi]-$holonomy invariant.
\end{proof}

\subsection{Proof of Theorem \ref{Thm:InvariancePrincipleForHiherRankAction}}

By Lemmas \ref{L:HolonomyInvarianceGiveInvariance} and \ref{L:HolonmyInvSuspGiveHolonomyInv} it suffices to show that the disintegration of the suspension of $\nu$ is essentially $[\chi]-$holonomy invariant. We will use the following general criteria for obtaining holonomy invariance, proved in \cite[Proposition 4.2]{AvilaViana2010} (or \cite[Corollary 4.3]{AvilaViana2010}).
\begin{theorem}\label{Thm:InvariancePrincipleHelpVersion}
Let $f:M\to M$ be a volume preserving diffeomorphism on a closed, smooth manifold with an invariant contracting smooth foliation $W$. Let $\mathcal{X}\to M$ be a Hölder fiber bundle with uniformly $C^{r}$ fibers and $F:\mathcal{X}\to\mathcal{X}$ a map covering $f$ such that $F$ is uniformly $C^{s}$ along the fibers of $\mathcal{X}\to M$. Assume that $W$ admits holonomies in $\mathcal{X}$, that is for every $y\in W(x)$ there is a map $h_{x,y}^{W}:\mathcal{X}_{x}\to\mathcal{X}_{y}$ satisfying $(sh1)$, $(sh2)$ and $(sh3)$ in \cite[Section 2.4]{AvilaViana2010}. Let $\nu$ be an $F-$invariant measure on $\mathcal{X}$ projecting onto volume. If the $\nu-$exponents of $F$ along the fibers of $\mathcal{X}\to M$ are $0$ then the disintegration of $\nu$ is essentially $W-$holonomy invariant. That is, there is a full volume set $Y\subset M$ such that $(h_{x,y}^{W})_{*}\nu_{x} = \nu_{y}$ for $x,y\in Y$.
\end{theorem}
\begin{proof}
The theorem would follow immediately from \cite[Proposition 4.2]{AvilaViana2010} if the foliation $W$ coincided with the stable foliation of $f$ (in the sense of \cite[Section 4.1]{AvilaViana2010}). However, following the proof, it suffices that $W$ is contracting. In fact, since $W$ is a contracting foliation it is standard to produce a measurable partition subordinate to $W$ (see for example \cite{LedrappierStrelcyn1982}) which simplifies the proof.
\end{proof}
\begin{proof}[Proof of Theorem \ref{Thm:InvariancePrincipleForHiherRankAction}]
By assumption we find $\mathbf{t}_{0}\in\mathbb{R}^{k}\setminus0$ such that $F_{\mathcal{S}}^{\mathbf{t}_{0}}$ has zero exponents along the fibers of $\Phi_{\mathcal{S}}$ and $[\chi](\mathbf{t}_{0}) < 0$. Let $f = \rho_{\mathcal{S}_{0}}^{\mathbf{t}_{0}}$, $F = F_{\mathcal{S}}^{\mathbf{t}_{0}}$, $W$ be the orbit foliation of $\Hat{T}_{[\chi]}$, and apply Theorem \ref{Thm:InvariancePrincipleHelpVersion} to conclude that $\nu_{\mathcal{S}}$ is essentially $[\chi]-$holonomy invariant. The theorem follows from Lemmas \ref{L:HolonomyInvarianceGiveInvariance} and \ref{L:HolonmyInvSuspGiveHolonomyInv}.
\end{proof}

\section{Invariant structure in the center direction}\label{Sec:InvariantStructureInCenter}

Let $\alpha:\mathbb{Z}^{k}\times X_{\Gamma}\to X_{\Gamma}$ be a smooth action satisfying the assumptions of Theorem \ref{Thm:MainTheorem1}. In this section we prove Theorem \ref{Thm:InvariantStructureInCenter}: $f$, and $\alpha$, have a unique measure of maximal entropy. Moreover, if $\mu$ is the measure of maximal entropy then $\Phi_{*}\mu = {\rm vol}$ and the disintegration of $\mu$ is invariant under stable and unstable holonomy. Equivalently \cite{AvilaViana2010,hertz_rodriguezhertz_tahzibi_ures_2012} we show that the $\mu-$center exponent vanish $\lambda_{\mu}^{c} = 0$. The proof of Theorem \ref{Thm:InvariantStructureInCenter} is by contradiction, so we assume that $\lambda_{\mu}^{c}(f)\neq 0$. The proof splits into two cases. First, we have a generic case when the kernel of $\lambda_{\mu}^{c}:\mathbb{Z}^{k}\to\mathbb{R}$ does not coincide with the kernel of some exponent $\chi:\mathbb{Z}^{k}\to\mathbb{R}$ of $\rho$. Second, we have an exceptional case when $\ker\lambda_{\mu}^{c} = \ker\chi$ for some exponent $\chi$ of $\rho$. The first case is dealt with by using Theorem \ref{Thm:InvariancePrincipleForHiherRankAction} and Lemma \ref{L:HigherRankActionIsCocycle}. The second, more technical, case is dealt with by studying the circle dynamics induced by the holonomy maps on the center leaves. Suppose that $\lambda_{\mu}^{c}$ has the same kernel as $[\chi]$. We begin by showing that the $\mathcal{P}_{[\chi]}-$action commute with the $\mathcal{Q}_{[\chi]}-$action. This implies that either $\mathcal{P}_{[\chi]}^{c}$ or $\mathcal{Q}_{[\chi]}^{c}$ act transitively on center leaves (see Lemma \ref{L:TransitiveActionCircle1}). If $\mathcal{Q}_{[\chi]}^{c}$ act transitively on center leaves then Theorem \ref{Thm:InvariancePrincipleForHiherRankAction} can be applied, as in the generic case. If $\mathcal{P}_{[\chi]}^{c}$ acts transitively on center leaves, then we show that $\mathcal{Q}_{[\chi]}$ acts minimally on $X_{\Gamma}$. We use the minimality of the $\mathcal{Q}_{[\chi]}-$action, and the fact that $\mathcal{Q}_{[\chi]}$ commute with $\mathcal{P}_{[\chi]}^{c}$ to produce a continuous $\mathbb{T}-$action preserving $W^{c}$ that commutes with $\alpha$. This shows that the exponent $\lambda_{\mu}^{c}$ must vanish, a contradiction.

Denote by $\mathcal{M}_{\rm vol}^{f}(X_{\Gamma})$ the $f-$invariant measures projecting to volume
\begin{align}
\Phi_{*}\mu = {\rm vol},\quad\mu\in\mathcal{M}^{f}(X_{\Gamma}).
\end{align}
Equivalently the measures $\mathcal{M}_{\rm vol}^{f}(X_{\Gamma})$ are precisely the measures of maximal entropy for $f$ \cite{hertz_rodriguezhertz_tahzibi_ures_2012}. From \cite{hertz_rodriguezhertz_tahzibi_ures_2012} the set $\mathcal{M}_{\rm vol}^{f}(X_{\Gamma})$ is finite so we may assume, after possibly dropping to a finite index subgroup, that $\mathcal{M}_{\rm vol}^{f}(X_{\Gamma})$ consist of $\alpha-$invariant measures. Assume for contradiction that $\lambda_{\mu}^{c}(f)\neq 0$ for some $\mu\in\mathcal{M}_{\rm vol}^{f}(X_{\Gamma})$. This implies $\lambda_{\nu}^{c}(f)\neq0$ for all $\nu\in\mathcal{M}_{\rm vol}^{f}(X_{\Gamma})$ \cite{hertz_rodriguezhertz_tahzibi_ures_2012}.
\begin{lemma}\label{L:InvariantManifoldsMMEs}
For any two $\mu,\nu\in\mathcal{M}_{\rm vol}^{f}(X_{\Gamma})$ we have $\ker\lambda_{\mu}^{c} = \ker\lambda_{\nu}^{c}$ where $\lambda_{\mu}^{c},\lambda_{\nu}^{c}:\mathbb{Z}^{k}\to\mathbb{R}$. Moreover, there is a ${\rm vol}-$full measure set $Y\subset\mathbb{T}^{d}$ such that for any $y\in Y$ we have $x_{1},...,x_{N}\in\Phi^{-1}(y)$ with $x_{1} < x_{2} < ... < x_{N}$ in the orientation of $\Phi^{-1}(y)$ such that $(x_{j-1},x_{j+1})$ is the stable or unstable manifold for some $\nu'\in\mathcal{M}_{\rm vol}(X_{\Gamma})$ in $\Phi^{-1}(y)$.
\end{lemma}
\begin{proof}
We sketch the construction of measures in \cite{hertz_rodriguezhertz_tahzibi_ures_2012}. Let $\mu\in\mathcal{M}_{\rm vol}^{f}(X_{\Gamma})$ and $(\mu_{x}^{c})_{x\in X_{\Gamma}}$ be the disintegration of $\mu$ with respect to the center foliation. The measures $\mu_{x}^{c}$ are atomic $\mu-$almost everywhere since $\lambda_{\mu}^{c}(f)\neq0$ and $W^{c}$ has $1-$dimensional leaves. Denote by $\mathcal{F}_{\mu}$ the $f-$invariant foliation (contracting or expanding) manifolds in $W^{c}$ defined $\mu-$almost everywhere. Let $Y\subset\mathbb{T}^{d}$ be such that $\mu_{x}^{c}$ exists and is atomic for each $x\in\Phi^{-1}(y)$. For $y\in Y$ let 
\begin{align}
\mu_{x}^{c} = (\delta_{p_{1}(x)}+...+\delta_{p_{k}(x)})/k.
\end{align}
Let $q_{j}(x)$ be the positively oriented endpoint of $\mathcal{F}_{\mu}(p_{j}(x))$. Define a new measure $\nu$ by $\nu_{y} = \delta_{q_{1}(x)}+...+\delta_{q_{k}(x)}$ and $\nu = \nu_{y}\otimes\intd{\rm vol}(y)$. Then $\Phi_{*}\nu = {\rm vol}$ and $\nu$ is ergodic, so $\nu\in\mathcal{M}_{\rm vol}^{f}(X_{\Gamma})$. For any $\mathbf{n}\in\mathbb{Z}^{k}$ such that $\lambda_{\mu}^{c}(\mathbf{n}) < 0$, $\mathcal{F}_{\mu}(p_{j}(x))$ is a stable manifold for $\alpha^{\mathbf{n}}$. If $\lambda_{\nu}^{c}(\mathbf{n}) < 0$, then any point in $(p_{j}(x),q_{j}(x))$ lie in the stable manifold for both $p_{j}(x)$ and $q_{j}(x)$, which is a contradiction. It follows that $\lambda_{\mu}^{c}(\mathbf{n}) < 0$ implies $\lambda_{\nu}^{c}(\mathbf{n})\geq 0$, or $\lambda_{\mu}^{c} = -c\lambda_{\nu}^{c}$ for some $c > 0$ (note that $c\neq0$ since there are no measure $\nu\in\mathcal{M}_{\rm vol}^{f}(X_{\Gamma})$ with zero center exponent \cite{hertz_rodriguezhertz_tahzibi_ures_2012}). The first part of the lemma follows for measures constructed as above. Since $\#\mathcal{M}_{\rm vol}(X_{\Gamma}) < \infty$ \cite{hertz_rodriguezhertz_tahzibi_ures_2012} the construction of new invariant measures outlined above can only produce new measures finitely many times. This proves the last part of the lemma, since if the invariant manifolds did not cover the center leaves, then we could proceed the construction. This also proves the first part of the lemma since the measures constructed above have invariant manifolds covering the entire center leaves for ${\rm vol}-$almost every $\Phi^{-1}(y)$.
\end{proof}
 
\subsection{Generic case of Theorem \ref{Thm:InvariantStructureInCenter}}

Let $\mu\in\mathcal{M}_{\rm vol}^{f}(X_{\Gamma})$. If $\ker\lambda_{\mu}^{c}$ does not coincide with $\ker[\chi]$ for any coarse exponent of $\rho$, then $\eta_{\mathcal{P}}|_{E_{0}^{[\chi]}}$ preserve $\mu$ for all coarse $[\chi]$ by Theorem \ref{Thm:InvariancePrincipleForHiherRankAction} and Lemma \ref{L:HigherRankActionIsCocycle}. So, $\mu$ is $\eta_{\mathcal{P}}-$invariant. Accessibility of $f$ implies that $\eta_{\mathcal{P}}$ acts transitively, which implies that the disintegration of $\mu$, $\mu_{x}^{c}$, is not atomic. This is a contradiction since we assumed $\lambda_{\mu}^{c}(f)\neq0$.

\subsection{Exceptional case of Theorem \ref{Thm:InvariantStructureInCenter}}

Now we deal with the exceptional case of Theorem \ref{Thm:InvariantStructureInCenter}. Fix $\mu\in\mathcal{M}_{\rm vol}^{f}(X_{\Gamma})$ and assume that $\lambda_{\mu}^{c} = r\chi$ for some Lyapunov exponent $\chi$ of $\rho$. Denote by $[\chi]$ the corresponding coarse Lyapunov exponent. We will need two preliminary result on circle maps.
\begin{lemma}\label{L:TransitiveActionCircle1}
If $G,H\subset{\rm Homeo}_{+}(\mathbb{T})$ are two path-connected groups such that $GH.x = HG.x = \mathbb{T}$ for all $x\in\mathbb{T}$ then either $G.x = \mathbb{T}$ or $H.x = \mathbb{T}$.
\end{lemma}
\begin{proof}
We have $G.x = (a,b)$, $G.x = (a,b]$, $G.x = [a,b)$ or $G.x = [a,b]$ since $G$ is path-connected. Write $I = G.x$. If $I\neq\mathbb{T}$ then $a$ and $b$ are fixed by $G$ since $g$ is orientation preserving (and therefore fix the endpoints of the $g-$invariant interval $I$). So, $G.a = a$ and $G.b = b$. If $H.x\subset I$ then $\mathbb{T} = GH.x = G.x = I\neq\mathbb{T}$, so $H.x\not\subset I$. Since $H$ is path-connected $J = H.x$ is an interval that contain $x$ and $J\not\subset I$, so $a\in H.x$ or $b\in H.x$. Assume that $a\in H.x$. Since $a$ is fixed by $G$ we have $\mathbb{T} = HG.a = H.a = H.x$ so $H.x = \mathbb{T}$.
\end{proof}
\begin{lemma}\label{L:SubsetOfCircleIsFiniteOrFull}
Let $K\subset\mathbb{T}$ be compact such that 
\begin{enumerate}[label = (\roman*)]
    \item there is a subgroup $\mathcal{G}\subset{\rm Homeo}_{+}(\mathbb{T})$, that acts transitively, with subgroup $\mathcal{G}_{K} = {\rm Stab}_{K}(\mathcal{G}) = \{g\in\mathcal{G}_{K}\text{ : }gK = K\}$,
    \item if $x\in K$ and $g\in\mathcal{G}$ satisfy $gx\in K$ then $g\in\mathcal{G}_{K}$,
    \item if $g\in\mathcal{G}_{K}$ satisfy $gx = x$ for $x\in K$ then $g|_{K} = {\rm id}_{K}$,
    \item there is a compact subset $\mathcal{G}_{0}\subset\mathcal{G}$ such that $\mathcal{G}_{0}x = \mathbb{T}$ for every $x\in\mathbb{T}$.
\end{enumerate}
Then $K = \mathbb{T}$ or $K$ is finite.
\end{lemma}
\begin{proof}
The group $\mathcal{G}_{K}$ act transitively and freely on $K$ by $(i)$, $(ii)$ and $(iii)$. Combining this with $(iv)$ it follows that $\mathcal{G}_{K}$ is a compact group, so $\mathcal{G}_{K}$ preserves a measure $\nu$ on $\mathbb{T}$. From Lemma \ref{L:AdditivityOfRotationNumber} it follows that the rotation number $\omega:\mathcal{G}_{K}\to\mathbb{T}$ is a homomorphism. If $\omega(g) = 0$ for $g\in\mathcal{G}_{K}$, then $g$ fix a point in $K$ by Lemma \ref{L:FixedPointsInCptInvSets}, so by $(iii)$ we have $g|_{K} = {\rm id}_{K}$. It follows that, if we view $\mathcal{G}_{K}\subset{\rm Homeo}(K)$, $\omega:\mathcal{G}_{K}\to\mathbb{T}$ is injective. The image $T = \omega(\mathcal{G}_{K})$ is compact, so either $T = \mathbb{T}$ or $\#T < \infty$. In the second case, since $\omega|_{\mathcal{G}_{K}}$ is injective, we have $\#K = \#\mathcal{G}_{K} < \infty$. In the first case, $\omega:\mathcal{G}_{K}\to\mathbb{T}$ is injective and surjective, so a homeomorphism. It follows that $K$ is homeomorphic to a circle, so $K\hookrightarrow\mathbb{T}$ is a local homeomorphism by invariance of domain. In particular, $K$ is both open and closed. Since $\mathbb{T}$ is connected it follows that $K = \mathbb{T}$.
\end{proof}
\begin{lemma}\label{L:RootActionsCommute}
The action of $\mathcal{P}_{[\chi]}$ commute with the action of $\mathcal{Q}_{[\chi]}$.
\end{lemma}
\begin{proof}
By induction it suffices to consider $v\in E_{0}^{\pm[\chi]}$ and $w\in E_{0}^{[\eta]}$ for $[\eta]$ independent of $[\chi]$. Assume $v\in E_{0}^{[\chi]}$, the other case is identical. Let $g = \eta_{\mathcal{P}}^{v}\eta_{\mathcal{P}}^{w}\eta_{\mathcal{P}}^{-v}\eta_{\mathcal{P}}^{-w}$. Fix
\begin{align*}
\mathcal{C}_{v,w}^{-} = \left\{\mathbf{n}\in\mathbb{Z}^{k}\text{ : }\lim_{\ell\to\infty}\frac{1}{\ell}\log\norm{\rho^{\ell\mathbf{n}}v},\frac{1}{\ell}\log\norm{\rho^{\ell\mathbf{n}}w} < 0\right\}
\end{align*}
which is a non-empty cone since $[\chi]$ and $[\eta]$ are independent. Let $Y\subset\mathbb{T}^{d}$ be the full measure set from Lemma \ref{L:InvariantManifoldsMMEs}. Given $y\in Y$ let $x\in\Phi^{-1}(y)$. We have
\begin{align*}
\alpha^{\ell\mathbf{n}}gx = & \alpha^{\ell\mathbf{n}}\eta_{\mathcal{P}}^{v}\eta_{\mathcal{P}}^{w}\eta_{\mathcal{P}}^{-v}\eta_{\mathcal{P}}^{-w}x = \\ &
\eta_{\mathcal{P}}^{\rho^{\ell\mathbf{n}}v}\eta_{\mathcal{P}}^{\rho^{\ell\mathbf{n}}w}\eta_{\mathcal{P}}^{-\rho^{\ell\mathbf{n}}v}\eta_{\mathcal{P}}^{-\rho^{\ell\mathbf{n}}w}\alpha^{\ell\mathbf{n}}x.
\end{align*}
The action $\eta_{\mathcal{P}}$ is Hölder by Lemma \ref{L:BasicPropertiesLiftedTranslations}, so from or choice of $\mathcal{C}_{v,w}^{-}$ there is uniform $\kappa > 0$ such that for any $\mathbf{n}\in\mathcal{C}_{v,w}^{-}$
\begin{align*}
\limsup_{\ell\to\infty}\frac{1}{\ell}\log\intd_{c}\left(\alpha^{\ell\mathbf{n}}x,\alpha^{\ell\mathbf{n}}gx\right)\leq & \kappa\min_{[\Tilde{\eta}] = [\eta],[\Tilde{\chi}] = [\chi]}\left(\Tilde{\chi}(\mathbf{n}) + \Tilde{\eta}(\mathbf{n})\right)\leq \\ &
\kappa\min_{[\Tilde{\eta}] = [\eta]}\Tilde{\eta}(\mathbf{n}).
\end{align*}
In particular, $gx$ lie in the stable manifold of $x$ for any $\alpha^{\mathbf{n}}$ with $\mathbf{n}\in\mathcal{C}_{v,w}^{-}$. Assume for contradiction that $gx\neq x$. Since $gx,x\in W^{c}(x)$ lie in the same stable manifold Lemma \ref{L:InvariantManifoldsMMEs} implies
\begin{align*}
\limsup_{\ell\to\infty}\frac{1}{\ell}\log\intd_{c}\left(\alpha^{\ell\mathbf{n}}x,\alpha^{\ell\mathbf{n}}gx\right) = \lambda_{\nu}^{c}(\mathbf{n})
\end{align*}
for some $\nu\in\mathcal{M}_{\rm vol}^{f}(X_{\Gamma})$. So, $\lambda_{\nu}^{c}(\mathbf{n})\leq\kappa\min_{[\Tilde{\eta}] = [\eta]}\Tilde{\eta}(\mathbf{n})$ for any $\mathbf{n}\in\mathcal{C}_{v,w}^{-}$. Since $[\chi]$ is independent of $[\eta]$ we find $\mathbf{n}_{j}\in\mathcal{C}_{v,w}^{-}$ such that $\Tilde{\chi}(\mathbf{n}_{j})\to 0$ for any $[\Tilde{\chi}] = [\chi]$ and $\Tilde{\eta}(\mathbf{n}_{j})\to-\infty$ for any $[\Tilde{\eta}] = [\eta]$. Since $\lambda_{\nu}^{c}$ has the same kernel as $[\chi]$ we also have $\lambda_{\nu}^{c}(\mathbf{n}_{j})\to 0$. This is a contradiction since
\begin{align*}
0 = \lim_{j\to\infty}\lambda_{\nu}^{c}(\mathbf{n}_{j})\leq\lim_{j\to\infty}\kappa\min_{[\Tilde{\chi}] = [\chi]}\Tilde{\eta}(\mathbf{n}_{j}) = -\infty.
\end{align*}
We conclude that $gx = x$ for $x\in\Phi^{-1}(Y)$ and since $Y$ has full volume we conclude that $gx = x$ holds on a dense set. The map $g$ is continuous so this implies $g = {\rm id}_{X_{\Gamma}}$.
\end{proof}
By Lemmas \ref{L:TransitiveActionCircle1} and \ref{L:RootActionsCommute} either $\mathcal{P}_{[\chi]}^{c}$ or $\mathcal{Q}_{[\chi]}^{c}$ acts transitively on $W^{c}(x)$ for every $x\in X_{\Gamma}$. We will prove that neither $\mathcal{P}_{[\chi]}^{c}$ or $\mathcal{Q}_{[\chi]}^{c}$ can act transitively on $W^{c}(x)$, which is a contradiction.
\begin{lemma}\label{L:RootActionIsNotTransitive2}
The group $\mathcal{Q}_{[\chi]}^{c}$ can not act transitively on any $W^{c}(x)$.
\end{lemma}
\begin{proof}
Let $Y\subset X_{\Gamma}$ be the subset where $\mathcal{Q}_{[\chi]}^{c}$ act transitively on $W^{c}(y)$. Since $\mathcal{Q}_{[\chi]}^{c}$ is normal in $\mathcal{Q}_{[\chi]}$ (it is the kernel of a homomorphism) $Y$ is $\mathcal{Q}_{[\chi]}-$invariant. But by Lemma \ref{L:RootActionsCommute} the set $Y$ is also $\mathcal{P}_{[\chi]}-$invariant. So $Y$ is $\mathcal{P}$-invariant (by Lemma \ref{L:RootActionsCommute}). Since $f$ is accessible it follows that either $Y = \emptyset$ or $Y = X_{\Gamma}$.

Assume for contradiction that $Y\neq\emptyset$, so $Y = X_{\Gamma}$. The remainder of the proof is similar to the proof of the generic case of Theorem \ref{Thm:InvariantStructureInCenter}. By Theorem \ref{Thm:InvariancePrincipleForHiherRankAction} and Lemma \ref{L:HigherRankActionIsCocycle} the action $\mathcal{Q}_{[\chi]}^{c}$ preserve $\mu$. Since $\mathcal{Q}_{[\chi]}^{c}$ act transitively on every $W^{c}(x)$ the disintegration of $\mu$, $\mu_{x}^{c}$, is not atomic. This is a contradiction, so we conclude that $Y = \emptyset$.
\end{proof}
Before proceeding we define the space
\begin{align}\label{Eq:ProjectionOfComplementaryGroup}
V = \Pi\mathcal{Q}_{[\chi]} = \bigoplus_{\eta\neq\pm[\chi]}E_{0}^{[\eta]}
\end{align}
and the associated translation action on $\mathbb{T}^{d}$ by $R_{v}(x) = x + v$. The action $R_{v}$ is minimal (Lemma \ref{L:Background2}).
\begin{lemma}\label{L:RootActionIsNotTransitive1}
The group $\mathcal{P}_{[\chi]}^{c}$ can not act transitively on any $W^{c}(x_{0})$.
\end{lemma}
The proof of Lemma \ref{L:RootActionIsNotTransitive1} is by contradiction. As in the proof of Lemma \ref{L:RootActionIsNotTransitive2} we may assume that $\mathcal{P}_{[\chi]}^{c}$ act transitively on every $W^{c}(x)$, we will do this in the remainder. We split the proof of Lemma \ref{L:RootActionIsNotTransitive1} into parts, beginning by proving that $\mathcal{Q}_{[\chi]}$ acts minimally.

For $x\in X_{\Gamma}$ define
\begin{align}
K_{x} := \overline{\{\eta_{\mathcal{P}}^{w}x\text{ : }w\in\mathcal{Q}_{[\chi]}\}}.
\end{align}
First, we show that $K_{x}$ is a minimal set for the $\mathcal{Q}_{[\chi]}-$action for every $x\in X_{\Gamma}$.
\begin{lemma}\label{L:MinimalityOfQ-Action}
For any $y\in K_{x}$ we have
\begin{align}
K_{x} = \overline{\{\eta_{\mathcal{P}}^{w}y\text{ : }w\in\mathcal{Q}_{[\chi]}\}},
\end{align}
that is $K_{x}$ is a minimal set for the $\mathcal{Q}_{[\chi]}-$action. In particular, the set $\{K_{x}\}_{x}$ form a partition of $X_{\Gamma}$.
\end{lemma}
\begin{proof}
Define $Y\subset X_{\Gamma}$ as those $x\in X_{\Gamma}$ such that $K_{x}$ is minimal for the $\mathcal{Q}_{[\chi]}-$action. A standard application of Zorn's lemma shows that $Y$ is non-empty. Given any $w\in\mathcal{Q}_{[\chi]}$ we have $\eta_{\mathcal{P}}^{w}K_{x} = K_{x}$, so $Y$ is $\mathcal{Q}_{[\chi]}-$invariant. By Lemma \ref{L:RootActionsCommute} the action of $\mathcal{P}_{[\chi]}$ commute with the action of $\mathcal{Q}_{[\chi]}$ so, given $w\in\mathcal{P}_{[\chi]}$, the map $\eta_{\mathcal{P}}^{w}:K_{x}\to\eta_{\mathcal{P}}^{w}K_{x}$ conjugates the $\mathcal{Q}_{[\chi]}-$action on $K_{x}$ to the $\mathcal{Q}_{[\chi]}-$action on $\eta_{\mathcal{P}}^{w}K_{x} = K_{\eta_{\mathcal{P}}^{w}x}$. It follows that $x\in Y$ if and only if $\eta_{\mathcal{P}}^{w}x\in Y$. Since $Y$ is $\mathcal{P}_{[\chi]}-$invariant and $\mathcal{Q}_{[\chi]}-$invariant Lemma \ref{L:RootActionsCommute} implies that $Y$ is $su-$saturated. Accessibility and the fact that $Y\neq\emptyset$ implies that $Y = X_{\Gamma}$.
\end{proof}
\begin{lemma}\label{L:MinimalityOfRestrictedTranslation}
The action of $\mathcal{Q}_{[\chi]}$ is minimal.
\end{lemma}
\begin{proof}
Since $K_{x}$ is $\mathcal{Q}_{[\chi]}-$invariant, Lemma \ref{L:PropertiesSUaction1} implies
\begin{align}
\Phi(K_{x}) = \Phi(\eta_{\mathcal{P}}^{w}K_{x}) = \Phi(K_{x}) + \Pi(w),\quad w\in\mathcal{Q}_{[\chi]}
\end{align}
so $\Phi(K_{x})$ is compact and invariant by $R_{v}$ from Equation \ref{Eq:ProjectionOfComplementaryGroup}. The action $R_{v}$ is minimal so $\Phi(K_{x}) = \mathbb{T}^{d}$. It follows that $K_{x}\cap W^{c}(y)\neq\emptyset$ for every $y\in X_{\Gamma}$.

If $\mathcal{G} = \eta_{\mathcal{P}}^{\mathcal{P}_{[\chi]}^{c}}$ and $K = K_{x}\cap W^{c}(y)$ then the assumptions in Lemma \ref{L:SubsetOfCircleIsFiniteOrFull} are satisfied (after identifying $W^{c}(y)\cong\mathbb{T}$). Indeed, $\mathcal{G}$ act transitively by Lemma \ref{L:PropertiesSUaction3}. Property $(ii)$ and $(iii)$ follows from Lemmas \ref{L:RootActionsCommute} and \ref{L:MinimalityOfQ-Action}. To show property $(iv)$ we follow \cite[Section 8.3]{AvilaVianaSantamaria2013}. For $N\in\mathbb{N}$, write 
\begin{align}
\mathcal{P}_{[\chi],N}^{c} = \{w = v_{1}v_{2}...v_{N}\in\mathcal{P}_{[\chi]}^{c}\text{ : }v_{j}\in E_{0}^{\pm[\chi]},\text{ }\norm{v_{j}}\leq N,\text{ }n\leq N\}.
\end{align}
That is $\mathcal{P}_{[\chi],N}^{c}$ consist of those words in $\mathcal{P}_{[\chi]}^{c}$ with at most $N$ letters, and each letter has length of at most $N$. The following is immediate from the definition
\begin{align}
\mathcal{P}_{[\chi],N+1}^{c}\supset\mathcal{P}_{[\chi],N}^{c},\quad\mathcal{P}_{[\chi]}^{c} = \bigcup_{N\geq 1}\mathcal{P}_{[\chi],N}^{c}.
\end{align}
Let $\mathcal{K}_{N} = \{\eta_{\mathcal{P}}^{w}\text{ : }w\in\mathcal{P}_{[\chi],N}^{c}\}$ and $K_{N} = \mathcal{K}_{N}x_{0} = \{k(x_{0})\text{ : }k\in\mathcal{K}_{N}\}$ for fixed $x_{0}\in W^{c}(y)$. Then $K_{N}\subset K_{N+1}$ is an ascending sequence of compact sets. If $K_{N_{0}} = W^{c}(y)$ and $x\in W^{c}(x)$ then $x_{0}\in\mathcal{K}_{N_{0}}x$ so $\mathcal{K}_{2N_{0}}x\supset\mathcal{K}_{N_{0}}x_{0} = W^{c}(y)$. That is, property $(iv)$ of Lemma \ref{L:SubsetOfCircleIsFiniteOrFull} holds if $K_{N_{0}} = W^{c}(y)$ for some $N_{0}$. Since
\begin{align}\label{Eq:AllLengthLegIsFull}
\bigcup_{N\geq1}K_{N} = \{\eta_{\mathcal{P}}^{w}x_{0}\text{ : }w\in\mathcal{P}_{[\chi]}^{c}\} = W^{c}(y),
\end{align}
Baire's category theorem implies that there is $N_{1}$ such that ${\rm int}(K_{N})\neq\emptyset$ for $N\geq N_{1}$ (where ${\rm int}(K_{N})$ is the interior in $W^{c}(y)$). Define $U_{N} = {\rm int}(K_{N})$ and $C_{N} = U_{N}^{c}$ (where the complement is in $W^{c}(y)$). If $N\geq N_{1}$ then $x_{0}\in\mathcal{K}_{N}U_{N}\subset U_{2N}$, so $K_{N}\subset\mathcal{K}_{N}U_{2N}\subset U_{3N}$. From Equation \ref{Eq:AllLengthLegIsFull} we obtain
\begin{align}
W^{c}(y) = \bigcup_{N\geq1}K_{N} = \bigcup_{N\geq N_{1}}K_{N}\subset\bigcup_{N\geq N_{1}}U_{3N} = \bigcup_{N\geq1}U_{N}
\end{align}
or
\begin{align}
\emptyset = \bigcap_{N\geq 1}C_{N}.
\end{align}
Since $C_{N+1}\subset C_{N}$ is a descending sequence of compact sets there is $N_{0}$ such that $C_{N_{0}} = \emptyset$. Equivalently $U_{N_{0}} = K_{N_{0}} = W^{c}(y)$ proving that property $(iv)$ of Lemma \ref{L:SubsetOfCircleIsFiniteOrFull} holds. Lemma \ref{L:SubsetOfCircleIsFiniteOrFull} implies that $\#K_{x}\cap W^{c}(y) < \infty$ or $K_{x}\cap W^{c}(y) = W^{c}(y)$ for every $x,y\in X_{\Gamma}$.

We claim that $\# K_{x}\cap W^{c}(y)$ is independent of $y$. Indeed, as in the proof of Lemma \ref{L:MinimalityOfQ-Action}, for $w\in\mathcal{P}$ we have $\eta_{\mathcal{P}}^{w}K_{x} = K_{\eta_{\mathcal{P}}^{w}x}$. So given $z\in K_{x}\cap W^{c}(y)$ we find, by accessibility, some $w\in\mathcal{P}$ that satisfy $\eta_{\mathcal{P}}^{w}x = z$ and in extension $\eta_{\mathcal{P}}^{w}K_{x} = K_{z} = K_{x}$ (since $z\in K_{x}$, and $\{K_{x}\}_{x\in X_{\Gamma}}$ is a partition, Lemma \ref{L:MinimalityOfQ-Action}). We also have $\eta_{\mathcal{P}}^{w}W^{c}(x) = W^{c}(y)$ (since $z\in W^{c}(y)$) so $\eta_{\mathcal{P}}^{w}(K_{x}\cap W^{c}(x)) = K_{x}\cap W^{c}(y)$. It follows that
\begin{align}
\# K_{x}\cap W^{c}(y) = \#K_{x}\cap W^{c}(x)
\end{align}
for all $y\in X_{\Gamma}$. If $\# K_{x}\cap W^{c}(x) = q < \infty$, then the fibers of $\Phi$ intersect $K_{x}$ precisely $q$ times. Given $z_{0}\in\mathbb{T}^{d}$ we define
\begin{align*}
\varepsilon_{0} = \inf_{y\neq y'\text{ }y,y'\in\Phi^{-1}(z_{0})\cap K_{x}}\intd(y,y') > 0.
\end{align*}
Let $\{y_{1},...,y_{q}\} = K_{x}\cap\Phi^{-1}(z_{0})$ and $\delta > 0$ small. For $z\in B_{\delta}(z_{0})$, we define $y_{j}(z)$ as the element in $\Phi^{-1}(z)\cap K_{x}$ that minimize $\intd(y_{j}(z),y_{j}(z_{0}))$ (if this does not define a unique point, choose an arbitrary minimizer). If $z_{n}\to z_{0}$ then any convergent subsequence of $y_{j}(z_{n})$ converges to an element of $K_{x}\cap\Phi^{-1}(z_{0})$. From our choice of $y_{j}(z)$ it is clear that any convergent subsequence converge to $y_{j}(z_{0})$. So, $y_{j}(z_{n})\to y_{j}(z_{0})$ and $y_{j}$ is continuous at $z_{0}$. Choose $\delta > 0$ small enough such that $\intd(y_{j}(z),y_{j}(z_{0})) < \varepsilon_{0}/100$ for all $z\in B_{\delta}(z_{0})$ and $j = 1,...,q$. For $z,z'\in B_{\delta}(z_{0})$ and $i\neq j$ the reverse triangle inequality implies
\begin{align*}
\intd(y_{i}(z),y_{j}(z'))\geq & \intd(y_{i}(z_{0}),y_{j}(z')) - \intd(y_{i}(z_{0}),y_{i}(z))\geq \\ &
\intd(y_{i}(z_{0}),y_{j}(z_{0})) - \intd(y_{j}(z_{0}),y_{j}(z')) - \intd(y_{i}(z_{0}),y_{i}(z)) > \\ &
\varepsilon_{0} - \frac{\varepsilon_{0}}{50} > \frac{49\varepsilon_{0}}{50} > 0.
\end{align*}
where the last inequality use the definition of $\varepsilon_{0}$. In particular, $y_{i}(z)\neq y_{j}(z)$ for every $z\in B_{\delta}(z_{0})$, so $K_{x}\cap\Phi^{-1}(z) = \{y_{1}(z),...,y_{q}(z)\}$. Moreover, we have
\begin{align}
\intd(y_{i}(z),y_{i}(z'))\leq\intd(y_{i}(z),y_{i}(z_{0})) + \intd(y_{i}(z_{0}),y_{i}(z')) < \frac{\varepsilon_{0}}{50}.
\end{align}
Since $\varepsilon_{0}/50 < 49\varepsilon_{0}/50$, $y_{i}(z')$ is the element $y\in K_{x}\cap\Phi^{-1}(z')$ that minimize $\intd(y_{i}(z),y)$. The same argument that showed continuity at $z_{0}$, now show that $y_{j}$ is continuous at any $z\in B_{\delta}(z_{0})$. So the functions $y_{1},...,y_{q}:B_{\delta}(z_{0})\to K_{x}$ are continuous. Note that $\Phi$ restricted to $K_{x}\cap\Phi^{-1}(B_{\delta}(z_{0}))\cap B_{\varepsilon_{0}}(y_{j}(z_{0}))$ has inverse given by $y_{j}$. So, $\Phi:K_{x}\to\mathbb{T}^{d}$ is a finite covering map. In particular, $K_{x}$ is homeomorphic to $\mathbb{T}^{d}$. It follows that the map
\begin{align}
(\Phi|_{K_{x}})_{*}:\pi_{1}\mathbb{T}^{d}\cong\mathbb{Z}^{d}\to\Gamma\cong\pi_{1}X_{\Gamma}
\end{align}
is injective with image of finite index. The map
\begin{align}
\mathbb{Z}^{d}\cong\pi_{1}K_{x}\xrightarrow{i_{K_{x}}}\pi_{1}X_{\Gamma} = \Gamma
\end{align}
is injective since its injective after composition with $\Phi_{*} = \pi_{*}$. The group $\Gamma' = {\rm Im}(i_{K_{x}})_{*}\times\ker\Phi_{*} = {\rm Im}(i_{K_{x}})_{*}\times[\Gamma,\Gamma]$ has finite index in $\Gamma$ and is abelian, but $\Gamma$ is not virtually abelian so this is a contradiction. We conclude that $K_{x}\cap W^{c}(y) = W^{c}(y)$ for every $y$ which implies $K_{x} = X_{\Gamma}$ for every $x\in X_{\Gamma}$.
\end{proof}
\begin{lemma}\label{L:RootActionIsNotTransitive1:1}
The map $r:\mathcal{P}_{[\chi]}^{c}\to\mathbb{T}$, $r(w) = \omega(\eta_{\mathcal{P}}^{w}|_{W^{c}(x)})$ is a well-defined surjective homomorphism. Moreover, $r(w) = 0$ if and only if $\eta_{\mathcal{P}}^{w} = {\rm id}_{X_{\Gamma}}$.
\end{lemma}
\begin{proof}
Fix $w\in\mathcal{P}_{[\chi]}^{c}$, $g = \eta_{\mathcal{P}}^{w}$, and define $\varphi:\mathbb{T}^{d}\to\mathbb{T}$ by 
\begin{align}
\varphi(y) := \omega(g|_{\Phi^{-1}(y)}).
\end{align}
The function $\varphi$ is continuous since rotation numbers vary continuously in the $C^{0}-$topology \cite[Proposition 11.1.6]{KatokHasselblatt1995}. Given $w\in\mathcal{Q}_{[\chi]}$ we have $g\eta_{\mathcal{P}}^{w} = \eta_{\mathcal{P}}^{w}g$ (by Lemma \ref{L:RootActionsCommute}) and $\Phi\eta_{\mathcal{P}}^{w} = R_{\Pi(w)}\Phi$ (by Lemma \ref{L:PropertiesSUaction1}). Let $V$ be the space defined in Equation \ref{Eq:ProjectionOfComplementaryGroup}. Given $v\in V$ let $w\in\mathcal{Q}_{[\chi]}$ be such that $\Pi(w) = v$. We have
\begin{align*}
\varphi(R_{v}(x)) = & \omega(g|_{\Phi^{-1}(R_{v}(x))}) = \omega(g|_{\eta_{\mathcal{P}}^{w}\Phi^{-1}(x)}) = \\ &
\omega(\eta_{\mathcal{P}}^{w}g\eta_{\mathcal{P}}^{w^{-1}}|_{\eta_{\mathcal{P}}^{w}\Phi^{-1}(x)}) = \omega(g|_{\Phi^{-1}(x)}) = \varphi(x)
\end{align*}
where the second to last equality uses that the rotation number is conjugacy invariant \cite[Proposition 11.1.3]{KatokHasselblatt1995}. Since $\varphi$ is invariant under $R_{v}$ and $R_{v}$ is minimal $\varphi$ is constant. That is, $g:W^{c}(x)\to W^{c}(x)$ has a rotation number independent of $x\in X_{\Gamma}$. This shows that $r:\mathcal{P}_{[\chi]}^{c}\to\mathbb{T}$ is well-defined. If $r(w) = 0$ then $\eta_{\mathcal{P}}^{w}$ fix some $x$, so Lemmas \ref{L:RootActionsCommute} and \ref{L:MinimalityOfRestrictedTranslation} implies that $\eta_{\mathcal{P}}^{w} = {\rm id}_{X_{\Gamma}}$.

Since $H = \eta_{\mathcal{P}}|_{\mathcal{P}_{[\chi]}^{c}}\cong\mathcal{P}_{[\chi]}^{c}/\{w\in\mathcal{P}_{[\chi]}^{c}\text{ : }\eta_{\mathcal{P}}^{w} = {\rm id}_{X_{\Gamma}}\}$ act transitively and freely on each $W^{c}(x)$, it follows that $H$ is homeomorphic to $W^{c}(x)\cong\mathbb{T}$. So $H$ is a compact topological group. Since $H$ is compact the $H-$action preserve a measure. That $r$ is a homomorphism follows from Lemma \ref{L:AdditivityOfRotationNumber}. The image $r(H)\subset\mathbb{T}$ is a compact subgroup, so it is either finite or all of $\mathbb{T}$. But $r(H)$ can not be finite since $r$ is injective and $H$ acts transitively.
\end{proof}
We can now prove Lemma \ref{L:RootActionIsNotTransitive1}.
\begin{proof}[Proof of Lemma \ref{L:RootActionIsNotTransitive1}]
Suppose that $\mathcal{P}_{[\chi]}^{c}$ act transitively on some $W^{c}(x_{0})$. By Lemma \ref{L:RootActionIsNotTransitive1:1} we can define an action $\eta_{c}:\mathbb{T}\times X_{\Gamma}\to X_{\Gamma}$ by $\eta_{c}^{s}x = \eta_{\mathcal{P}}^{w}x$ with $w\in\mathcal{P}^{c}$ chosen such that $r(w) = s$. The action $\eta_{c}$ is free and transitive on each center leaf $W^{c}(x)$. Moreover, for $\mathbf{n}\in\mathbb{Z}^{k}$ and $w\in\mathcal{P}_{[\chi]}^{c}$ we have
\begin{align}
\eta_{\mathcal{P}}^{\rho^{\mathbf{n}}w} = \alpha^{\mathbf{n}}\eta_{\mathcal{P}}^{w}\alpha^{-\mathbf{n}}
\end{align}
so by conjugacy invariance of the rotation number $r(\rho^{\mathbf{n}}w) = r(w) = s$. It follows that
\begin{align}
\eta_{c}^{s} = \eta_{\mathcal{P}}^{\rho^{\mathbf{n}}w} = \alpha^{\mathbf{n}}\eta_{\mathcal{P}}^{w}\alpha^{-\mathbf{n}} = \alpha^{\mathbf{n}}\eta_{c}^{s}\alpha^{-\mathbf{n}}
\end{align}
which proves that $\alpha$ commute with $\eta_{c}$. This implies that $\alpha$ have vanishing Lyapunov exponent along $W^{c}$, which is a contradiction. We conclude that $\mathcal{P}_{[\chi]}^{c}$ can not act transitively on $W^{c}(x_{0})$.
\end{proof}
We finish the proof of Theorem \ref{Thm:InvariantStructureInCenter} in the exceptional case.
\begin{proof}[Proof of Theorem \ref{Thm:InvariantStructureInCenter} in exceptional case]
By Lemma \ref{L:RootActionsCommute} the action of $\mathcal{P}_{[\chi]}^{c}$ commute with the action of $\mathcal{Q}_{[\chi]}^{c}$. Since $f$ is accessible $\eta_{\mathcal{P}}^{\mathcal{P}_{[\chi]}^{c}}\eta_{\mathcal{P}}^{\mathcal{Q}_{[\chi]}^{c}}x = W^{c}(x)$ for all $x\in X_{\Gamma}$. By Lemma \ref{L:TransitiveActionCircle1} either $\mathcal{P}_{[\chi]}^{c}$ or $\mathcal{Q}_{[\chi]}^{c}$ act transitively on $W^{c}(x)$. Lemma \ref{L:RootActionIsNotTransitive1} shows that $\mathcal{P}_{[\chi]}^{c}$ does not act transitively on $W^{c}(x)$ so $\mathcal{Q}_{[\chi]}^{c}$ must act transitively on $W^{c}(x)$. This is a contradiction by Lemma \ref{L:RootActionIsNotTransitive2}.
\end{proof}

\section{Compatible algebraic structure}\label{Sec:CompatibleAlgebraicStructure}

In this section we produce a topological conjugacy between $\alpha$ and an affine action $\alpha_{0}$ (Theorem \ref{Thm:ExistenceBiHölderConj}). Let $\alpha:\mathbb{Z}^{k}\times X_{\Gamma}\to X_{\Gamma}$ and $f:X_{\Gamma}\to X_{\Gamma}$ satisfy the assumptions of Theorem \ref{Thm:MainTheorem1}. By Theorem \ref{Thm:InvariantStructureInCenter} there is a unique $\alpha-$invariant measure $\mu$, with center exponent $\lambda_{\mu}^{c} = 0$, that projects to volume.

We begin by using the measure $\mu$ to construct a circle action $\eta_{c}:\mathbb{T}\times X_{\Gamma}\to X_{\Gamma}$ that commutes with $\alpha$ and $\eta_{\mathcal{P}}$, preserves center leaves and acts transitively and freely on each center leaf. Moreover, if $w\in\mathcal{P}^{c}$ then $\eta_{\mathcal{P}}^{w} = \eta_{c}^{t}$ for some $t\in\mathbb{T}$. This shows that the action of $\mathcal{P}$ on $X_{\Gamma}$ factor through a nilpotent group with base $E_{0}^{u}\oplus E_{0}^{s}$ and center isomorphic to $\mathbb{T}$. Lifting this action to $G$ we obtain a transitive free action of a $2-$step nilpotent group, $N$, on $G$.

Using that $\lambda_{\mu}^{c} = 0$, the following two lemmas are proved identically as Lemmas \ref{L:RootActionsCommute} and \ref{L:RootActionIsNotTransitive1:1}.
\begin{lemma}\label{L:CompAlgStruc1}
Let $[\chi]$ and $[\eta]$ be independent. If $w_{1}\in\mathcal{P}_{[\chi]}$ and $w_{2}\in\mathcal{P}_{[\eta]}$ then 
\begin{align}
\eta_{\mathcal{P}}^{w_{1}}\eta_{\mathcal{P}}^{w_{2}} = \eta_{\mathcal{P}}^{w_{2}}\eta_{\mathcal{P}}^{w_{1}}.
\end{align}
\end{lemma}
\begin{lemma}\label{L:CompAlgStruc2}
For any $w\in\mathcal{P}^{c}$ the rotation number of 
\begin{align}
\eta_{\mathcal{P}}^{w}:W^{c}(x)\to W^{c}(x)
\end{align}
is independent of $x$. The map $r(w)$ mapping $w\in\mathcal{P}^{c}$ to the rotation number of $\eta_{\mathcal{P}}^{w}$ is a homomorphism with kernel $\{w\in\mathcal{P}^{c}\text{ : }\eta_{\mathcal{P}}^{w} = {\rm id}_{X_{\Gamma}}\}$.
 \end{lemma}
Using Lemmas \ref{L:CompAlgStruc1} and \ref{L:CompAlgStruc1} we produce an $N-$action, with $N$ nilpotent, on $G$.
\begin{lemma}\label{L:CompAlgStruc3}
The action of $\mathcal{P}$ on $G$ factor through a nilpotent Lie group $N$ that acts transitively and freely on $G$, the action of $N$ descends to an action on $X_{\Gamma}$.
\end{lemma}
\begin{proof}
Denote by $N$ the image of $\mathcal{P}$ in ${\rm Homeo}(G)$ and $\Tilde{N}$ the image of $\mathcal{P}$ in ${\rm Homeo}(X_{\Gamma})$. We begin by showing that the $\mathcal{P}^{c}-$action on $X_{\Gamma}$ factor through an abelian group. As in the proof of Lemma \ref{L:RootActionIsNotTransitive1} we define $\eta_{c}^{s}:X_{\Gamma}\to X_{\Gamma}$ by $\eta_{c}^{s} = \eta_{\mathcal{P}}^{w}$ for any $w\in\mathcal{P}^{c}$ that satisfies $\omega(\eta_{\mathcal{P}}^{w}|_{W^{c}(x)}) = s$. By Lemma \ref{L:CompAlgStruc2} $\eta_{c}$ is well-defined, continuous, and acts transitively and freely on each center leaf $W^{c}(x)$.

As in the proof of Lemma \ref{L:RootActionIsNotTransitive1}, $\alpha$ commutes with $\eta_{c}$. Let $r$ be from Lemma \ref{L:CompAlgStruc2}. Given $w\in\mathcal{P}$ and $w_{c}\in\mathcal{P}^{c}$ we have $r(ww_{c}w^{-1}) = r(w_{c})$ since $\eta_{\mathcal{P}}^{w}\eta_{\mathcal{P}}^{w_{c}}(\eta_{\mathcal{P}}^{w})^{-1}$ has the same rotation number as $\eta_{\mathcal{P}}^{w_{c}}$. It follows that
\begin{align}
\eta_{c}^{r(w_{c})} = \eta_{\mathcal{P}}^{w}\eta_{\mathcal{P}}^{w_{c}}(\eta_{\mathcal{P}}^{w})^{-1} = \eta_{\mathcal{P}}^{w}\eta_{c}^{r(w_{c})}(\eta_{\mathcal{P}}^{w})^{-1}
\end{align}
which shows that $\eta_{c}$ commute with $\eta_{\mathcal{P}}$. Let $\Tilde{r}:\mathcal{P}^{c}\to\mathbb{R}$ be defined by $\Tilde{r}(w) = \Tilde{\omega}(\eta_{\mathcal{P}}^{w})$ where $\Tilde{\omega}(\eta_{\mathcal{P}}^{w})\in\mathbb{R}$ is the rotation number of $\eta_{\mathcal{P}}^{w}$ with respect to the natural lift to $G$ from Definition \ref{Def:ActionOfsuPathGroup}. Lift $\eta_{c}$ to $G$, then $\eta_{\mathcal{P}}^{w} = \eta_{c}^{\Tilde{r}(w)}$. It is immediate from the analogous properties of $\eta_{c}:\mathbb{T}\times X_{\Gamma}\to X_{\Gamma}$ that $\eta_{c}:\mathbb{R}\times G\to G$ commute with $\alpha$ and $\eta_{\mathcal{P}}$.

Let $N^{c}\subset N$ be the image of $\mathcal{P}^{c}$. Since $\eta_{c}$ commute with $\eta_{\mathcal{P}}$ and since $\eta_{\mathcal{P}}^{w} = \eta_{\mathcal{P}}^{\Tilde{r}(w)}$ for $w\in\mathcal{P}^{c}$ the subgroup $N^{c}$ is central in $N$. Moreover, Lemma \ref{L:CommutatorOfFreeProduct} implies that
\begin{align}
\mathcal{P}^{c} = \left[\mathcal{P},\mathcal{P}\right]
\end{align}
so $N^{c} = [N,N]$. It follows that $N$ is $2-$step nilpotent. By Lemma \ref{L:CommutatorOfFreeProduct} we can write any $w\in\mathcal{P}$ as $w = w_{c}v_{s}v_{u}$ with $v_{\sigma}\in E_{0}^{\sigma}$, $\sigma = s,u$, and $w_{c}\in\mathcal{P}^{c}$. It follows that $\eta_{\mathcal{P}}^{w} = \eta_{c}^{\Tilde{r}(w_{c})}\eta_{s}^{v_{s}}\eta_{u}^{v_{u}}$, so $T:E_{0}^{s}\times E_{0}^{u}\times\mathbb{R}\to N$, $T(v_{s},v_{u},t) = \eta_{c}^{t}\eta_{s}^{v_{s}}\eta_{u}^{v_{u}}$ is surjective. If $\eta_{c}^{t}\eta_{s}^{v_{s}}\eta_{u}^{v_{u}} = \eta_{c}^{t'}\eta_{s}^{v_{s}'}\eta_{u}^{v_{u}'}$ then we apply $\Phi$ and obtain $v_{s} = v_{s}'$, $v_{u} = v_{u}'$, Lemma \ref{L:PropertiesSUaction1}. If $v_{s} = v_{s}'$, $v_{u} = v_{u}'$ then we can simplify and obtain $\eta_{c}^{t} = \eta_{c}^{t'}$. Since $\eta_{c}$ is a free action this implies $t = t'$. So $T:E_{0}^{s}\times E_{0}^{u}\times\mathbb{R}\to N$ is bijective. It is also clear that $T$ is continuous. Let $(v_{s},v_{u},t)$ and $(v_{s}',v_{u}',t')$ be such that $T(v_{s},v_{u},t)$ is close to $T(v_{s}',v_{u}',t')$. After applying $\Phi$, it is immediate that $v_{s}$ is close to $v_{s}'$ and $v_{u}$ is close to $v_{u}'$. Writing $h = T(v_{s},v_{u},t)T(v_{s}',v_{u}',t')^{-1}$ we obtain
\begin{align*}
h = \eta_{s}^{v_{s}}\eta_{u}^{v_{u}-v_{u}'}\eta_{s}^{v_{s}'}\eta_{c}^{t-t'},\quad\eta_{c}^{t-t'} = \eta_{s}^{-v_{s}'}\eta_{u}^{-v_{u}+v_{u}'}\eta_{s}^{-v_{s}}h
\end{align*}
so $\eta_{c}^{t-t'}$ is close to ${\rm id}_{G}$. Since $\eta_{c}$ act freely $t - t'$ is close to $0$, so $T$ has a continuous inverse and is therefore a homeomorphism. Since $N$ is a topological group homeomorphic to $E_{0}^{s}\times E_{0}^{u}\times\mathbb{R}$, $N$ has a unique structure as a Lie group \cite{MontgomeryZippin1952,Gleason1952}. Finally, $N$ is $2-$step nilpotent as an abstract group so $N$ is a $2-$step nilpotent Lie group.
\end{proof}
\begin{theorem}\label{Thm:ExistenceBiHölderConj}
The diffeomorphism $f:X\to X$ is bi-Hölder conjugate to some affine map $f_{0}:X\to X$ where $f_{0}(x) = L(x)z_{0}$ with $z_{0}\in G^{c}$.
\end{theorem}
\begin{proof}
Let $N$ be the group from Lemma \ref{L:CompAlgStruc3} and $F$ a lift of $f$. For $w\in\mathcal{P}$ we have $F\eta_{\mathcal{P}}^{w}F^{-1} = \eta_{\mathcal{P}}^{L_{su}w}$ by Lemma \ref{L:SolvableRelations}. It follows that $\Tilde{L}(n) = FnF^{-1}\in N$ for every $n\in N$. So, $\Tilde{L}:N\to N$ is a continuous autormophism such that $Fn = \Tilde{L}(n)F$ for $n\in N$.

Define $\Lambda\subset N$ by $\lambda\in\Lambda$ if $\lambda(e)\in\Gamma$ for the identity element $e\in G$. Any two $\lambda,\lambda'\in\Lambda$ are lifted from $X_{\Gamma}$ so there is $\gamma\in\Gamma$ such that $\lambda'(\lambda(e)) = \lambda'(\gamma) = \gamma\lambda'(e)\in\Gamma$. Similarly if $\lambda\in\Lambda$ then $\lambda^{-1}(\lambda(e)) = \lambda(e)\lambda^{-1}(e)$ so $\lambda^{-1}(e) = (\lambda(e))^{-1}\in\Gamma$. It follows that $\Lambda$ is a subgroup of $N$. Moreover, $\Gamma$ is closed in $G$ so $\Lambda\leq N$ is closed.

The map $q:N\ni h\mapsto h(e)\in G$ is a homeomorphism by invariance of domain since $q$ is bijective. Let $h\in N$ and $\lambda\in\Lambda$, then 
\begin{align*}
q(h\lambda) = h(\lambda(e)) = \lambda(e)h(e) = \lambda(e)q(h) = q(\lambda)q(h)
\end{align*}
so $q$ descends to a map $N/\Lambda\to\Gamma\setminus G = X_{\Gamma}$. If $q(h) = \gamma q(h')$ for $h,h'\in N$ and $\gamma\in\Gamma$ then we find $\lambda\in\Lambda$ such that $\gamma = \lambda(e)$ (since the action of $N$ on $G$ is transitive). It follows that
\begin{align*}
q(h) = \lambda(e)q(h') = \lambda(e)h'(e) = h'(\lambda(e)) = q(h'\lambda)
\end{align*}
or $h(e) = h'(\lambda(e))$. The action of $N$ is free, so $h = h'\lambda$. That is, we have a diagram
\begin{center}
        \begin{tikzcd}[row sep=large,column sep=huge]
            N\arrow{r}{q}\arrow{d} & G\arrow{d} \\
            N/\Lambda\arrow{r}{q} & X_{\Gamma}
        \end{tikzcd}
\end{center}
where both horizontal maps are homeomorphisms. The group $N$ is a nilpotent Lie group by Lemma \ref{L:CompAlgStruc3}, so $N/\Lambda$ is a nilmanifold. It follows that $N$ is isomorphic to $G$ as a Lie group and $\Lambda$ is isomorphic to $\Gamma$ under this map $N\to G$ \cite{DiscreteSubgroupsOfLieGroups1972}.

We claim that the induced map $\Tilde{f}:N/\Lambda\to N/\Lambda$, $\Tilde{f} = q^{-1}\circ f\circ q$ is affine. Equivalently, $\Tilde{F}:N\to N$, $\Tilde{F} = q^{-1}\circ F\circ q$ is affine. From the relation $Fn = (\Tilde{L}n)F$, $n\in N$, we obtain
\begin{align*}
q\left(\Tilde{F}(n)\right) = F(q(n)) = F(n(e)) = \Tilde{L}(n)(Fe) = \Tilde{L}(n)\left(n_{0}(e)\right) = q\left(\Tilde{L}(n)\cdot n_{0}\right)
\end{align*}
where $n_{0}\in N$ is chosen such that $n_{0}(e) = Fe$. It follows that $\Tilde{F}$, and therefore also $\Tilde{f}$, is affine.

Denote by $H:X_{\Gamma}\to X_{\Gamma}$ the conjugacy such that $H(fx) = f_{0}(h(x)) = L(H(x))z_{0}$. After conjugating with a translation we may assume that $z_{0}\in G^{c} = [G,G]$. To finish the proof we show that $H$ is bi-Hölder. By uniqueness of $\Phi$ we have $\Phi(x) = \pi(H(x)) = H(x)G^{c}$. Since $\Phi:X_{\Gamma}\to\mathbb{T}^{d}$ is Hölder and the inverse of $\Phi$ restricted to stable and unstable leaves is Hölder (Lemma \ref{L:FrankManningCoordinates2}) $H$ is bi-Hölder along $W^{s}$ and $W^{u}$. So it suffices to show that $H$ is Hölder along $W^{c}(x)$ and $H^{-1}$ is Hölder along $W_{0}^{c}(x)$. Write
\begin{align}
H(x) = xe^{-h(x)} = xe^{-\left(h_{s}(x) + h_{c}(x) + h_{u}(x)\right)}
\end{align}
with $h:X_{\Gamma}\to\mathfrak{g}$ and $h_{\sigma}:X_{\Gamma}\to E_{0}^{\sigma}$, $\sigma = s,c,u$. Since $\pi(H(x)) =\Phi(x)$ and $\Phi$ is Hölder, both $h_{s}$ and $h_{u}$ are Hölder. The functional equation for $H$ implies
\begin{align*}
H(fx) = & (fx)e^{-h(fx)} = L(x)e^{-v(x)}e^{-h(fx)} = L(x)e^{-v(x) - h(fx) + \frac{[v(x),h(fx)]}{2}} = \\ &
L(H(x))z_{0} = L(x)e^{-Lh(x) + Z_{0}}
\end{align*}
or
\begin{align}
h(fx) - Lh(x) = \frac{[v(x),h(fx)]}{2} - v(x) + Z_{0}.
\end{align}
We obtain an equation for $h_{c}(x)$
\begin{align*}
h_{c}(fx) - h_{c}(x) = \frac{[v_{s}(x),h_{u}(fx)] + [v_{u}(x),h_{s}(fx)]}{2} - v_{c}(x) + Z_{0} = w(x)
\end{align*}
with $w:X_{\Gamma}\to E_{0}^{c}$ Hölder. It follows that $h_{c}(x)$ is Hölder \cite{Wilkinson2008}, so $H$ is Hölder. Next we show that $H^{-1}$ is Hölder along $W_{0}^{c}(x)$. Fix $X\in E_{0}^{s}$, $Y\in E_{0}^{u}$ such that $[X,Y] = Z\neq0$ with $Z\in E_{0}^{c}$. Given $t\geq 0$ the Baker–Campbell–Hausdorff formula implies $e^{tX}e^{tY}e^{-tX}e^{-tY} = e^{t^{2}[X,Y]}$, so
\begin{align*}
H^{-1}\left(xe^{t^{2}Z}\right) - H^{-1}\left(x\right) = & H^{-1}\left(xe^{tX}e^{tY}e^{-tX}e^{-tY}\right) - H^{-1}\left(x\right) = \\ &
H^{-1}\left(xe^{tX}e^{tY}e^{-tX}e^{-tY}\right) - H^{-1}\left(xe^{tX}e^{tY}e^{-tX}\right) + \\ &
H^{-1}\left(xe^{tX}e^{tY}e^{-tX}\right) - H^{-1}\left(xe^{tX}e^{tY}\right) + \\ &
H^{-1}\left(xe^{tX}e^{tY}\right) - H^{-1}\left(xe^{tX}\right) + \\ &
H^{-1}\left(xe^{tX}\right) - H^{-1}\left(x\right).
\end{align*}
Letting $t$ be close to $0$ and using that $H^{-1}$ is Hölder along $W_{0}^{s}$ and $W_{0}^{u}$ we obtain
\begin{align*}
\intd\left(H^{-1}\left(xe^{t^{2}Z}\right),H^{-1}\left(x\right)\right)\leq Ct^{\theta}.
\end{align*}
Since $\intd(xe^{t^{2}Z},x)\geq ct^{2}$ for small $t$ we have
\begin{align*}
\intd\left(H^{-1}\left(xe^{t^{2}Z}\right),H^{-1}\left(x\right)\right)\leq C\intd\left(x,xe^{t^{2}Z}\right)^{\theta/2},
\end{align*}
so for $y\in W_{0}^{c}(x)$ close. It follows that $\intd(H^{-1}(x),H^{-1}(y))\leq C\intd(x,y)^{\theta/2}$.
\end{proof}

\section{Rigidity: Smoothness of the bi-Hölder conjugacy}\label{Sec:ImprovedRegularity}

In this section, we prove Theorem \ref{Thm:ImprovedRegularity1}. Let $\alpha$, $\alpha_{0}$, $f = \alpha^{\mathbf{n}_{0}}$, $f_{0} = \alpha_{0}^{\mathbf{n}_{0}}$, and $H$ be as in Theorem \ref{Thm:ImprovedRegularity1}. We prove Theorem \ref{Thm:ImprovedRegularity1} under the assumption that $f$ is accessible, this is done for two reasons. First, the proof simplifies because we can apply \cite{Wilkinson2008} to obtain regularity of the conjugacy along the center direction. In particular, there is no loss of generality in assuming that the center of $f_{0}$ coincides with the joint center of $\alpha_{0}$ (see Remark \ref{Rmk:JointCenter}). Second, if $f$ is accessible then $\alpha$ naturally preserves a volume form $\mu$ and $H_{*}\mu = \mu_{\Gamma}$ (see Lemma \ref{L:ImprovedRegularity3}).

If $\alpha$ is assumed to preserve a volume form $\mu$ such that $H_{*}\mu = \mu_{\Gamma}$ then a result similar to Theorem \ref{Thm:ImprovedRegularity1} still holds, without accessibility. We give a brief sketch of the proof. Let $H:G\to G$ be a lift of the conjugacy. By the argument below $H$ is uniformly $C^{\infty}$ along $W^{\sigma}$ and $D_{x}H:E^{\sigma}(x)\to E_{0}^{\sigma}$ is invertible at each $x$ for $\sigma = s,u$. Define
\begin{align}
\Hat{H}:G\xrightarrow{H}G\to G/G^{cs}
\end{align}
then the fibers of $\Hat{H}$ coincides with $\Hat{W}^{cs}$. In particular, $\Hat{H}$ is uniformly $C^{\infty}$ along $\Hat{W}^{cs}$. Moreover, $H$ is uniformly $C^{\infty}$ along $W^{u}$ so $\Hat{H}$ is uniformly $C^{\infty}$ along $W^{u}$. By Journé's lemma \cite{Journe1988} $\Hat{H}$ is smooth. The map $D_{x}\Hat{H}|_{E^{u}(x)}$ is invertible, so $\Hat{H}$ is a submersion. It follows that the leaves of $\Hat{H}$ form a $C^{\infty}-$foliation, so $W^{cs}$ is a $C^{\infty}-$foliation. Similarly, $W^{cu}$ is a $C^{\infty}-$foliation. Once we know that $W^{c}$ is a $C^{\infty}-$foliation the arguments in \cite{Rodriguez-HertzWang2014,FisherKalininSpatzier2013}, to obtain regularity of $H$ along $W^{s}$ and $W^{u}$, can also be used along $W^{c}$ (note that, a priori, the assumptions of \cite[Theorem A.1]{Rodriguez-HertzWang2014} are not satisfied for $W^{c}$).

\subsection{Dynamical coherence and regularity of center leaves}

We begin by proving that $f$ is dynamically coherent with $W^{cs}$, $W^{cu}$ and $W^{c}$ all uniquely integrable. We use this without mention in the remainder. 
\begin{lemma}
The map $f$ is dynamically coherent with $W^{cs}$, $W^{cu}$ and $W^{c}$ all uniquely integrable with uniformly $C^{\infty}$ leaves.
\end{lemma}
\begin{proof}
Since $f$ is conjugated to $f_{0} = \alpha_{0}^{\mathbf{n}_{0}}$, and $f_{0}$ is uniformly subexponetial along its center $E_{0}^{c}$, $f$ is also uniformly subexponential along its center $E^{c}$. Indeed, for any $f-$invariant measure $\nu$ the stable Pesin manifold $W_{\nu}^{s}(x)$ maps into the stable manifold, $W_{0}^{s}(H(x))$, of $f_{0}$ under $H$. Since $H$ is invertible $\dim W_{\nu}^{s}(x)\leq\dim W_{0}^{s}(H(x))$. It follows that no $\nu-$exponents are negative along $E^{c}$. By exchanging $f$ for $f^{-1}$, no $\nu-$exponent is positive along $E^{c}$. By \cite{Schreiber1998} $f$ is uniformly subexponential along $E^{c}$. Fix lifts $F$, $F_{0}$, and $H$ to $G$ such that $H(Fx) = F_{0}(H(x))$. If $\gamma:I\to G$ is a $C^{1}-$curve tangent to $E^{cs}$ then the length $|F^{n}\circ\gamma|$ satisfies
\begin{align}
|F^{n}\circ\gamma|\leq C_{\varepsilon}e^{n\varepsilon},
\end{align}
for any $\varepsilon > 0$. The conjugacy $H:G\to G$ can be written $H(x) = xh(x)^{-1}$ with $h:G\to G$ being $\Gamma-$invariant. We estimate
\begin{align*}
\intd(H(x),H(y)) = &\intd(xh(x)^{-1},yh(y)^{-1})\leq \\ \leq &
\intd(x,xh(x)^{-1}) + \intd(y,yh(y)^{-1}) + \intd(x,y) = \\ &
\intd(e,h(x)^{-1}) + \intd(e,h(y)^{-1}) + \intd(x,y)\leq\intd(x,y) + C
\end{align*}
for some uniform $C$. With $x = H^{-1}(z)$ and $y = H^{-1}(w)$ we obtain 
\begin{align*}
\intd(H^{-1}(z),H^{-1}(w))\geq\intd(z,w) - C.
\end{align*}
Now, $\intd(F^{n}(\gamma(0)),F^{n}(\gamma(1)))\leq C_{\varepsilon}e^{n\varepsilon}$, so
\begin{align*}
C_{\varepsilon}e^{n\varepsilon}\geq & \intd(F^{n}(\gamma(0)),F^{n}(\gamma(1))) = \intd(H^{-1}HF^{n}(\gamma(0)),H^{-1}HF^{n}(\gamma(1)))\geq \\ &
\intd(HF^{n}(\gamma(0)),HF^{n}(\gamma(1))) - C = \intd(F_{0}^{n}H(\gamma(0)),F_{0}^{n}(H(\gamma(1)))) - C.
\end{align*}
With $\varepsilon$ sufficiently small it follows that $H(\gamma(1))\in W_{0}^{cs}(H(\gamma(1)))$. We conclude that $E^{cs}$ is uniquely integrable with leaves given by 
\begin{align}
W^{cs}(x) = H^{-1}(W_{0}^{cs}(H(x))).
\end{align}
Similarly, $E^{cu}$ is uniquely integrable with $W^{cu}(x) = H^{-1}(W_{0}^{cu}(H(x)))$. Since $f$ is uniformly subexponential along $E^{c}$ each foliation $W^{cs}$, $W^{cu}$ and $W^{c} = W^{cs}\cap W^{cu}$ have uniformly $C^{\infty}$ leaves, see \cite[Theorem 7]{DamjanovicWilkinsonXu2021} (or \cite{HirschShubPugh1977}).
\end{proof}

\subsection{Volume preservation and smoothness along the center direction}

To apply arguments using exponential mixing we need to show that $\alpha$ preserves a smooth volume form. We assume that the action $\alpha$ is accessible so we show that $H$ is $C^{\infty}$ along $W^{c}$ without using the exponential mixing argument from \cite{FisherKalininSpatzier2013} and instead rely on results from \cite{Wilkinson2008}.
\begin{lemma}\label{L:ImprovedRegularity3}
Let $\alpha$ be as in Theorem \ref{Thm:ImprovedRegularity1}, then $\alpha$ preserve a smooth volume form $\mu$. Moreover, the conjugacy $H:X_{\Gamma}\to X_{\Gamma}$ is volume preserving in the sense that $H_{*}\mu = \mu_{\Gamma}$ where $\mu_{\Gamma}$ is the Haar measure on $X_{\Gamma}$.
\end{lemma}
\begin{remark}
By Moser's trick, there is no loss of generality if we assume that $\mu = \mu_{\Gamma}$.
\end{remark}
\begin{proof}
Existence of an invariant Hölder continuous volume form, $\mu$, follows from \cite[Theorem 1.5]{GorodnikSpatzier2015}. We assume accessibility, so smoothness of $\mu$ follows from \cite[Theorem A, case (IV)]{Wilkinson2008}. As in \cite[Proposition 2.4]{FisherKalininSpatzier2011}, $H_{*}\mu$ is a measure of maximal entropy for the $\alpha_{0}-$action. From \cite{Wang2018} it follows that $H_{*}\mu = \mu_{\Gamma}$.
\end{proof}
We show that $H$ is uniformly $C^{\infty}$ along $W^{c}$ following \cite{Wilkinson2008}.
\begin{lemma}\label{L:ImprovedRegularity4}
The restriction $H:W^{c}(x)\to W_{0}^{c}(H(x))$ is uniformly $C^{\infty}$.
\end{lemma}
\begin{remark}\label{Rmk:JointCenter}
Using Lemma \ref{L:ImprovedRegularity4} we may assume, without loss of generality, that the center $E_{0}^{c}$ of $f_{0}$ coincides with the joint center of $\alpha$.
\end{remark}
\begin{proof}
Let $M = W^{c}(x)\times W_{0}^{c}(H(x))$ and let $N\subset M$ be the graph of $H$:
\begin{align}
N = \{(y,H(y))\text{ : }y\in W^{c}(x)\}.
\end{align}
Given any two $z,w\in W^{c}(x)$ we fix a $su-$path $\gamma$ from $z$ to $w$, such a path always exists by accessibility. Denote by $h_{z,w}^{\gamma}:W^{c}(x)\to W^{c}(x)$ the composition of holonomy maps along $\gamma$ ($z,w\in W^{c}(x)$ so $W^{c}(z) = W^{c}(w) = W^{c}(x)$). Since $f$ is uniformly subexponential along $W^{c}$, $f$ is $\infty-$bunching so $h_{z,w}^{\gamma}$ is $C^{\infty}$ \cite{PughShubWilkinson1997}. Since $H$ map $W^{\sigma}$, $\sigma = s,u$, onto $W_{0}^{\sigma}$ we have $Hh_{z,w}^{\gamma} = h_{H(z),H(w)}^{H\gamma,0}H$ where $h_{H(z),H(w)}^{H\gamma,0}:W_{0}^{c}(H(x))\to W_{0}^{c}(H(x))$ is the composition of holonomy maps along $H\gamma$. Define
\begin{align}
\Hat{h}_{z,w}:M\to M,\quad\Hat{h}_{z,w}(p,q) = \left(h_{z,w}^{\gamma}(p),h_{H(z),H(w)}^{H\gamma,0}(q)\right).
\end{align}
Since $Hh_{z,w}^{\gamma} = h_{H(z),H(w)}^{H\gamma,0}H$ we have $\Hat{h}_{z,w}(N) = N$. Moreover, $\Hat{h}_{z,w}(z,H(z)) = (w,H(w))$, $z$ and $w$ were arbitrary, and $\Hat{h}_{z,w}$ is smooth so $N$ is $C^{\infty}-$homogeneous (see \cite{Wilkinson2008} for definitions). By \cite[Corollary 1.3]{Wilkinson2008} $N$ is a $C^{\infty}$ submanifold. The graph of $H:W^{c}(x)\to W^{c}(H(x))$ is $C^{\infty}$ and it follows that $H:W^{c}(x)\to W_{0}^{c}(H(x))$ is also $C^{\infty}$. Finally, $H$ is uniformly $C^{\infty}$ along $W^{c}$ since it intertwines the holonomies of $f$ with the holonomies of $f_{0}$.
\end{proof}

\subsection{Smoothness of coarse components along the stable foliation}

Let $\mathfrak{g} = E_{0}^{s}\oplus E_{0}^{c}\oplus E_{0}^{u}$ be the splitting of $\mathfrak{g}$ with respect to $f_{0}$. Denote by $G^{\sigma}$, $\sigma = s,c,u,cs,cu$, the subgroup associated to $E_{0}^{\sigma}$. We have a $C^{\infty}-$diffeomorphism $G^{s}\times G^{c}\times G^{u}\to G$ defined by $(g^{s},g^{c},g^{u})\mapsto g_{s}g_{c}g_{u}$. Write $H(x) = xh(x)^{-1}$ with $h:X_{\Gamma}\to G$ satisfying $h(\gamma x) = h(x)$ for all $\gamma\in\Gamma$. We decompose $h(x)$ with respect to the map $G^{s}\times G^{c}\times G^{u}\to G$
\begin{align}
h(x) = h_{s}(x)h_{c}(x)h_{u}(x).
\end{align}
It is immediate that each $h_{\sigma}$ is Hölder. Given a coarse exponent $[\chi]$ we denote by $G^{[\chi]}$ the subgroup associated to $E_{0}^{[\chi]}$.

Let $[\chi]$ be a coarse exponent along $E_{0}^{s}$. We decompose $h_{s}(x)$ further as
\begin{align}\label{Eq:StableDecompositionOfConj}
h_{s}(x) = h_{ss}(x)h_{[\chi]}(x)
\end{align}
where $h_{[\chi]}(x)$ is the component of $h_{s}(x)$ along $G^{[\chi]}$ and $h_{ss}(x)$ is the component of $h_{s}(x)$ along the complementary group $G^{ss}$ (see \cite[Lemma 3.1]{Rodriguez-HertzWang2014}). The following is proved in \cite[Section 3]{Rodriguez-HertzWang2014}.
\begin{lemma}\label{L:ImprovedRegularity5}
Let $\alpha$ be as in Theorem \ref{Thm:ImprovedRegularity1} with $\alpha^{\mathbf{n}_{0}} = f$ partially hyperbolic and $\ker[\chi]$, $[\chi](\mathbf{n}_{0}) < 0$, a wall of the Weyl chamber that contains $\mathbf{n}_{0}$. The map $h_{[\chi]}(x)$ in the Equation \ref{Eq:StableDecompositionOfConj} is $C^{\infty}$ along $W^{s}$, with all derivatives along $W^{s}$ uniformly Hölder.
\end{lemma}
\begin{proof}
The proof follows as in \cite[Section 3]{Rodriguez-HertzWang2014} once we note that we do not need $f$ to be Anosov. Indeed, once we restrict to the foliation $W^{s}$ the argument only requires the action $\alpha$ to be exponentially mixing with respect to volume. Exponential mixing with respect to volume follows from Lemma \ref{L:ImprovedRegularity3} and \cite{GorodnikSpatzier2015}.
\end{proof}

\subsection{New partially hyperbolic elements: passing the chamber wall}

Let $\mathbf{n}_{0}$ be in a Weyl chamber $\mathcal{C}$ and let $\ker[\chi]$ be a chamber wall for $\mathcal{C}$. Now we start the work of passing the Weyl chamber wall $\ker[\chi]$ by constructing a partially hyperbolic element in the chamber $\mathcal{C}'$ adjacent to $\mathcal{C}$ through $\ker[\chi]$. We initially follow \cite{Rodriguez-HertzWang2014}, but change the argument from Section $4$ in \cite{Rodriguez-HertzWang2014}, by not relying on smooth ergodic theory. If $x\in X_{\Gamma}$ then the map $H:W^{s}(x)\to W_{0}^{s}(H(x)) = H(x)G^{s}$ is a homeomorphism. Define
\begin{align}
H_{s,x}:W^{s}(x)\to G^{s},\quad H(y) = H(x)\left(H_{s,x}(y)\right)^{-1}.
\end{align}
If $\alpha_{0}^{\mathbf{n}}(x) = \rho^{\mathbf{n}}(x)\eta_{\mathbf{n}}^{-1}$ with $\rho^{\mathbf{n}}\in{\rm Aut}(X_{\Gamma})$ and $\eta_{\mathbf{n}}\in G$ then $H_{s,x}$ satisfy
\begin{align*}
& H_{s,\alpha^{\mathbf{n}}x}(\alpha^{\mathbf{n}}y) = \eta_{\mathbf{n}}\rho^{\mathbf{n}}\left(H_{s,x}(y)\right)\eta_{\mathbf{n}}^{-1}.
\end{align*}
For $y\in W^{s}(x)$ we write $y = xg_{x}(y)^{-1}$ where $(x,y)\mapsto g_{x}(y)$ is chosen continuously and such that $g_{x}(x) = e$. Then $g_{x}:W^{s}(x)\to G^{s}$ is $C^{\infty}$. With this notation
\begin{align*}
xg_{x}(y)^{-1}h(y)^{-1} = yh(y)^{-1} = H(y) = H(x)H_{s,x}(y)^{-1} = xh(x)^{-1}H_{s,x}(y)^{-1}
\end{align*}
or
\begin{align*}
H_{s,x}(y) = & h(y)g_{x}(y)h(x)^{-1} = \\ &
h_{s}(y)h_{cu}(y)\left[(h(x)g_{x}(y)^{-1})_{s}(h(x)g_{x}(y)^{-1})_{cu}\right]^{-1} = \\ &
h_{s}(y)\left[h_{cu}(y)\left(h(x)g_{x}(y)^{-1}\right)_{cu}^{-1}\right]\left(h(x)g_{x}(y)^{-1}\right)_{s}^{-1}.
\end{align*}
If we multiply on the left by $h_{s}(y)^{-1}$ and on the right by $(h(x)g_{x}(y))_{s}$, then the right-hand side of the equality lie in $G^{cu}$, but the left-hand side lies in $G^{s}$. It follows that both sides of the equality are identity, so for $y\in W^{s}(x)$ we have
\begin{align}\label{Eq:ImprovedRegularity1}
& h_{cu}(y) = \left(h(x)g_{x}(y)^{-1}\right)_{cu}, \\
& H_{s,x}(y) = h_{s}(y)\left(h(x)g_{x}(y)^{-1}\right)_{s}^{-1}.
\end{align}
Using the fact that the map $a\mapsto a_{\sigma}$, $\sigma = s,cu$, and the map $g_{x}:W^{s}(x)\to G$ are both smooth the following is immediate, see also \cite[Corollary 3.14]{Rodriguez-HertzWang2014}.
\begin{lemma}\label{L:ImprovedRegularity6}
The map $h_{cu}$ is uniformly $C^{\infty}$ along $W^{s}$.
\end{lemma}
Define $H_{s,x}(y) = (H_{s,x}(y))_{ss}H_{s,x}^{[\chi]}(y)$ with $H_{s,x}^{[\chi]}(y)\in G^{[\chi]}$. From our definitions it is immediate
\begin{align*}
\{y\in W^{s}(x)\text{ : }H_{s,x}^{[\chi]}(y) = e\} = H^{-1}(W_{0}^{ss}(H(x))),\quad W_{0}^{ss}(y) = yG^{ss}.
\end{align*}
So if we prove that $H_{s,x}^{[\chi]}$ is a (local) $C^{\infty}$ submersion for every $x$, then the fibers of $H_{s,x}^{[\chi]}$ defines a smooth foliation (within $W^{s}$). This shows that the, a priori only Hölder, foliation $W^{ss}(x) = H^{-1}(W_{0}^{ss}(H(x)))$ is a Hölder foliation with uniformly smooth leaves.
\begin{lemma}\label{L:ImprovedRegularity7}
The map $H_{s,x}^{[\chi]}:W^{s}(x)\to G^{[\chi]}$ is uniformly $C^{\infty}$.
\end{lemma}
\begin{proof}
Since $G^{ss}$ is normal in $G^{s}$ \cite[Lemma 3.1]{Rodriguez-HertzWang2014} we have
\begin{align*}
H_{s,x}(y) = & h_{s}(y)\left(h(x)g_{x}(y)^{-1}\right)_{s}^{-1} = \\ &
h_{ss}(y)h_{[\chi]}(y)\left[\left(h(x)g_{x}(y)^{-1}\right)_{s}^{-1}\right]_{ss}\left[\left(h(x)g_{x}(y)^{-1}\right)_{s}^{-1}\right]_{[\chi]} = \\ &
a_{ss}\cdot h_{[\chi]}(y)\left[\left(h(x)g_{x}(y)^{-1}\right)_{s}^{-1}\right]_{[\chi]}
\end{align*}
for some $a_{ss}\in G^{ss}$. That is, we obtain the formula
\begin{align}
H_{s,x}^{[\chi]}(y) = h_{[\chi]}(y)\left[\left(h(x)g_{x}(y)^{-1}\right)_{s}^{-1}\right]_{[\chi]}
\end{align}
so $H_{s,x}^{[\chi]}$ is uniformly $C^{\infty}$ along $W^{s}$ since $h_{[\chi]}$ and $g_{x}$ are (Lemma \ref{L:ImprovedRegularity5}).
\end{proof}
\begin{lemma}\label{L:ImprovedRegularity8}
The map $D_{x}H_{s,x}^{[\chi]}:E^{s}(x)\to E_{0}^{[\chi]}$ is surjective at every $x\in X_{\Gamma}$. In particular, the foliation $W^{ss}$ has uniformly $C^{\infty}$ leaves.
\end{lemma}
\begin{proof}
Denote by $K = \{D_{x}H_{s,x}^{[\chi]}\text{ : is not surjective}\}$, then $K$ is compact and $\alpha-$invariant. Our goal is to show that $K = \emptyset$. Since $H:W^{s}(x)\to W_{0}^{s}(H(x))$ is surjective it follows by Sard's theorem that $K\neq X_{\Gamma}$. We will prove that if $K$ is non-empty then $K$ contains a $W^{ss}-$leaf, which is a contradiction since every $W^{ss}-$leaf is dense (Lemma \ref{L:ImprovedRegularity1}).

In the remainder, fix some background metric $\langle\cdot,\cdot\rangle$ and calculate determinants with respect to the top form induced by $\langle\cdot,\cdot\rangle$. Assume for contradiction that $K\neq\emptyset$. Fix $\mathbf{n}$ close to the kernel $\ker[\chi]$ (we specify later how close) such that $\alpha^{\mathbf{n}}$ contract $W^{ss} := H^{-1}(H(x)G^{ss})$ and $W^{[\chi]}$. Fix $y\in W^{ss}(x)$ and some subspace $V_{0}\subset E^{s}(y)$ of dimension $\dim(E_{0}^{[\chi]})$. We also define $V_{n} = D_{y}\alpha^{n\mathbf{n}}V_{0}$. The relation $H_{s,\alpha^{\mathbf{m}}x}^{[\chi]}(\alpha^{\mathbf{m}}y) = \eta_{\mathbf{m}}\rho^{\mathbf{m}}\left(H_{s,x}^{[\chi]}(y)\right)\eta_{\mathbf{m}}^{-1}$ implies
\begin{align}\label{Eq:FormulaForStablePinch1}
\det(D_{y}H_{s,y}^{[\chi]}|_{V_{0}}) = \det(\rho^{-n\mathbf{n}}|_{E_{0}^{[\chi]}})\det(D_{\alpha^{n\mathbf{n}}y}H_{s,\alpha^{n\mathbf{n}}y}^{[\chi]}|_{V_{n}})\det(D_{y}\alpha^{n\mathbf{n}}|_{V_{0}}).
\end{align}
If $\chi$ is a representative of $[\chi]$ then $\det(D\rho^{-n\mathbf{n}})\leq Ce^{-rn\chi(\mathbf{n})}$ for some uniform $r > 0$ depending on the dimension of $E_{0}^{[\chi]}$. If $\ell = \dim(E_{0}^{[\chi]})$ and ${\rm Gr}_{\ell}(E^{s})$ is the $\ell-$grassmannian bundle of $E^{s}$ then ${\rm Gr}_{\ell}(E^{s})\ni(V,x)\mapsto\det(D_{x}H_{s,x}^{[\chi]}|_{V})$ is uniformly smooth along $W^{s}$. For $x\in K$ we have $\det(D_{x}H_{s,x}^{[\chi]}|_{V}) = 0$ since $D_{x}H_{s,x}^{[\chi]}$ is not surjective. It follows that there is some constant $C$ such that
\begin{align}
|\det(D_{y}H_{s,y}^{[\chi]}|_{V})|\leq C\intd_{s}(y,K)
\end{align}
for $y\in W^{s}(K)$. If $y\in W^{ss}(x)$, then
\begin{align}\label{Eq:FormulaForStablePinch2}
\det(D_{\alpha^{n\mathbf{n}}y}H_{\alpha^{n\mathbf{n}}y}^{[\chi]}|_{V_{n}})\leq &  C\intd(\alpha^{n\mathbf{n}}y,K)\leq C\intd(\alpha^{n\mathbf{n}}y,\alpha^{n\mathbf{n}}x)\leq C'e^{-\lambda n}
\end{align}
where $\lambda > 0$ can be chosen independently of $\mathbf{n}$ since we let $\mathbf{n}$ be close to the kernel of $[\chi]$. Equations \ref{Eq:FormulaForStablePinch1} and \ref{Eq:FormulaForStablePinch2} implies
\begin{align*}
|\det(D_{y}H_{s,y}^{[\chi]}|_{V_{0}})|\leq Ce^{-n(\lambda + r\chi(\mathbf{n}))}|\det(D_{y}\alpha^{n\mathbf{n}}|_{V_{0}})|.
\end{align*}
Since $\alpha_{0}^{\mathbf{n}}$ contract the foliations $W^{ss}$ and $W^{[\chi]}$ we have a uniform bound $|\det(D_{y}\alpha^{n\mathbf{n}}|_{V_{0}})|\leq C$. We obtain an estimate
\begin{align*}
|\det(D_{y}H_{s,y}^{[\chi]}|_{V_{0}})|\leq Ce^{-n(\lambda + r\chi(\mathbf{n}))}.
\end{align*}
With $\mathbf{n}$ be sufficiently close to $\ker[\chi]$ we have $\lambda + r\chi(\mathbf{n}) > 0$, so if $n\to\infty$ then
\begin{align*}
\det(D_{y}H_{s,y}^{[\chi]}|_{V_{0}}) = 0.
\end{align*}
The point $y$ was arbitrary so $W^{ss}(x)\subset K$ for $x\in K$. The foliation $W^{ss}$ is minimal (Lemma \ref{L:ImprovedRegularity1}) and $K$ is closed so $K = X_{\Gamma}$. This contradicts Sard's theorem.
\end{proof}
To show that $E^{[\chi]}$ exists as a continuous bundle we will apply a linear graph transform argument. Denote by $E^{ss}(x)$ the $D\alpha^{\mathbf{n}}-$invariant, continuous subbundle tangent to $W^{ss}(x)$. Let $F(x)$ be any continuous bundle that is complementary to $E^{ss}$ within $E^{s}$. Fix $\mathbf{n}\in\mathbb{Z}^{k}$ such that $[\chi](\mathbf{n}) < 0$ and $\alpha^{\mathbf{n}}$ expands $W^{ss}$ (that is, we chose $\mathbf{n}$ such that $-\mathbf{n}$ have passed the chamber wall $\ker[\chi]$ from the chamber that contains $\mathbf{n}_{0}$). Write $g = \alpha^{\mathbf{n}}$. With respect to the splitting $E^{s} = F\oplus E^{ss}$ we write
\begin{align*}
D_{x}g(u,v) = \left(A(x)u,B(x)v + C(x)u\right)
\end{align*}
where $A(x)u$ is $D_{x}g(u)$ projected onto $F(gx)$ and $C(x)$ is $D_{x}g(u)$ projected onto $E^{ss}(gx)$.
\begin{lemma}\label{L:ImprovedRegularity9}
We have $\norm{A}_{C^{0}} < 1$.
\end{lemma}
\begin{proof}
Note that $D_{x}H_{s,x}^{[\chi]}D_{x}g = D_{e}\alpha_{0}^{\mathbf{n}}D_{x}H_{s,x}^{[\chi]}$. There is $\mu < 1$, depending on our choice of $\mathbf{n}\in\mathbb{Z}^{k}$, such that for any $u\in F(x)$ we have
\begin{align}\label{Eq:NormEstimateInGraphTransform1}
\norm{D_{gx}H_{s,gx}^{[\chi]}\left(A(x)u + B(x)u\right)} = \norm{D_{e}\alpha_{0}^{\mathbf{n}}D_{x}H_{s,x}^{[\chi]}(u)}\leq\mu\norm{D_{x}H_{s,x}^{[\chi]}(u)}.
\end{align}
Since $B(x)u\in E^{ss}(gx)$ and $\ker D_{x}H_{s,x}^{[\chi]} = E^{ss}(x)$ (the fibers of $H_{s,x}^{[\chi]}$ are the foliation $W^{ss}(x)$) it follows that $D_{gx}H_{s,gx}^{[\chi]}B(x) = 0$. That is, Equation \ref{Eq:NormEstimateInGraphTransform1} simplifies
\begin{align}\label{Eq:NormEstimateInGraphTransform2}
\norm{D_{gx}H_{s,gx}^{[\chi]}A(x)u}\leq\mu\norm{D_{x}H_{s,x}^{[\chi]}(u)}.
\end{align}
Since $D_{x}H_{s,x}^{[\chi]}:F(x)\to E_{0}^{[\chi]}$ is an isomorphism (of vector bundles after identifying $E_{0}^{[\chi]}$ with the trivial bundle $X_{\Gamma}\times E_{0}^{[\chi]}\to X_{\Gamma}$) the lemma follows (after either changing the norm used or exchanging $\mathbf{n}$ for $N\mathbf{n}$ with $N$ sufficiently large).
\end{proof}
\begin{lemma}\label{L:ImprovedRegularity10}
The map $T:\Gamma^{0}({\rm Hom}(F,E^{ss}))\to\Gamma^{0}({\rm Hom}(F,E^{ss}))$ defined by
\begin{align*}
(TP)(x) = B(x)^{-1}\left(P(gx)A(x) - C(x)\right)
\end{align*}
has a unique fixed point.
\end{lemma}
\begin{proof}
Since $\norm{A}_{C^{0}},\norm{B^{-1}}_{C^{0}} < 1$ the lemma follows from Banach's fixed point theorem.
\end{proof}
\begin{lemma}\label{L:ImprovedRegularity11}
There exists an $\alpha-$invariant continuous subbundle $E^{[\chi]}\subset E^{s}$ such that $E^{s} = E^{[\chi]}\oplus E^{ss}$.
\end{lemma}
\begin{proof}
Let $P\in\Gamma^{0}({\rm Hom}(F,E^{ss}))$ be the unique $T-$fixed point from Lemma \ref{L:ImprovedRegularity10}. Define $E^{[\chi]}(x) := {\rm Graph}(P(x)) = \{(u,P(x)u)\text{ : }u\in F(x)\}$. It is immediate that $E^{[\chi]}\oplus E^{ss} = E^{s}$. Given $u\in F(x)$
\begin{align*}
D_{x}g(u,P(x)u) = & \left(A(x)u,B(x)P(x)u + C(x)u\right) = \\ &
\left(A(x)u,P(gx)A(x)u - C(x)u + C(x)u\right) = \\ &
\left(A(x)u,P(gx)A(x)u\right)\in{\rm Graph}(P(gx))
\end{align*}
so $D_{x}gE^{[\chi]}(x)\subset E^{[\chi]}(gx)$. Or $D_{x}gE^{[\chi]}(x) = E^{[\chi]}(gx)$ since $D_{x}g$ is invertible. That $E^{[\chi]}(x)$ is $\alpha^{\mathbf{m}}-$invariant for all $\mathbf{m}\in\mathbb{Z}^{k}$ follows by applying the graph transform of $\alpha^{\mathbf{m}}$ on the element $P\in\Gamma^{0}({\rm Hom}(F,E^{ss}))$. This defines a $T-$fixed point, since $g$ commute with $\alpha^{\mathbf{m}}$, and the $T-$fixed point is unique (Lemma \ref{L:ImprovedRegularity10}) so the $\alpha^{\mathbf{m}}-$graph transform of $P$ is $P$. Equivalently, $E^{[\chi]}$ is $D_{x}\alpha^{\mathbf{m}}-$invariant.
\end{proof}
Since $D_{x}H_{s,x}^{[\chi]}:E^{[\chi]}(x)\to E_{0}^{[\chi]}$ conjugates $D_{x}\alpha^{\mathbf{n}}$ to ${\rm Ad}(\eta_{\mathbf{n}})D\rho^{\mathbf{n}}$ the following lemma follows by induction.
\begin{lemma}\label{L:ImprovedRegularity12}
Every element $\mathbf{n}\in\mathbb{Z}^{k}\setminus0$ defines a partially hyperbolic $\alpha^{\mathbf{n}}:X_{\Gamma}\to X_{\Gamma}$ where the center of $\alpha^{\mathbf{n}}$ has the same dimension as the center of $\alpha_{0}^{\mathbf{n}}$.
\end{lemma}
\begin{remark}
Lemma \ref{L:ImprovedRegularity12} lets us pass the Weyl chamber containing $\mathbf{n}_{0}$ and produce new partially hyperbolic elements in adjacent chambers. 
\end{remark}

\subsection{Finishing the proof of Theorem \ref{Thm:ImprovedRegularity1}}

We finish the proof of Theorem \ref{Thm:ImprovedRegularity1} by showing that $H$ is a $C^{\infty}$ diffeomorphism. We begin by proving that $h_{\sigma}$, $\sigma = s,u$, are $C^{\infty}$.
\begin{lemma}\label{L:ImprovedRegularity13}
The maps $h_{\sigma}:X_{\Gamma}\to G^{\sigma}$, $\sigma = s,u$, are $C^{\infty}$.
\end{lemma}
\begin{proof}
It suffices to consider the case $\sigma = s$. By Lemma \ref{L:ImprovedRegularity6} and Journé's lemma \cite{Journe1988} it suffices to show that $h_{s}$ is uniformly smooth along $W^{s}$. Number all coarse exponents along $W^{s}$ by $[\chi_{1}],...,[\chi_{N}]$. The function $h_{s}$ can be decomposed with respect to the map $G^{[\chi_{1}]}\times...\times G^{[\chi_{N}]}\to G^{s}$ (see \cite[Lemma 3.2]{Rodriguez-HertzWang2014}) as
\begin{align}\label{Eq:CoarseDecompositionOfConjugacy}
h_{s}(x) = h_{[\chi_{1}]}(x)...h_{[\chi_{N}]}(x).
\end{align}
Using Lemmas \ref{L:ImprovedRegularity12} and \ref{L:ImprovedRegularity5} it is immediate that $h_{[\chi_{N}]}(x)$ is $C^{\infty}$. Moreover, we use Lemmas \ref{L:ImprovedRegularity12} and \ref{L:ImprovedRegularity5} to show that $h_{[\chi_{j}]}^{(2)}(x) = h_{[\chi_{j}]}(x)/G^{(2)}$ is smooth for all $j = 1,...,N$. If we define $h_{[\chi_{j}]}^{(i)}(x) = h_{[\chi_{j}]}(x)/G^{(i)}$ and assume that $h_{[\chi_{j}]}^{(i-1)}(x)$ is smooth, then we can change the order of products in \ref{Eq:CoarseDecompositionOfConjugacy} modulo a polynomial in $h_{[\chi_{j}]}^{(i-1)}$. Applying Lemmas \ref{L:ImprovedRegularity12} and \ref{L:ImprovedRegularity5} once more, $h_{[\chi_{j}]}^{(i)}$ is $C^{\infty}$ for each $j = 1,...,N$. For $i$ large enough, using that $G$ is nilpotent, it follows that $h_{[\chi_{j}]}$ is $C^{\infty}$ along $W^{s}$.
\end{proof}
We are now ready to prove Theorem \ref{Thm:ImprovedRegularity1}.
\begin{proof}[Proof of Theorem \ref{Thm:ImprovedRegularity1}]
If we show that $H$ is $C^{\infty}$ then $H$ is automatically a diffeomorphism since the Jacobian can never vanish, this would contradict volume preservation of $\alpha$ (Lemma \ref{L:ImprovedRegularity3}).

By Lemma \ref{L:ImprovedRegularity4} the map $H$ is uniformly $C^{\infty}$ along $W^{c}$. From equation \ref{Eq:ImprovedRegularity1} we have $H_{s,x}(y) = h_{s}(y)\left(h(x)g_{x}(y)^{-1}\right)_{s}^{-1}$. Since $g_{x}(y)$ and $h_{s}(y)$ are both uniformly $C^{\infty}$ along $W^{s}$ byLemma \ref{L:ImprovedRegularity13}, it follows that $H_{s,x}$ is uniformly $C^{\infty}$ along $W^{s}$, so $H$ is uniformly $C^{\infty}$ along $W^{s}$. Similarly $H$ is uniformly $C^{\infty}$ along $W^{u}$. By using Journé's lemma along $W^{c}$ and $W^{s}$ it follows that $H$ is uniformly $C^{\infty}$ along $W^{cs}$. Using Journé's lemma once more along $W^{u}$ and $W^{cs}$ it follows that $H$ is $C^{\infty}$.
\end{proof}

\section{Proof of main theorems}\label{Sec:ProofOfMainTheorems}

We are now ready to complete the proofs of Theorems \ref{Thm:MainTheorem1} and \ref{Thm:MainTheorem2}.
\begin{proof}[Proof of Theorem \ref{Thm:MainTheorem1}]
By Theorem \ref{Thm:ExistenceBiHölderConj} the action $\alpha$ is bi-Hölder conjugated to some affine action $\alpha_{0}$. We produce the conjugacy $H$ for the special element $f = \alpha^{\mathbf{n}_{0}}$, but it conjugates the full action into an affine action. This is immediate from the construction of $H$, but also follows from an argument as in \cite{AdlerPalais1965}. Theorem \ref{Thm:MainTheorem1} now follows from Theorem \ref{Thm:ImprovedRegularity1}.
\end{proof}
\begin{proof}[Proof of Theorem \ref{Thm:MainTheorem2}]
Let $f_{0}\in{\rm Aff}(X_{\Gamma})$ and $f\in{\rm Diff}^{\infty}(X_{\Gamma})$ be $C^{1}-$close to $f_{0}$. Write $L_{su}:\mathbb{T}^{d}\to\mathbb{T}^{d}$ for the induced hyperbolic automorphism. Note that $f$ satisfy assumptions $(i)$ and $(ii)$ of Theorem \ref{Thm:MainTheorem1} since it is close to $f_{0}$ (see also \cite[Lemma A.2]{LeeSandfeldt2024}). In particular, we obtain $\Phi$ from Theorem \ref{Thm:PropertiesOfPHdiffeos}. Denote by $Z^{\infty}(f)$ the $C^{\infty}-$centralizer of $f$ and let $Z_{c}^{\infty}(f)\subset Z^{\infty}(f)$ be the center fixing, normal subgroup of $Z^{\infty}(f)$ from Equation \ref{Eq:CenterFixingCentralizer}. Define the quotient 
\begin{align}
Z_{su}^{\infty}(f) = Z^{\infty}(f)/Z_{c}^{\infty}(f).
\end{align}
If $g\in Z_{c}^{\infty}(f)$ then $\Phi(gx) = \Phi(x)$ so the induced map on $H_{1}X_{\Gamma}$ satisfy $\Phi_{*}g_{*} = \Phi_{*}$. Since $\Phi_{*}:H_{1}X_{\Gamma}\to H_{1}\mathbb{T}^{d}$ is an isomorphism it follows that $g_{*} = {\rm id}$. Conversely, if $g_{*} = {\rm id}$ then $\Phi(gx) = \Phi(x)$ (Theorem \ref{Thm:PropertiesOfPHdiffeos}) so $g\in Z_{c}^{\infty}(f)$. It follows that each non-trivial $g\in Z_{su}^{\infty}(f)$ represent an element of $Z^{\infty}(f)$ that project onto a non-trivial automorphism on $\mathbb{T}^{d}$. In particular, if ${\rm rank}(Z_{su}^{\infty}(f)) > 1$ then the image of $Z^{\infty}(f)$ in $Z_{\rm Aut}(L_{su})$ contains a subgroup isomorphic to $\mathbb{Z}^{2}$. Irreducibility of $L_{su}$ implies that this $\mathbb{Z}^{2}-$subgroup in $Z^{\infty}(f)$ is higher rank, so the action of $Z^{\infty}(f)$ is $C^{\infty}-$conjugate to some affine action by Theorem \ref{Thm:MainTheorem1}. If ${\rm rank}(Z_{su}^{\infty}(f)) = 1$ and $\# Z_{c}^{\infty}(f) = \infty$ then \cite[Corollary 18]{DamjanovicWilkinsonXu2023} implies case $(ii)$ of Theorem \ref{Thm:MainTheorem2}. Finally, if ${\rm rank}(Z_{su}^{\infty}(f)) = 1$ and $\# Z_{c}^{\infty}(f) < \infty$ then $Z^{\infty}(f)$ is virtually $\mathbb{Z}$ so case $(i)$ of Theorem \ref{Thm:MainTheorem2} holds.
\end{proof}

\newpage

\appendix

\section{Some algebraic lemmas}\label{Sec:Appendix}

In this appendix, we show some basic properties of higher rank algebraic actions on nilmanifolds. The first two lemmas, \ref{L:Background1} and \ref{L:Background2}, are stated in Section \ref{Sec:BackgroundNilmanifolds}.
\begin{proof}[Proof of Lemma \ref{L:Background1}]
If the conclusion does not hold then there is a decomposition $\mathfrak{g} = E_{0}^{s}\oplus E_{0}^{c}\oplus E_{0}^{u}$ so that every $\alpha_{0}^{\mathbf{n}}$ is subexponential along $E_{0}^{c}$ and either contract or expand $E_{0}^{s}$ and $E_{0}^{u}$. If the projected action on the base is rank$-1$, then the whole action is rank$-1$ (and has a rank$-1$ factor). For any $n\in\mathbb{N}$ there is $\varepsilon_{n} > 0$ such that if $L\in{\rm GL}(n,\mathbb{Z})$ satisfies that the eigenvalues of $L$ with modulus larger than one have a product bounded by $1+\varepsilon_{n}$, then $L$ have no eigenvalues of modulus larger than $1$.

Let $W\subset\mathbb{R}^{k}$ be the kernel of the unique pair of negatively proportional exponents of $\alpha_{0}$. If $\mathbf{n}\in\mathbb{Z}^{k}$ is sufficiently close to $W$ then all eigenvalues of $\alpha_{0}^{\mathbf{n}}\in{\rm GL}(n,\mathbb{Z})$ will be close to the unit circle, which implies that all eigenvalues of $\alpha_{0}^{\mathbf{n}}$ lie on the unit circle. It follows that any $\mathbf{n}\in\mathbb{Z}^{k}$ sufficiently close to $W$ lies in $W$. In particular
\begin{align}
{\rm rank}(\mathbb{Z}^{k}\cap W) = \dim(W) = k-1.
\end{align}
Elements in $W$ have all eigenvalues on the unit circle so after dropping to a finite index subgroup of $\mathbb{Z}^{k}$, we may assume that all eigenvalues are $1$. After taking a quotient to remove Jordan blocks the action of $W$ is trivial. Since $\dim(W) = k-1$ the action $\alpha_{0}$ factor through a $\mathbb{Z}-$action.
\end{proof}
\begin{proof}[Proof of Lemma \ref{L:Background2}]
The translation action of $V$ is minimal if and only if the induced translation action on the base is minimal, so we assume without loss of generality that $X_{\Gamma}$ is a torus. Let $W$ be the rational span of $V$, then $V$ acts minimally if and only if $W = \mathbb{R}^{d}$. If $W\neq\mathbb{R}^{d}$ then $\mathbb{T}^{\ell}\cong\mathbb{T}^{d}/W$ with $\ell\geq 1$ and $\alpha_{0}$ descend to $\mathbb{T}^{\ell}$. We have quotiened out all coarse exponents except for one negatively proportional pair, so the factor $\mathbb{T}^{\ell}$ has only one pair of negatively proportional exponents. By Lemma \ref{L:Background1} $\mathbb{T}^{\ell}$ is a rank$-1$ factor, which is a contradiction. We conclude that $W = \mathbb{R}^{d}$ so $V$ act minimally.
\end{proof}
We will need a lemma like Lemma \ref{L:Background2}, but only considering coarse directions that lie in the same stable direction for some element of the action. The following lemma is a consequence of Lemma \ref{L:Background1}.
\begin{lemma}\label{L:ImprovedRegularity1}
Let $\alpha_{0}:\mathbb{Z}^{k}\to{\rm Aff}(X_{\Gamma})$ be higher rank. We say that $\mathbf{n}\in\mathbb{Z}^{k}$ is regular if the center of $\alpha_{0}^{\mathbf{n}}$ coincide with the joint center of $\alpha_{0}$. Let $\mathbf{n}_{0}$ be regular and $E_{0}^{s}$ be the stable space associated to $\alpha_{0}^{\mathbf{n}_{0}}$. Let $[\chi]$ be a coarse exponent such that $\ker[\chi]$ is a wall for the Weyl chamber that contains $\mathbf{n}_{0}$ and $E_{0}^{[\chi]}\subset E_{0}^{s}$. Either $E_{0}^{s} = E_{0}^{[\chi]}$ or the complementary subspace $E_{0}^{ss}$, $E_{0}^{ss}\oplus E_{0}^{[\chi]} = E_{0}^{s}$, defines a minimal foliation in $X_{\Gamma}$.
\end{lemma}
\begin{proof}
After projecting to the base, we assume without loss of generality that $X_{\Gamma}$ is a torus $\mathbb{T}^{d}$. If $E_{0}^{[\chi]} = E_{0}^{s}$ for $[\chi](\mathbf{n}_{0}) < 0$ then there is nothing to prove, so assume that there is at least one coarse exponent $[\eta]$ satisfying $[\eta]\neq[\chi]$ and $[\eta](\mathbf{n}_{0}) < 0$. Write $E_{0}^{ss}$ for the complementary subspace defined by
\begin{align}
E_{0}^{ss} = \bigoplus_{\substack{[\eta](\mathbf{n}_{0}) < 0 \\ [\eta]\neq[\chi]}}E_{0}^{[\eta]}.
\end{align}
Let $W$ be the rational closure of $E_{0}^{ss}$, then $W$ is $\alpha_{0}-$invariant and rational. Each $\alpha_{0}^{\mathbf{n}}|_{W}$ preserves the lattice $W\cap\mathbb{Z}^{d}$ so $\det(\alpha_{0}^{\mathbf{n}}|_{W}) = \pm1$. Given a coarse exponent $[\eta]$, define $\eta' = r(\eta')\eta$ for $[\eta'] = [\eta]$ and $d_{[\eta]}^{W}\in[0,1]$ by
\begin{align}
d_{[\eta]}^{W} = \frac{\sum_{\eta'\in[\eta]}r(\eta')\dim(E_{0}^{\eta'}\cap W)}{\sum_{\eta'\in[\eta]}r(\eta')\dim(E_{0}^{\eta'})}.
\end{align}
From this definition it is immediate that
\begin{align}
\left|\det(\alpha_{0}^{\mathbf{n}}|_{E_{0}^{[\eta]}\cap W})\right| = \left|\det(\alpha_{0}^{\mathbf{n}}|_{E_{0}^{[\eta]}})\right|^{d_{[\eta]}^{W}},\quad\mathbf{n}\in\mathbb{Z}^{k}.
\end{align}
Our choice of $W$ implies that $d_{[\eta]}^{W} = 1$ if $[\eta](\mathbf{n}_{0}) < 0$ and $[\eta]\neq[\chi]$. Rewrite
\begin{align*}
1 = & \left|\det(\alpha_{0}^{\mathbf{n}}|_{W})\right| =  \prod_{[\eta]}\left|\det(\alpha_{0}^{\mathbf{n}}|_{E_{0}^{[\eta]}})\right|^{d_{[\eta]}^{W}} = \\ &
\left(\prod_{\substack{[\eta](\mathbf{n}_{0}) < 0 \\ [\eta]\neq[\chi]}}|\det(\alpha_{0}^{\mathbf{n}}|_{E_{0}^{[\eta]}})|\right)\cdot\left|\det(\alpha_{0}^{\mathbf{n}}|_{E_{0}^{[\chi]}})\right|^{d_{[\chi]}^{W}}\cdot\left(\prod_{[\eta](\mathbf{n}_{0}) > 0}|\det(\alpha_{0}^{\mathbf{n}}|_{E_{0}^{[\eta]}})|^{d_{[\eta]}^{W}}\right).
\end{align*}
Using $|\det(\alpha_{0}^{\mathbf{n}})| = 1$
\begin{align*}
\prod_{\substack{[\eta](\mathbf{n}_{0}) < 0 \\ [\eta]\neq[\chi]}}|\det(\alpha_{0}^{\mathbf{n}}|_{E_{0}^{[\eta]}})| = \frac{1}{\left|\det(\alpha_{0}^{\mathbf{n}}|_{E_{0}^{[\chi]}})\right|}\cdot\prod_{[\eta](\mathbf{n}_{0}) > 0}\frac{1}{\left|\det(\alpha_{0}^{\mathbf{n}}|_{E_{0}^{[\eta]}})\right|}
\end{align*}
and combining estimates
\begin{align}
1 = \left|\det(\alpha_{0}^{\mathbf{n}}|_{E_{0}^{[\chi]}})\right|^{d_{[\chi]}^{W} - 1}\cdot\prod_{[\eta](\mathbf{n}_{0}) > 0}|\det(\alpha_{0}^{\mathbf{n}}|_{E_{0}^{[\eta]}})|^{d_{[\eta]}^{W}-1}.
\end{align}
If $\mathbf{n}_{j}\in\mathbb{Z}^{k}$ is a sequence in the same Weyl chamber as $\mathbf{n}_{0}$ such that $\mathbf{n}_{j}\to\ker[\chi]$, then
\begin{align*}
\left|\det(\alpha_{0}^{\mathbf{n}_{j}}|_{E_{0}^{[\chi]}})\right|^{d_{[\chi]}^{W} - 1}\to 1.
\end{align*}
Each $[\eta]\neq\pm[\chi]$ satisfies $|\det(\alpha_{0}^{\mathbf{n}_{j}}|_{E_{0}^{[\eta]}})|\geq\mu > 1$ for some uniform $\mu$. Letting $j\to\infty$, $d_{[\eta]}^{W} = 1$ for $[\eta](\mathbf{n}_{0}) > 0$, $[\eta]\neq-[\chi]$. It follows that $E_{0}^{[\eta]}\subset W$ for $[\eta]\neq\pm[\chi]$. If $\mathbb{T}^{d}/W$ is non-trivial then the projected action on $\mathbb{T}^{d}/W$ does not have two independent coarse Lyapunov exponents. By Lemma \ref{L:Background1} $\mathbb{T}^{d}/W$ is a rank$-1$ factor of $\alpha_{0}$. This contradicts the assumption that $\alpha_{0}$ is higher rank so $W = \mathbb{R}^{d}$ which proves the lemma.
\end{proof}
\begin{lemma}\label{L:RankVersusNumberOfLyapunovExp}
Let $A\subset{\rm GL}(n,\mathbb{Z})$ be a free abelian subgroup with Lyapunov exponents ${\rm Lyap}(A)$. Let $N$ be the maximal number of linearly independent Lyapunov exponents (so $N = \dim({\rm span}({\rm Lyap}(A)))$). If the intersection of the kernels of all $\chi\in{\rm Lyap}(A)$ is trivial in $A$ then ${\rm rank}(A)\leq N$.
\end{lemma}
\begin{proof}
Let $p(t)\in\mathbb{Z}[t]$ be of degree $d$, monic and with constant term $\pm1$. Let $\lambda_{1},...,\lambda_{d}$ be the roots of $p(t)$ (possibly with multiplicity). There is a constant $\mu_{d} > 1$ such that either $p(t)$ has only roots on the unit circle or
\begin{align}
M(p(t)) = \prod_{j = 1}^{d}\max(1,|\lambda_{j}|)\geq\mu_{d}
\end{align}
see for example \cite{Dobrowolski1979}. We number $\{\chi_{1},...,\chi_{n}\} = {\rm Lyap}(A)$. Let $d_{j}$ be the dimension of the Lyapunov space associated to $\chi_{j}$. Given $a\in A$ let $p_{a}(t)$ be the corresponding characteristic polynomial. We obtain
\begin{align}
\log M(p_{a}(t)) = \sum_{\chi_{j}(a) > 0}d_{j}\chi_{j}(a).
\end{align}
Suppose that ${\rm rank}(A) > N$ for contradiction. Let $\chi_{1},...,\chi_{N}$ be chosen such that every $\chi_{j}$ lie in ${\rm span}(\chi_{1},...,\chi_{N})$. If $a_{n}\in A$ is such that $\chi_{1}(a_{n}),...,\chi_{N}(a_{n})\to0$ then $\chi_{j}(a_{n})\to0$ for all $j = 1,...,n$. The intersection
\begin{align}
V := \bigcap_{j = 1}^{N}\ker\chi_{j}\subset A\otimes\mathbb{R}
\end{align}
is non-trivial since ${\rm rank}(A) > N$. The set $A$ is a lattice in $A\otimes\mathbb{R}$ so we find a sequence $a_{n}\in A$ such that $\intd(a_{n},V)\to 0$ as $n\to\infty$ but $a_{n}\not\to e$ (where $e$ is the identity in $A\subset{\rm GL}(n,\mathbb{Z})$). It follows that
\begin{align}
\log M(p_{a_{n}}(t)) = \sum_{\chi_{j}(a) > 0}d_{j}\chi_{j}(a_{n})
\end{align}
tends to $0$ as $n\to\infty$. With $n$ big enough $\log M(p_{a_{n}}(t)) < \log\mu_{d}$ which implies that $p_{a_{n}}$ has all roots on the unit circle. It follows that $a_{n}$ lie in the kernel of all $\chi_{j}$, so $a_{n} = e$ by assumption. This implies that $a_{n}\to e$ which is a contradiction.
\end{proof}
\begin{lemma}\label{L:RankOfCentralizerInSymplecticGroup}
Let $A\in{\rm Sp}(d,\mathbb{Z})$ be hyperbolic with irreducible characteristic polynomial. If $r_{1}(A)$ denotes the number of real eigenvalues of $A$ and $r_{2}(A)$ the number of pairs of complex eigenvalues of $A$ then ${\rm rank}(Z_{{\rm Sp}(d,\mathbb{Z})}(A)) = r_{1}(A)/2 + r_{2}(A)/2$.
\end{lemma}
\begin{proof}
Fix $A\in{\rm Sp}(d,\mathbb{Z})$ with irreducible characteristic polynomial. Let 
\begin{align}
\Lambda_{2}(\mathbb{R}^{d}) = \mathbb{R}^{d}\wedge\mathbb{R}^{d}
\end{align}
be the vector space of $2-$vectors. Let $A\wedge A:\Lambda_{2}(\mathbb{R}^{d})\to\Lambda_{2}(\mathbb{R}^{d})$ be the induced map on $2-$vectors. Note that $\Lambda_{2}(\mathbb{Z}^{d})$ is a $A\wedge A-$invariant lattice in $\Lambda_{2}(\mathbb{R}^{d})$. Denote by $W\leq\Lambda_{2}(\mathbb{R}^{d})$ the eigenspace of $1$ for $A\wedge A$. Write $\Gamma := W\cap\Lambda_{2}(\mathbb{Z}
^{d})$, since $W$ is a rational subspace of $\Lambda_{2}(\mathbb{R}^{d})$ the subgroup $\Gamma\leq W$ is a lattice. 

Given $B\in Z_{{\rm GL}(d,\mathbb{Z})}(A)$ the wedge $B\wedge B$ preserves $W$ and stabilize $\Gamma$ in $W$. So, after identifying $W\cong\mathbb{R}^{d/2}$ and $\Gamma\cong\mathbb{Z}^{d/2}$ we obtain a map $\Psi:Z_{{\rm GL}(d,\mathbb{Z})}(A)\to{\rm GL}(d/2,\mathbb{Z})$ defined by $\Psi(B) := (B\wedge B)|_{W}$. It is immediate that $\Psi$ is a homomorphism. Fix eigenvectors $e_{1},...,e_{d/2},\Tilde{e}_{1},...,\Tilde{e}_{d/2}\in\mathbb{C}^{d}$ such that $Ae_{j} = \lambda_{j}e_{j}$ and $A\Tilde{e}_{j} = \Tilde{e}_{j}/\lambda_{j}$. We identify (the complexification of) $W$ by $W = {\rm span}(e_{1}\wedge\Tilde{e}_{1},...,e_{d/2}\wedge\Tilde{e}_{d/2})$ so for $B\in Z_{{\rm GL}(d,\mathbb{Z})}(A)$ we have
\begin{align}\label{Eq:LyapunovOfMorphism}
B\wedge B(e_{j}\wedge\Tilde{e}_{j}) = \mu_{j}(B)\Tilde{\mu}_{j}(B)e_{j}\wedge\Tilde{e}_{j}
\end{align}
where $Be_{j} = \mu_{j}(B)e_{j}$ and $B\Tilde{e}_{j} = \Tilde{\mu}_{j}(B)\Tilde{e}_{j}$. If $\Psi(B) = e$ then $\mu_{j}(B)\Tilde{\mu}_{j}(B) = 1$ which implies that $B$ preserve the symplectic form that $A$ preserve (note that $e_{1},...,e_{d/2},\Tilde{e}_{1},...,\Tilde{e}_{d/2}$ can be chosen such that the symplectic form can be written $e^{1}\wedge\Tilde{e}^{1}+...+e^{d/2}\wedge\Tilde{e}^{d/2}$). Conversely, if $B\in Z_{{\rm Sp}(d,\mathbb{Z})}(A)$ then $\Psi(B) = e$. Equation \ref{Eq:LyapunovOfMorphism} also implies that the Lyapunov exponents of ${\rm Im}(\Psi)$ are given by
\begin{align}
\chi_{j}^{\Psi}(\Psi(B)) = \log|\mu_{j}(B)| + \log|\Tilde{\mu}_{j}(B)|.
\end{align}
It follows that ${\rm Im}(\Psi)$ has $r_{1}(A)/2 + r_{2}(A)/2$ Lyapunov exponents. Indeed, if $\mu_{j}(B),\Tilde{\mu}_{j}(B)$ takes values in $\mathbb{R}$ then $\log|\mu_{j}(B)| + \log|\Tilde{\mu}_{j}(B)|$ defines one Lyapunov exponent. If $\mu_{j}(B),\Tilde{\mu}_{j}(B)$ takes values in $\mathbb{C}\setminus\mathbb{R}$ then $\overline{\mu_{j}(B)} = \mu_{j'}(B)$ and $\overline{\Tilde{\mu}_{j}(B)} = \Tilde{\mu}_{j'}(B)$ for some $j\neq j'$, so
\begin{align}
\chi_{j}^{\Psi}(B) = \chi_{j'}^{\Psi}(B).
\end{align}
By Lemma \ref{L:RankVersusNumberOfLyapunovExp} ${\rm rank}({\rm Im}(\Psi))\leq r_{1}(A)/2 + r_{2}(A)/2 - 1$ since $|\det(\Psi(B))| = 1$ implies that, at least, one Lyapunov exponent of ${\rm Im}(\Psi)$ can be written as a combination of the other exponents. By \cite{KatokKatokSchmidt2002} the rank of $Z_{{\rm GL}(d,\mathbb{Z})}(A)$ is $r_{1}(A) + r_{2}(A) - 1$, so
\begin{align*}
{\rm rank}(\ker(\Psi)) = & {\rm rank}(Z_{{\rm GL}(d,\mathbb{Z})}(A)) - {\rm rank}({\rm Im}(\Psi))\geq \\ &
r_{1}(A) + r_{2}(A) - 1 - \left[\frac{r_{1}(A) + r_{2}(A)}{2} - 1\right] = \\ &
\frac{r_{1}(A) + r_{2}(A)}{2}
\end{align*}
or since $\ker(\Psi) = Z_{{\rm Sp}(d,\mathbb{Z})}(A)$, ${\rm rank}(Z_{{\rm Sp}(d,\mathbb{Z})}(A))\geq(r_{1}(A) + r_{2}(A))/2$. The converse inequality is clear from Lemma \ref{L:RankVersusNumberOfLyapunovExp} since if $\chi$ is a Lyapunov exponent of $Z_{\rm Sp(d,\mathbb{Z})}(A)$ then so is $-\chi$.
\end{proof}

%\input{FrenchTitleAndName}

%------------------------------------------------------------%

%\input{OldLemmas}

%\input{TEST}

%\input{ShortIntroduction}

%\input{InvariancePrincipleHigherRank}

\newpage
\bibliography{main.bib}{}
\bibliographystyle{abbrv}

\end{document}